\newif\ifSiam
\newif\ifArxiv

\Arxivtrue

\ifSiam
\documentclass{siamart220329}
\usepackage{amsmath}
\usepackage{amssymb}
\usepackage{enumitem}
\newtheorem{example}{Example}[section]
\newtheorem{remark}{Remark}

\fi

\ifArxiv
\documentclass{article}
\usepackage{geometry}
 \usepackage{authblk}
 \usepackage{bbding}

\usepackage{amsmath}
\usepackage{amsthm}
\numberwithin{equation}{section}
\usepackage{amssymb}
\newtheorem{theorem}{Theorem}[section]
\newtheorem{example}{Example}[section]
\newtheorem{corollary}[theorem]{Corollary}

\newtheorem{proposition}[theorem]{Proposition}
\newtheorem{remark}{Remark}
\usepackage[colorlinks=true]{hyperref}

\hypersetup{
    linkcolor = {teal},
    citecolor={blue}
    }
\fi

\usepackage{subcaption}
\usepackage{ upgreek }
\usepackage{dsfont}
\usepackage{graphicx} 
\usepackage{xcolor}
\usepackage{soul}
\usepackage{comment}
\usepackage{enumitem}

\newcommand{\R}{\mathbb{R}}
\newcommand{\N}{\mathbb{N}}

\newcommand{\XJ}{J}
\newcommand{\eA}{\mathbf{e}_A}

\newcommand{\KA}{K}
\newcommand{\hE}{\hat{E}}
\newcommand{\bE}{\bar{E}}

\newcommand{\hXJ}{\check{J}}
\newcommand{\hXA}{\check{A}}
\newcommand{\bXJ}{\hat{J}}
\newcommand{\bXA}{\hat{A}}
\newcommand{\scXJ}{\bar{J}}
\newcommand{\scXA}{\bar{A}}
\newcommand{\tXJ}{\tilde{J}}
\newcommand{\tXA}{\tilde{A}}
\newcommand{\tbeta}{\tilde{\beta}}

\newcommand{\CP}{\mathbf{C}_p}

\newcommand{\XA}{A}
\newcommand{\sA}{\sigma_A}
\newcommand{\phiA}{\phi_A}
\newcommand{\phiJ}{\phi_J}
\newcommand{\sJ}{\sigma_J}
\newcommand{\sE}{\sigma_E}
\newcommand{\sP}{\sigma_P}
\newcommand{\ZZ}{{\bf X}}
\newcommand{\K}{{\bf K}}

\newcommand{\tZZ}{{\bf \tilde{Z}}}

\newcommand{\ind}[1]{{\bf 1}_{#1}}
\newcommand{\cuteps}{\Psi_{\varepsilon}}
\newcommand{\cut}{\Psi}
\newcommand{\cutDt}{\Psi_{\Delta t}}
\newcommand{\FJ}{F^{(J)}}
\newcommand{\FA}{F^{(A)}}

\newcommand{\kappaA}{\kappa_{A}}

\newcommand{\kappaJ}{\kappa_{J}}
\newcommand{\deltaWA}[2]{\updelta \WA_{#1,#2}}
\newcommand{\deltaWJ}[2]{\updelta \WJ_{#1,#2}}
\newcommand{\deltaWE}[2]{\updelta \WE_{#1,#2}}
\newcommand{\kA}{K}
\newcommand{\rhoA}{\rho_A}
\newcommand{\rhoJ}{\rho_J}
\newcommand{\mJ}{m^{(J)}}
\newcommand{\mA}{m^{(A)}}
\newcommand{\rJ}{r^{(J)}}
\newcommand{\rrJ}{\mathbf{\rho}_{J}}
\newcommand{\WA}{W^{(A)}}
\newcommand{\WJ}{W^{(J)}}
\newcommand{\WE}{W^{(E)}}
\newcommand{\WP}{W^{(P)}}
\newcommand{\Fmaxj}{F_{\max}^{(J)}}
\newcommand{\FmaxA}{F_{\max}^{(A)}}
\newcommand{\PjA}{\rho_{A}}
\newcommand{\XX}{{\bf X}}
\newcommand{\bXX}{{\bf \bar{X}}}
\newcommand{\bZZ}{\bar{\ZZ}}
\newcommand{\zz}{{\bf z}}

\newcommand{\E}{{\mathbb E}}
\newcommand{\F}{\mathcal{F}}
\newcommand{\tauU}{\tau_U}

\renewcommand{\P}{{\bf P}}
\newcommand{\dt}{\Delta t}
\newcommand{\ddt}{{\updelta}(t)}
\newcommand{\DXA}{\Delta \XA^\delta}
\newcommand{\DXJ}{\Delta \XJ^\delta}
\newcommand{\DE}{\Delta E^\delta}
\newcommand{\mc}[1]{\mathcal{#1}}
\newcommand{\err}{\mathcal{E}}

\newcommand{\yy}{\mathbf{y}}
\newcommand{\Np}{\mathcal{N}}
\newcommand{\Cmult}{C}
\newcommand{\Cadd}{C}

\usepackage{mdframed}
  {\begin{mdframed}[backgroundcolor=yellow]}%
  {\end{mdframed}}

\title{\bf Age-structured stochastic populations under dynamic  harvesters' behavior: well-posedness, asymptotic stability and numerically-amenable approximations}\usepackage[perpage, symbol*]{footmisc}

\author{M. Isidora \'{A}vila-Thieme \thanks{Center for Resilience, Adaptation and Mitigation (CReAM), Universidad Mayor, Chile}  \thanks{Escuela de negocios, Facultad de Ciencias Sociales y Artes, Universidad Mayor, Chile} \thanks{Center of Applied Ecology and Sustainability (CAPES), Chile}\thanks{Instituto Milenio en Socio-Ecolog\'ia Costera (SECOS), Chile}, Kerlyns Mart\'{i}nez\thanks{Departamento de Ingeniería Matemática, Universidad de Concepción, Chile}, H\'{e}ctor Olivero\thanks{CIMFAV - Instituto de Ingenier\'ia Matem\'atica, Universidad de Valpara\'iso, Chile}, Mauricio Tejo\thanks{CIMFAV - Instituto de Estad\'istica, Universidad de Valpara\'iso, Chile}, Leonardo Videla\thanks{Departamento de Matem\'atica y Ciencia de la Computaci\'on, Universidad de Santiago de Chile, Chile}}

\begin{document}
\maketitle


\begin{abstract}
In this paper we study a model of age-structured ecological populations in continuous interaction with a community of harvesters. We propose an individual-based model for this feedback interactions and prove its convergence to a system of coupled stochastic differential equations. This limit process is complemented with socio-ecological and environmental factors.

For the resulting system we establish well-posedness, positivity, and the non-attainability of extinction boundaries, together with its long-term behavior and numerical approximation, relying on ad-hoc arguments to handle the super-linear growth and non-Lipschitz coefficients. Finally, we illustrate some of our findings with numerical experiments. 

\end{abstract}


\section{Introduction}\label{sec:introduction}
In natural ecosystems, illegal, unreported and unregulated fishing, along with non-compliance with fisheries regulations, are central drivers of global fisheries collapse and a persistent stumbling block for managers and fishers \cite{arias2015understanding,boonstra2014chain, gomez2024moving,  sumaila2019carding}. They generate severe socio-ecological impacts, undermine food security, and cause substantial economic losses that destabilize the social systems relying on fisheries for livelihoods and governance \cite{agnew2009estimating, nahuelhual2020social}. The Chilean kelp fishery, one of the largest wild seaweed fisheries in the world, plays a key socioeconomic role for coastal communities, while also sustaining biodiversity and ecosystem health \cite{ fellous2017modulatory,  kefi2015network, navarrete2023monitoring, navarrete1990resource,vasquez2016brown, vasquez2024brown}. However, weak enforcement, as perceived by fishers, and volatile market prices contribute to over-exploitation and non-compliance. This shift causes populations to move from mature forests to stands dominated by juveniles, or even to severely depleted states \cite{avila2025exploring, gonzalez2021exploring, gonzalez2023artisanal,vasquez2024brown}. These impacts are also influenced by low-frequency environmental perturbations, such as El Niño-Southern Oscillation (ENSO) events, whose warm phase (El Niño) can dramatically reduce kelp abundance and even trigger local extinctions, while its cold phase (La Niña) can enhance the kelp population growth \cite{camus1994recruitment, martinez1999latitudinal}. In this context, non-compliance represents a wicked problem within complex socio-ecological systems, where enforcement or punishment alone is often insufficient \cite{avila2025exploring, nahuelhual2023reframing}. Deterrent strategies, together with fishers’ participation in sustainability programs that offer price premiums for compliance, can help promote responsible fishing behavior \cite{ jones2023enrollment, potts2007international}. However, integrative approaches capable of capturing cross-scale feedbacks are needed to achieve efficient and sustainable fisheries management. One way to support such approaches is through mathematical modeling of harvesting in age-structured populations.

Mathematical modeling of harvesting in age-structured populations has been assessed from different approaches in previous literature. From a deterministic point of view, ordinary differential equations, integro-partial differential equations and partial differential equations of McKendrick-type have been proposed, in which a population depending on the continuously-varying age is posed, and where  the aim is to establish an optimal harvesting policy (see, e.g., \cite{belyakov2016optimal,gurtin1981optimal,  he2019optimal,murphy1990optimal}). In the stochastic setup, and in a similar way to us, dynamics of  population classes have been modeled by including random fluctuations or environmental noises  (see, e.g., \cite{ chowdhury2001stochastic, engen2018spatial, kuusela2020modeling,ludwig1995theory, wikstrom2012role}). Environmental noises are considered as perturbations in fitness, generally set as an additive Brownian motion in per capita growth rates, so that the diffusion coefficient presents a linear dependence. Less frequent but more severe events, such as catastrophes or ENSO events, can also be incorporated through jump processes (typically Poissonian) with variable magnitude \cite{ludwig1995theory}. However, much of this literature focuses on establishing an optimal harvesting rate that maintains the population above a positive stable level, in order to ensure persistence of the overall population, but ignoring the explicit dynamics of the harvesters, whose level of extraction will depend on different socio-ecological factors. 

To the best of our knowledge, the explicit inclusion of the harvester behavior, with its regulations, within a mathematical model for age-structured populations is an issue that has not been widely addressed. A recent attempt made in \cite{wszola2023integrating}, although from a deterministic setting,  considers the dynamics of a population divided into adults and juveniles, where the harvesting for each subpopulation depend upon an intrinsic per-capita capture rates combined with probabilities of keeping legal/illegal captures. Such rates and probabilities are quantified depending on different regulatory scenarios. Nevertheless, these compliance/non-compliance ratios are not dynamic.

In the present paper, we model harvesting dynamics while considering their interaction with an age-structured ecological population as a social-ecological phenomenon. We propose a coupled system describing the feedback between the ecological population and a population of interacting harvesters represented by the vertices of a network. Harvesters may comply or not with regulatory policies, and their behavior is influenced by both market prices and the ecological state. Under a suitable scaling, the model is ultimately represented by a system of stochastic differential equations (SDEs, for short) that incorporate social and ecological factors, including price fluctuations and environmental variability.

Our system is driven by multiple sources of stochasticity: Lévy-type environmental perturbations affecting adult and juvenile kelp, random fluctuations in harvester compliance, and irregular exogenous price behavior. The analysis of harvesting effects therefore corresponds to studying the dynamics of harvester behavior, particularly the conditions under which compliance levels allow the ecological population to persist or go extinct. All of this sets a complex and challenging mathematical analysis. In particular, the influence of noises on the ecological populations and the harverster's behavior are of different nature: The diffusion coefficient is linear in the ecological population and of square-root type for the compliance level. Far from arbitrariness, the difference between these effects will be explained through rigorous modeling choices issued from microscopic (or more precisely, individual-based) considerations relying upon mathematical and social-ecological-theoretic aspects. 

The purpose of this paper is to study the theoretical properties of this model, which introduces several interesting and challenging mathematical problems. To address these challenges, we adopt a two-pronged approach. First, we establish the convergence of an individual-based model for the harvester population to a mean-field stochastic differential equation model, which provides a solid micro-founded basis for our macroscopic description. Second, we employ a variety of analytical techniques to study the resulting system of SDEs. This includes overcoming several mathematical difficulties posed by the super-linear growth of the coefficients and their non-Lipschitzian nature, which is a key feature of the socio-ecological interactions. Importantly, we explore the long-term behavior of the system, including its asymptotic stability, to understand the conditions under which the kelp population can persist over time. Finally, we propose a suitable exponential-truncated numerical scheme for the model, which provides an approximation that is admissible in the state space and strongly convergent. This reliable numerical approximation allows us to explore how harvester's behavior, influenced by incentives such as price premiums for compliance, interacts with ecological dynamics, thereby supporting the design feedback mechanisms for efficient and sustainable fisheries management.

\section{Model description: from abstract individual's behavior to a system of parameterized SDEs} \label{sec:model_description}

In this section, we introduce the modeling assumptions that lead to the system of SDEs. We consider a large group of harvesters interacting with age-structured resources, where each harvester may either comply with regulations or not. Initially, the compliance state of each harvester should depend on the state of other harvesters and the current amount of available resources. In turn, the resources are affected by their natural dynamics (birth and death rates with density-dependent behavior) and by the number of harvesters that do not comply with regulations. Consequently, the ecological population is affected by extraction rates depending on the state of the set of harvesters. In order to provide a mathematical framework for this structure, we proceed as follows. 

Assume that for each $n\in \N$ a certain (possibly random) connected graph $G_n=(V_n, F_n)$ with vertices in $V_n$ and edges in $F_n\subseteq V_n\times V_n$ is given. We assume that $(G_n)_{n\ge 1}$ is a growing sequence of graphs on $n$ vertices, in the sense that $\vert F_n\vert  \to \infty$. We regard the vertices of $G_n$ as the individual fishers having one of two possible opinions: to comply ($1$) or not to comply ($0$) with the regulations. The edges in $F_n$ represent the contact relationships between individuals: for $i,j\in V_n$, if $\{i,j\}\in F_n$, then the opinion of $i$ regarding compliance can be directly influenced by that of $j$ (and the other way around). 

Now, the decisions of the fishers (i.e. the state of the compliance graph described above) should exhibit a dynamic behavior in continuous feedback with the ecological population. So, we prescribe a coupled dynamics for both, the ecological population, which has to stages (juveniles and adults), and the fisher's states through a $\mathbb{Z}^2_{+} \times \{0,1\}^{V_n}$-valued process $(\K^{(n)}_t, \xi^{(n)}_t: t \ge 0)$.  The $\mathbb{Z}^{2}_{+}$-valued component $\K^{(n)}= (J^{(n)}, A^{(n)})$ represents the number of units of biomass of the juvenile and adult population of the resource, meanwhile for $i \in V_n$, $\xi_i^{(n)}(t) \in \{0,1\}$, represents the compliance behavior of individual $i$ at time $t$ in the graph $G_n$. 

Let $E^{(n)}_t:=\sum_{i}\xi^{(n)}_i(t)$ be the counting process associated to the compliance process, i.e. the number of individuals that comply with regulations at time $t$. At the ticks of independent and exponentially distributed random clocks with rate $\gamma > 0$, the individual $i$ changes its opinion to that of one of the fishers chosen uniformly at random from among its neighbors on $G_n$. This change in the opinion can be regarded as the neutral part of the dynamics: the agents are more prone to comply if many of its direct acquaintances are already complying. Superimposed over this neutral change of opinion, the individual can reconsider their behavior further, but this time the mutation is based upon the information individual $i$ can obtain from the dynamics of the population process $\K^{(n)}= (J^{(n)}, A^{(n)})$ and from the assessment they have of the overall opinions of the harvesters $\xi^{(n)}$. More precisely, assume that individual $i$ was chosen to resample his/her opinion at the random time $T$. If, after looking at his/her neighbors (after the neutral step), the new opinion is $1$, with probability $\dfrac{1}{n}\beta_0 (\K^{(n)}_{T-}/n,E^{(n)}_{T-}/n)$, the individual decides not to comply. On the contrary, if its current state is $0$, then with probability $\dfrac{1}{n}\beta_1 (\K^{(n)}_{T-}/n,E^{(n)}_{T-}/n)$, he/she decides to comply. 

Coupled with the above dynamics, $\K^{(n)}$  evolves by jumps of size $1$ over each component of the lattice. If $\K^{(n)}_t = (j, a)$ we say that there are $j$ juveniles individuals and $a$ adult individuals by time $t$. This component evolves as follows. Independently of every other individual, and also independently from the compliance process, each juvenile individual becomes an adult at rate $\rho_A/n$. Conditionally on $\xi^{(n)}_t$ and $\K^{(n)}_t$, each adult individual, independently of the other individuals, gives birth to a juvenile on the time interval $[t, t+dt]$ with probability $dt \dfrac{1}{n}\rho_J (\K^{(n)}_t/n, E^{(n)}_t/n) + o(dt)$; conditionally on $\xi^{(n)}_{t}$, every adult (resp. juvenile) dies, independently of all the other individuals, at rate $\dfrac{1}{n}\kappa_A (E^{(n)}_t/n)$ (resp. $\dfrac{1}{n}\kappa_J (E^{(n)}_t/n)$).

The connection of the voter model with certain Wright-Fisher type SDEs (see for example the seminal \cite{chenchoicox2016}), suggests that under various hypotheses on the sequence of graphs $(G_n)_{n \ge 1}$, it is reasonable to expect that as $n \to \infty$, the sequence of time-scaled processes $((\dfrac{1}{n}\K^{(n)}_{nt}, \dfrac{1}{n}E^{(n)}_{nt }): t \ge 0)_{n \ge 2}$ (which are not Markov in general) converges (in law with probability $1$ with respect to the randomness from which the sequence of graphs was sampled) to a Markov process $\ZZ:=(J, A, E)$ satisfying the system of equations:

\begin{equation}\label{eq:simplified-model}
\begin{aligned}
    \XJ_t
    =\XJ_0&+\int_0^t \rrJ(\XJ_s,\XA_s,E_s)\XA_s  -  \kappaJ(E_s)\XJ_sds  \\
    \XA_t
    = \XA_0 &+\int_0^t  \PjA \XJ_sds - \int_0^t \kappaA(E_s)\XA_sds \\
    E_t =E_0 &+ \int_0^t\tilde \beta_1(\ZZ_s)(1-E_s)  -\tilde \beta_0(\ZZ_s)E_sds +\int_0^t\tilde \sE\sqrt{E_s(1-E_s)}dW^{(E)}_s,
\end{aligned} 
\end{equation}
for certain coefficients $\tilde \beta_0, \tilde \beta_1, \tilde \sE$ which can be explicitly determined, and which depend on global features of the social structure given by the sequence of graphs and the rate of resampling $\gamma$.

In Section \ref{sec:derivation}, we demonstrate that under mild hypotheses, when $(G_n)_{n \ge 1}$ represents the sequence of complete graphs, the claimed convergence in law holds.

Several ecological and socio-ecological considerations naturally lead us to explore further enrichments of the dynamics represented by \eqref{eq:simplified-model}. These enrichments follow two main directions.
\begin{enumerate}
\item The decisions of the harvesters can be affected by the trading prices whereby the fishers are more prone to extract as the price gets higher. 

\item The instantaneous rate of change of the populations can be perturbed by different sources of noise, accounting for environmental stochasticity, as opposed to \emph{demographic noise} (which is usually observed with variances of the order of the population).
\end{enumerate}

Regarding the first point, we will assume that kelp's price will act as an exogenous process (i.e., this process emulates a ``market price'' which is not affected by our local dynamics of kelp population and the compliance of the harvesters). 
Irrespective of any particular choice, we assume hereafter that we are provided with a right-continuous, $\R_+$-valued process $(P_t: t \ge 0)$ which accounts for the observed price of the kelp. Now, it is expected that if the price is high, the fishers will be prone to harvest with scarce care for the regulations. It is natural to expect that fisher behavior changes when prices are very high or very low. Such effects should be incorporated in the parametrization of the functions $\beta_1$ and $\beta_0$ in \eqref{eq:simplified-model} so that the rate of change of harvesters' opinions at time $t$ depends explicitly on the current price $P_t$.    

As for the second point, the introduction of random perturbations in fitness can be explained by the necessity to account for external random effects -whether biotic or abiotic- that may influence the resources' ability to reproduce. 

Naturally, thermal or luminic variations, which are perceived on the short-time scale, can be introduced through multiplicative Gaussian noise. 
To introduce perturbations  whose natural time-scale is much slower, it is natural to regard these perturbations as pulses affecting the average fitness of the populations under consideration. Indeed, such perturbations are ubiquitous in marine ecosystems. As examples we can cite: sudden migration of resident species or incoming flow of biomass through ecosystem boundaries, etc. In fact, multiple other stressors associated with climate (such as ENSO events) can compromise the capacity of these socio-ecological systems to respond to perturbations.
In this context, upwelling or cold ENSO phases may dampen the negative impacts of non-compliance on algal populations, while the warm ENSO phase may intensify them. In our model, these perturbations (as well as others: intense coastal storm, sudden swells, etc.) are assumed to be introduced through some (well-behaved) jump Markov processes. This explain the stochastic integral with respect to a Poisson random measure in \eqref{eq:full-model}.

In this way, we arrive to the target model, which this article is mostly devoted to. In the sequel we consider a fixed filtered probability space $(\Omega, \mc{F}, \mathbb{P}, (\F_t: t \ge 0))$, with a filtration $(\F_t: t \ge 0)$ obeying the usual assumptions of right-continuity and completeness. Our target model is made up of four components denoted by $(\XJ,\XA, E, P)$,  where ($\XJ,\XA$) represents the biological dynamics of kelp population (in terms of biomass) separated in juvenile and adults, $E$ represent the fishers behavior around compliance, and $P$ the price of kelp. As explained above, this last component is completely exogenous, and its dynamic does not impact in the results derived below as long as some very mild assumptions are fulfilled. Putting $\K=(\XJ,\XA)$ and $\ZZ=(\XJ,\XA, E, P)$, the model for the components of $\ZZ$ can be written as follows:
\begin{equation}\label{eq:full-model}\tag{M}
\left\{\begin{aligned}
    \XJ_t
    =\XJ_0&+\int_0^t\rrJ(\K_s,E_s)\XA_s  - \kappaJ(E_s)\XJ_sds + \int_{0}^t\sJ \XJ_s dW^{(J)}_s  + \int_{0}^{t^+} \int_{\R\setminus\{0\}} \phiJ (\XJ_{s-}, z) \Np(dz\otimes ds)\\
    \XA_t= \XA_0 &+\int_0^t \PjA\XJ_sds - \int_0^t\kappaA(E_s)\XA_sds + \int_{0}^t\sA \XA_s dW^{(A)}_s  + \int_{0}^{t^+}\int_{\R\setminus\{0\}} \phiA (\XA_{s-}, z) \Np(dz\otimes ds) \\
    E_t =E_0 &+ \int_0^t\beta_1(\ZZ_s)(1-E_s)  -\beta_0(\ZZ_s)E_sds +\int_0^t\sE\sqrt{E_s(1-E_s)}dW^{(E)}_s\\
    P_t = P_0 &+ \int_0^t \mu(P_s)ds + \int_0^t \sP(P_s)d\WP_s,
\end{aligned} 
\right.
\end{equation}
where $\WJ,\WA,\WE$ and $\WP$ are independent standard Brownian motions defined on $(\Omega, \mc{F}, \mathbb{P}, (\F_t: t \ge 0))$, and independent of $\Np (dz\otimes dt)$, a Poisson random measure with intensity given by $\nu(dz)\cdot\lambda dt$, with  $\nu$ a probability measure on $(\mathbb{R}, \mathcal{B} (\mathbb{R}))$ and $\lambda>0$. Notice that $\pi(dz):=\lambda\nu(dz)$ represents a L\'evy measure. 

We explain now the functional form of the coefficients appearing in \eqref{eq:full-model}. To do so, we introduce first a family of $1$-Lipschitz cut-off function
\begin{equation}\label{eq:def-psi}   
\cuteps(z) = (\varepsilon \lor z) \wedge (1-\varepsilon),
\end{equation}
with the convention $\cut: = \cut_0$.

\paragraph*{The ecological parameters.} The functions appearing in the dynamics of the kelp biomass are given by
\begin{align*}
\rrJ(\K, E)&:= \rJ(E)\Psi\left(1-\frac{\XA+\XJ}\KA\right),\\
\kappaJ(E)&:=\left(\PjA + \Fmaxj(E) + \mJ\right),\\
\kappaA(E)&:=\left(\mA +\FmaxA(E)\right),
\end{align*}
where $\KA$ represents the carrying capacity of adult kelp. At carrying capacity $\KA$, new juveniles can be produced, but cannot grow successfully due to shading, whiplash by adult plants, lack of space and/or nutrients \cite{navarrete2023monitoring}. The constants $\mJ$ and $\mA$ represent the natural mortality rate for juveniles and adults, respectively, and $\rhoA$ denotes the maturity rate, namely the rate at which juveniles become adults. Here, the functions $\kappaJ,\kappaA$, depending on the compliance process $E$, describe the biomass loss and may also depend on the price $P$ through the extraction rates $\Fmaxj,\FmaxA$.

\paragraph*{The parameters of the compliance process $E$.} \label{sec:motivation_process_E}

The process $E$, representing the fraction of fishers that comply with harvesting regulation, is determined by the functions $\beta_0, \beta_1$, which are given by: 
\begin{align}
 \beta_0(\ZZ) 
 &= \left(\Psi\left(\frac{\XJ+\XA}{\KA}\right) +[1- \mathcal{I}^{\eta}_{P_{\min}}(P)+\mathcal{I}^{\eta}_{P_{\max}}(P)] + 1\right)\tauU(1-E) + \bar\beta_{0}\label{eq:def_beta_1}\\
 \beta_1(\ZZ) 
 &= \left( 1- \Psi\left(\frac{\XJ+\XA}{\KA}\right) +\mathcal{I}^{\eta}_{P_{\min}}(P + s)+1 \right)\tauU E + \bar\beta_{1},\label{eq:def_beta_2}
\end{align}
where 
\begin{equation}\label{eq:activation-function}
    \mathcal{I}^{\eta}_{p_0}(p):=\frac{1}{1+\exp(-\eta(p-p_0))},
\end{equation}
is a regularized version of an indicator $p \mapsto \mathds{1}_{[p_0,+\infty)}(p)$. 

As explained above, the functions $\beta_0$ (resp. $\beta_1$) represents the rate of passage from compliance to non-compliance behavior (resp., non-compliance to compliance). Each rate function is determined by the following three decision rules (for further insights on this decision-making-based design, see \cite{avila2025}):

\begin{itemize}

\item The first decision rule depends on population size. As kelp abundance increases, fishers are more likely to shift from compliance to non-compliance, since they perceive their impact would be negligible (first term in \eqref{eq:def_beta_1}). Conversely, when kelp declines, the chances of shifting from non-compliance to compliance increases, as they recognize the resource is depleting and require recovery (first two terms in \eqref{eq:def_beta_2}). Here, the use of $\cut$ intends to restrict the range of values between $0$ and $1$ to avoid mathematical complexities.

\item The second decision rule is associated with the kelp price $P$, which considers two thresholds: a  price perceived by fishers as sufficient to encourage compliance ($P_{\min}$), and a price perceived as high enough to encourage non-compliance ($P_{\max}$). Prices within the range [$P_{\min}$, $P_{\max}$]  promote the change from non-compliance to compliance (second term in \eqref{eq:def_beta_2}), while prices outside this range encourage the transition from compliance to non-compliance. Furthermore, participation in sustainability programs provides an alternative strategy for compliance by offering an incentive ($s$) on top of the kelp price, increasing fishers incomes. If the combined price and incentive exceed the compliance threshold $P_{\min}$, the chances of fishers changing from non-compliance to compliance increases. In all cases, the activation function has a sigmoid shape, whose steepness is defined by the parameter $\eta$. 

\item The third decision rule is associated with the global behavior of fisher's aggregate (as a function of the proportion of compliant harvesters $E$ or non-compliant harvesters $(1-E)$). Fishers are usually organized through a syndicate, with whom they share internal rules. In case of non-compliance with the syndicate's rules, the syndicate could exert some pressure or punishment on its fishers. Consequently, under the pressure of not being reprimanded by their syndicate, the chances that a fisher's behavior will change depend on the behavior of others (third terms in \eqref{eq:def_beta_1}-\eqref{eq:def_beta_2}). Here, the parameter $\tau_U$ is a measure of the fraction of fishers that are associated with a syndicate.

\end{itemize}

\paragraph*{The price factor $P$}
The process $P$ represents the kelp price, which, historically, is highly fluctuating. In some cases, the price is so low that the active extraction might decrease to almost zero. Therefore, in the biological model, $P$ influences $\FmaxA,\Fmaxj$ and $\rJ$ explicitly through an activation equation. However, the kelp price is an exogenous force for fishers. Indeed, most harvested kelp is exported to the Asian market, with Chile being the world's main exporter of wild kelp. The international market dynamics determine both the demand and pricing of kelp, irrespective of local resource availability \cite{porras2020extractivismo}. Therefore, we will assume here that the resource availability does not influence directly the price. 
For our purposes, it is enough to consider the price as the solution of a general model, where, for simplicity, we assume sufficiently smooth drift and diffusion coefficients (see Assumption \ref{H4}).

\subsection*{Aim and organization of this article}
The aim of this paper is to investigate the theoretical properties of the model \eqref{eq:full-model}. These properties are crucial for gaining a deeper understanding of the underlying dynamics, and could assist in the evaluation of socio-political mechanisms that benefit of sustainable economic activities. In this, we not only contribute to the theory of SDEs with non-Lipschitz coefficients but also offer a robust framework for understanding the intricate feedback dynamics between human behavior and ecological systems.

The rest of this paper is organized as follows. In Section \ref{sec:main_results} we present a brief description of the main results that the reader will find thereafter. In Section \ref{sec:derivation}, we first show that the sequence of time-scaled processes $((\dfrac{1}{n}\XX^{(n)}_{nt}, E^{(n)}_{nt }): t \ge 0)_{n \ge 2}$ converges to the proxy system \eqref{eq:simplified-model}. Later in the same section, we study conditions that guarantee the well-posedness of the full system \eqref{eq:full-model}.

Section \ref{sec:inaccessibility_origin_invariant} is devoted to the stability of the system and the identification of condition under which we can guarantee invariant distributions that persist over time. Finally,  the strong convergence analysis of a suitable time-approximation of the solution of the system \eqref{eq:full-model} and simulation study is presented in Section \ref{sec:numerical_scheme}.
Open paths for further research are briefly discussed in Section \ref{sec:closing}. Some technical proofs are gathered in Appendix \ref{sec:appendix}.

\paragraph*{Notations} 
The notation $\R^k_+$ and $\R^k_{++}$ stand for the closed (respectively open) cone of non-negative (resp. strictly positive) real $k$-tuples. For a metric space $S$, the notation $\mathcal{B}(S)$ stands for the Borel $\sigma$-field of $S$.

In what follows, $a \lesssim b$ means $a\leq Cb$ where $C$ is a constant depending only on parameters of the model and the geometry of $\R^4$. We denote by $\texttt{diag}(\mathbf{x})$ the diagonal matrix with diagonal given by the vector $\mathbf{x}$. For a Brownian motion $W$, we write $\updelta W_{s,t}$ to denote its increment $W_t -W_s$.

\section{Overview of main results}\label{sec:main_results}

In this section, we collect the main mathematical aspects of the model \eqref{eq:full-model}. The proofs are deferred to the following sections. 
From now on, let $\mathcal{S}:=\mathbb{R}^2_+ \times [0,1] \times \mathbb{R}_+$ be the state-space of the state-vector $\ZZ = (\XJ,\XA, E, P)$, endowed with the product topology.

In order to prove that the model \eqref{eq:full-model} determine a strong Markov process taking values on $\mathcal{S}$,
we consider the following set of hypotheses.

\paragraph{Assumptions}

\begin{enumerate}[start=1,label={(\bfseries H\arabic*)}]
    \item \label{H1} The functions $\rJ$, $\FmaxA$, $\Fmaxj:[0,1]\rightarrow [0,\infty)$ are bounded and globally Lipschitz continuous.
    \item \label{H2} The jump coefficients $\phiA, \phiJ: [0,\infty)\times \R\rightarrow\R$ are chosen such that there exist non-negative functions $g^{(\XA)}, g^{(\XJ)}$ in $L^2(\nu)$ for which:
    \begin{enumerate}[start=1,label={(\bfseries H2.\roman*)}]
    \item\label{hip:linear-growth-phi} Linear growth: For all $(x,z)\in[0,+\infty)\times \R$, we have:
    \begin{align*}
        |\phiA(x,z)|&\leq |x|\,g^{(\XA)}(z),\\
        |\phiJ(x,z)|&\leq |x|\,g^{(\XJ)}(z).
    \end{align*}
    \item\label{hip:lipschitz-condition-phi}  Lipschitz condition: For all $(x,z),(x',z)\in[0,+\infty)\times \mathbb{R}$, we have
     \begin{align*}
     \int_{-\infty}^\infty |\phiJ(x,z)-\phiJ(x',z)|^2\nu(dz)&\leq L_{\phiJ}|x-x'|^2,\\
       \int_{-\infty}^\infty |\phiA(x,z)-\phiA(x',z)|^2\nu(dz)&\leq L_{\phiA}|x-x'|^2,
    \end{align*}
     for some constants $L_{\phi_{(\cdot)}}>0$ depending on the $L^2(\nu)$-norm of $g^{(\cdot)}$.
     \item \label{hip:jumps_positiveness}  Lower control of jumps:
\[\left(x+\phiJ(x,z)\right)\wedge\left(x+\phiA(x,z)\right)\geq 0,\quad \forall \, (x,z)\in [0,\infty)\times \R.\]

    \end{enumerate}
    \item \label{H3} The constants $\overline{\beta}_0$ and $\overline{\beta}_1$ in \eqref{eq:def_beta_1}-\eqref{eq:def_beta_2} are non-negative and satisfy 
    $\overline{\beta}_0\wedge\overline{\beta}_1 > \sE^2/2$.
    \item \label{H4} Coefficients $\mu$ and $\sigma_P$ in Equation \eqref{eq:full-model} are globally Lipschitz continuous and such that $P_t(\omega)\in\R_+$ for all $(t,\omega)\in [0,+\infty)\times \Omega$.
\end{enumerate}

Hypotheses \ref{H1}, \ref{hip:linear-growth-phi} and \ref{H4} are considered in the analysis of the existence of a weak solution for the system \eqref{eq:full-model} and control of moments, while conditions \ref{hip:lipschitz-condition-phi}-\ref{hip:jumps_positiveness} are included to guarantee the existence of strong solutions and non-negativity of the components of process $\K$. Condition \ref{H3} allows us to give a criterion of non-absorption within complete or incomplete compliance. Notice that Assumption \ref{hip:linear-growth-phi} can be easily generalized with a linear form in $x$ and not just affine. However, we have considered an affine form for the growth condition of $\phiA$ and $\phiJ$ because the jump term represents the effect of extreme phenomena that may, for example, annihilate just a fraction of the existing population. 

We now state the existence and uniqueness of a strong admissible solution to model \eqref{eq:full-model}.
The proof of this result is postponed to Section \ref{sec:well-posedness}.

\begin{theorem}\label{th:well-posedness}
Consider a given filtered probability space $(\Omega, \mathcal{F}, (\mathcal{F}_t:t\ge0),\mathbb{P})$ where we can define the independent standard Brownian motions $W^{(J)}$ , $W^{(A)}$, $W^{(E)}$, $ W^{(P)}$ and the independent Poisson random measure $\Np$ with intensity measure $\nu(dz)\cdot \lambda dt$, being $\nu$ a probability measure on $(\mathbb{R}, \mathcal{B} (\mathbb{R}))$. Then 
there exists a $\mathcal{S}$-valued and $\mc{F}_t$-adapted Markov process $(\ZZ_t: (\K_t, E_t, P_t), 0 \le t < \zeta)$ pathwise solution of \eqref{eq:full-model}. If $\ZZ_0\in \mathcal{S}$ and  $\mathbb{P}_{\bf x}$ denotes the probability measure under which $\ZZ_0 = {\bf x}$, then $\mathbb{P}_{\bf x} (\zeta < \infty)=0$.

Moreover, for any non-negative integer $p$ such that $g^{(\cdot)}\in L^p(\nu)$ and $\ZZ_0\in L^p(\Omega; [0,\infty)^2\times [0,1]\times [0,\infty))$, we have that the $p$th-moment -in finite time intervals- of the process $\ZZ$ is finite.

\end{theorem}

By reinforcing the assumptions in Theorem \ref{th:well-posedness}, we demonstrate that the origin $(0,0)$ is an inaccessible state for the populations $\K$ and that there exist
an invariant probability measure. This result, summarized in the following statement, underpins the long-term analysis of the system. For more details see Section \ref{sec:inaccessibility_origin_invariant}.

\begin{theorem}\label{th:assymptotics}
Assume the hypotheses of Theorem~\ref{th:well-posedness} hold, together with additional conditions \ref{ass:no_jump_near_0_supp}-\ref{ass:no_jump_near_0_no_excess}. Then, we have that \[
\mathbb{P}_{\bf x} \left( \exists t > 0 : \K_t = (0,0) \right) = 0,
\]
and the process $\left(\K_t; t\geq0\right)$ admits an invariant probability measure on \((0,\infty)^2\).

\end{theorem}

Having established the existence of an invariant measure under suitable conditions, we turn to the regime in which the population undergoes asymptotic extinction. The following result identifies mild sufficient conditions for long-term extinction, thus characterizing a stability regime.

\begin{theorem}
Assume the conditions of Theorem~\ref{th:assymptotics} hold. Then, under mild condition on the parameters of the model (see details in Section \ref{prop:extinction}) the total population $(\XJ_t+\XA_t;t\geq0)$  gets asymptotically extinct exponentially fast. 
\end{theorem}

\medskip

In this work, we also provide a numerical scheme to simulate \eqref{eq:full-model} and prove its strong convergence. The proof of the following theorem, and the details of its construction are presented in Section \ref{sec:numerical_scheme} and Appendix  \ref{sec:delayed-proofs-ns}.

Since the process $P$ is exogenous, we assume that there is a numerical scheme $\bar{P}$, with bounded moments, such that
\begin{equation}\label{eq:scheme-for-P}
\begin{aligned}
\sup_{0\leq t\leq T}\E\left[(P_t - \bar P_t)^{2p} \right] &\leq  C\dt^p.
\end{aligned}
\end{equation}

\begin{theorem}\label{thm:convergence-numerical-scheme} For a fix time horizon $T>0$, let $\dt = T/N $ with $N\in \N$ and consider a temporal grid $t_k=k\dt$, $\eta(t):=\sup\{t_k:t_k< t\}$, $\ddt:=t-\eta(t)$, and $\bZZ_{0}:=\ZZ_0=(\XA_0,\XJ_0,E_0,P_0)$, the initial condition of \eqref{eq:full-model}. Given $(\scXJ_{\eta(t)},\scXA_{\eta(t)},\bE_{\eta(t)})$, for any $t\in[0,T]$, we define $(\scXJ_{t},\scXA_{t},\bE_{t})$ iteratively by

\begin{equation}\label{eq:Numerical-Scheme}
    \begin{aligned}
    \scXJ_{t}  &=\left(\scXJ_{\eta(t)}  -\ddt \kappaJ(\bE_{\eta(t)})\scXJ_{\eta(t)} \right) \exp\left(-\frac{\sJ^2}{2}\ddt + \sJ\deltaWJ{\eta(t)}{t}\right) \\
&\qquad+    \rhoJ(\scXJ_{\eta(t)} ,\scXA_{\eta(t)},\bE_{\eta(t)})\scXA_{\eta(t)} \exp\left(-\frac{\sJ^2}{2}\ddt + \sJ\deltaWJ{\eta(t)}{t}\right) +   \ind{\{\eta (t) \le \tau^{1}_{\eta(t)} < t\}} \phiJ \left(\scXJ_{\eta(t)}, \xi_1\right),\\
\scXA_{t} & = 
  \left(\scXA_{\eta(t)}-\ddt \kappaA(\bE_{\eta(t)})\scXA_{\eta(t)}  \right) \exp\left(-\frac{\sA^2}{2}\ddt + \sA\deltaWA{\eta(t)}{t}\right) \\
&\qquad+    \rhoA\scXJ_{\eta(t)}\exp\left(-\frac{\sA^2}{2}\ddt + \sA\deltaWA{\eta(t)}{t}\right)  +   \ind{\{\eta(t) \le \tau^{1}_{\eta(t)} < t\}} \phiA \left(\scXA_{\eta(t)}, \xi_1\right),\\
          \bE_t &= \cutDt\Bigg(\bE_{\eta(t)} +\left[\beta_1(\bZZ_{\eta(t)})(1-\bE_{\eta(t)})  -\beta_0(\bZZ_{\eta(t)})\bE_{\eta(t)} \right]\ddt         +\sE\sqrt{\bE_{\eta(t)}(1-\bE_{\eta(t)})}\deltaWE{\eta(t)}{t} \Bigg),
    \end{aligned}
\end{equation}
where  $\tau^1_{\eta(t)}$ is an exponential random variable of rate $\lambda$ and $\xi_1$ is a random variable with law $\nu$. \\
  
Under the hypotheses \ref{H1}-\ref{H4}, it holds:
\begin{equation}\label{eq:convergence-strong-error}
\begin{aligned}
 \lim_{\dt\to0}\sup_{0\leq t\leq T}\left\{\E\left[\left(\XJ_t-\scXJ_t \right)^2\right] + \E\left[\left(\XA_t-\scXA_t \right)^2\right] + \E\left[\left(E_t-\bE_t \right)^2\right] \right\} =0.
\end{aligned}
\end{equation}

\end{theorem}

\section{Derivation of the model and well-posedness }\label{sec:derivation}

The aim of this section is twofold: to fill the gap between the idealized harvester-population dynamics schematically depicted in the introduction and the proxy system \eqref{eq:simplified-model}, and to prove the well-posedness of its enriched version \eqref{eq:full-model}. The first result is, in a sense, quite independent of the properties we will prove about the SDE \eqref{eq:full-model} but, as the reader will verify, it depends on the well-posedness of the system \eqref{eq:simplified-model}.

\subsection{From individual dynamics to an intermediate SDE: mean-field approximation}

We consider here the pure jump process 
$(\ZZ^{(n)}_t: t \ge 0):=(\dfrac{1}{n}\K^{(n)}_t, \dfrac{1}{n}E^{(n)}_{t}: t \ge 0),$ as described in Section \ref{sec:model_description}, with the sequence of graphs being the deterministic sequence of complete graphs. In this case, it is clear that the process $\ZZ^{(n)}$ is in fact a Markov process. If $S^{(n)}:= \mathbb{Z}^{2}_+(n) \times \dfrac{\{0, 1, \ldots, n\}}{n}$ is the state-space of this process and a typical element of it is denoted $\zz = (j, a, e)$, the transitions of $\ZZ^{(n)}$ are summarized as follows. From the state $(j, a, e)$, the process will jump to:
\begin{itemize}[leftmargin=15pt]
    \item $ (j-1/n, a+1/n, e) \text{ at rate } j \rho_A(\zz)$.
    \item $(j+1/n, a, e) \text{ at rate } a \rho_J(\zz)$.
 \item $(j, a-1/n, e) \text{ at rate } a \kappa_A(\zz)$.
\item  $(j, a-1/n, e) \text{ at rate } j \kappa_J(\zz)$.
\item $(j, a, e+1/n)$ at rate 
$\gamma n (1-e) \left ( \dfrac{n e }{n-1} (1- \beta_0 (\zz)/n)+  \dfrac{n(1-e)-1}{n-1} \beta_1 (\zz)/n )   \right).$
\item $(j, a, e-1/n) \text{ at rate } \gamma n e \left ( \dfrac{n (1-e)}{n-1} (1- \beta_1 (\zz)/n)+  \dfrac{ne-1}{n-1} \beta_0 (\zz)/n )   \right)$. 

\end{itemize}

The next theorem states the convergence in distribution of the speed-up process $\tilde \ZZ^{(n)}$ given by $\tilde \ZZ^{(n)}_t:= \ZZ^{(n)}_{nt}$, to the unique solution of the system \eqref{eq:simplified-model}.  

\begin{theorem}
Assume that \ref{H1}-\ref{H2}-\ref{H3}-\ref{H4} hold and the network structure is, for each $n$, the complete graph on $n$ vertices. Then the sequence of probability measures $(\textbf{Law}(\tilde \ZZ^{(n)}))_{n \ge 1}$ on the Skorokhod space $\mathbb{D}([0, +\infty); \R^2_+ \times [0,1])$, converges weakly to the unique martingale solution of \eqref{eq:simplified-model}. 
\end{theorem}

\begin{proof}
For a bounded function, $f: S^{(n)} \mapsto \R$, the generator associated with the Markov process $\ZZ^{(n)}$ acts as:
\begin{align*} \label{eq:n-generator}
\mathcal{L}^{(n)} f(\zz) & = j \rho_A  (\zz)\Big(f(j-1/n, a+1/n, e)-f(\zz)\Big) + a \rho_J  (\zz) \Big(f(j+1/n, a, e)-f(\zz)\Big) \\
& \quad + a \kappaA(\zz)  \Big(f(j, a-1/n, e)-f(\zz)\Big) + j \kappaJ(\zz) \Big(f(j-1/n, a, e)-f(\zz)\Big) \\
& \quad +  \gamma n (1-e) \left ( \dfrac{n e }{n-1} (1- \beta_0 (\zz)/n)+  \dfrac{n(1-e)-1}{n-1} \beta_1 (\zz)/n )   \right)\\
&\quad+ \gamma n e \left ( \dfrac{n (1-e)}{n-1} (1- \beta_1 (\zz)/n)+  \dfrac{ne-1}{n-1} \beta_0 (\zz)/n )   \right)  \times \Big(f(j,a, e-1/n)-f(\zz)\Big). 
\end{align*}
Thus, for suitable functions $f$ (for example, $f\in \mathcal{C}^\infty(\R^2 \times [0,1]; \R)$, the space of infinitely-differentiable functions with compact support), the generator of $\tilde \ZZ^{(n)}$ acts as
 \(\tilde{ \mathcal{L}}^{(n)}f = n\mathcal{L}^{(n)}f.\)

A simple computation shows that, at least formally, if the sequence of processes $(\tilde \ZZ^{(n)})$ converge in law, then the limit process should have the generator $\mathcal{L}$ associated to \eqref{eq:simplified-model}. Under the assumptions \ref{H1}-\ref{H2}-\ref{H3}-\ref{H4} we will see below, as a consequence of Theorem \ref{th:well-posedness}, that the $(\mathcal{L}, \mu)$ martingale problem is well-posed for any initial measure $\mu$ with compact support. So, in order to prove the claimed convergence, it remains to check:
\begin{itemize}
\item The tightness of the family of laws of the processes $(\tilde \ZZ^{(n)})$ (as probability measures on the Skorokhod space).
\item The convergence of the finite dimensional distributions of the processes $\tilde \ZZ^{(n)}$ to those of the process $\ZZ$.  
\end{itemize}
Incidentally, the hypotheses of Th. 4.21, Chapter IX, pag. 558, of \cite{jacodshiryaev2003} hold in our case. So, it remains to prove the convergence of the characteristics of the sequence of semimartingales. To this end, consider the kernel $K^{(n)}$ given by:
\begin{align*}
K^{(n)}(\zz, \cdot)&= n j \rho_A  (\zz) \delta_{-1/n, 1/n, 0}+ n a \rho_J  (\zz)\delta_{1/n, 0 , 0} + n a \kappaA(\zz)  \delta_{0, -1/n, 0} + n j \kappaJ(\zz) \delta_{-1/n, 0, 0}\\
& \quad +  n^2\gamma  (1-e) \left ( \dfrac{n e }{n-1} (1- \beta_0 (\zz)/n)+  \dfrac{n(1-e)-1}{n-1} \beta_1 (\zz)/n )   \right) \delta_{0,0, 1/n} \\
& \quad + \gamma n^2 e \left ( \dfrac{n (1-e)}{n-1} (1- \beta_1 (\zz)/n)+  \dfrac{ne-1}{n-1} \beta_0 (\zz)/n )   \right) \delta_{0, 0, -1/n}.
\end{align*}
Hence:
\begin{align*}
\tilde{\mathcal{L}}^{(n)} f (\zz) = \int_{\R^3} K^{(n)} (\zz, dy) (f(\zz+ y)- f(\zz)).
\end{align*}
Given $\varepsilon> 0$, for every $n \ge \sqrt{3}\lfloor \varepsilon^{-1}\rfloor+1$, we have for every $R > 0$:
\begin{align*}
\sup_{\vert \zz\vert < R }\int_{\R^3} K^{(n)} (\zz, d\yy) \vert \yy\vert^2 \mathbf{1}_{\vert \yy \vert > \varepsilon}=0,
\end{align*}
and thus condition (ii) of the cited Theorem hold. Moreover, we compute easily:
\begin{align*}
\int_{\R^3} K^{(n)} (\zz, d\yy) y_1 & = a \rho_J (\zz) - j (\kappaJ (\zz)+ \rho_A (\zz)). \\
\int_{\R^3} K^{(n)} (\zz, d\yy) y_2 & = j \rho_A (\zz)-a \kappaA(\zz).\\
\int_{\R^3} K^{(n)} (\zz, d\yy) y_3 & = - \gamma \beta_0(\zz) e + \gamma \beta_1 (\zz) (1-e).\\
\int_{\R^3} K^{(n)} (\zz, d\yy) y_iy_k & = O (1/n) \quad i \text{ or } k \text{ equals } 1 \text{ or } 2, \text{ uniformly on compacts.}\\
\int_{\R^3} K^{(n)} (\zz, d\yy) y^2_3 & = 2 \gamma e (1-e)+ O(1/n) \quad \text{ uniformly on compacts.}\\
\end{align*}

Taking the limit as $n \to \infty$, we obtain the convergence to the corresponding characteristics of the semimartingale determined by the unique solution of \eqref{eq:simplified-model}. Thus, condition (i) of the cited Theorem hold, and the claim follows.
\end{proof}

\begin{remark}
The previous computations indicate that the biomass and harvesters' behavior, described by the macroscopic model \eqref{eq:simplified-model}, represent the limit behavior of two sets of discrete particles as they scale proportionally. This fact, which cannot be easily derived just on the basis of our target SDE, becomes important if we aim at applying the limiting model to some real-life situation, since it indicates that, in order to be suitable as a prediction tool, we must somehow ensure (possibly based upon socio-ecological considerations) that harvesters and biomass availability are approximately in linear relation, which is assumed hereafter.     

\end{remark}

\subsection{Well-posedness and positivity of the process}\label{sec:well-posedness}

As announced in Section \ref{sec:model_description}, after modeling consideration, we arrive to the full system \eqref{eq:full-model}, for which we need to establish the well-posedness. This is the aim of this section. For technical reasons we will first consider a slightly modified version of \eqref{eq:full-model}, namely we substitute $E$ (resp. $1-E$) by $\Psi(E)$ (resp. by $\Psi(1-E)$). Notice that the coefficient $\beta_1(\ZZ)(1-E)  -\beta_0(\ZZ)E$ can be written as:
\begin{align*}
\beta_1(\ZZ)(1-E)  -\beta_0(\ZZ)E 
&=     \Delta\tbeta(\K,P)\, E(1-E) + \bar\beta_{1} -(\bar\beta_{0}+\bar\beta_{1})E, 
\end{align*}
where
\begin{equation}\label{eq:beta_tilde}
\begin{aligned}
    \Delta\tilde{\beta}(\K,P)
    &:= \tauU \mathcal{I}^{\eta}_{P_{\min}}(P + s) + \tauU  \mathcal{I}^{\eta}_{P_{\min}}(P) - \tauU \mathcal{I}^{\eta}_{P_{\max}}(P)  - \tauU \Psi\left(\frac{\XJ+\XA
    }{\KA}\right) - \tauU \Psi\left(\frac{\XJ+\XA}{\KA}\right).
\end{aligned} 
\end{equation} 
Now, if ${\bf x} = (x,y,e,p)$ stands for the canonical state vector, we consider the auxiliary system
\begin{equation*}
\left\{\begin{aligned}
    d\XJ_t &= \left[\rrJ(\K_t, E_t)\XA_t - \kappaJ(\Psi(E_t))\XJ_t \right]dt  + \sJ\XJ_tdW^{(J)}_t  + \int_{\R\setminus\{0\}} \phiJ (\XJ_t, z)\Np(dz\otimes dt)  \\
    d\XA_t &= \left[\PjA \XJ_t - \kappaA(\Psi(E_t))\XA_t \right]dt + \sA\XA_t dW^{(A)}_t+ \int_{\R\setminus\{0\}} \phiA (\XA_t, z)\Np(dz\otimes dt) \\
    dE_t &= \Psi(E_t)\Psi(1-E_t)\,\Delta\tbeta(\K_t,P_t)dt   +\Psi(1-E_t)\overline{\beta}_2dt-\Psi(E_t)\overline{\beta}_1dt\\
    &\quad+  \sE\sqrt{\Psi(E_t)\Psi(1-E_t)}dW^{(E)}_t \\
    dP_t &= \mu(P_t)dt +\sP(P_t)  dW_t^{(P)}, 
\end{aligned} 
\right.
\end{equation*}
with initial condition $\ZZ_0 = (\XJ_0,\XA_0,E_0,P_0)\in L^2(\Omega; \R_+^2\times [0,1]\times \R_+)$, and $\Delta\tilde{\beta}(x,y,p)$ as defined in \eqref{eq:beta_tilde}.

\medskip
We notice that $x\mapsto \Psi(x)$ is a Lipschitz continuous map with Lipschitz constant $1$. Similarly, $\mathcal{I}^{\eta}_{P_0}(\cdot)$ is Lipschitz in $\mathbb{R}_+$, which made $\Delta \tilde{\beta}$  also Lipschitz and bounded by $2$. 

For simplicity of notation, we write our auxiliary system in vector form
\begin{equation}\label{eq:system_welldef}\tag{\~M}
d\ZZ_t = b(\ZZ_t)dt +\sigma(\ZZ_t)dW_t + \int_{-\infty}^{\infty} \Phi(\ZZ_t,u)\Np(du\otimes dt),
\end{equation}
where $W_t = [ \WJ, \WA,\WE, \WP]$ is a $4$-dimensional standard Brownian motion, the diffusion matrix is given by $\sigma(\mathbf{x}) = \texttt{diag}(\sJ,\sA,\sE \sqrt{\Psi(e)\Psi(1-e)}, \sP(p))$, $\Phi(\mathbf{x}) =(\phiJ(x,z),\phiA(x,z),0,0)^t $, and
\begin{equation*}
b(\mathbf{x}) = \left[\begin{aligned}
\rrJ(x,y,\Psi(e))y - \kappaJ(\Psi(e))x&\\
\PjA x - \kappaA(\Psi(e))y&\\
\Psi(e)\Psi(1-e)\Delta \tilde{\beta}(x,y,p) +\Psi(1-e)\overline{\beta}_2 - \Psi(e)\overline{\beta}_1&\\
\mu(p) \end{aligned} \right].
\end{equation*}

\paragraph*{A priori control of moments.}

Notice that under assumption \ref{H4}, the price process has finite moments of every order. The proofs of the following preliminary results are rather classical so we postpone them to the Appendix \ref{sec:proof-sec-well-posedness}.

\begin{proposition}\label{prop:moments}
Assume there is a strong solution to the system \eqref{eq:system_welldef} and assumptions \ref{H1},\ref{H2} and \ref{H4} hold. Then, provided that $g^{(\cdot)}\in L^p(\nu)$ and $\ZZ_0\in L^p(\Omega; [0,\infty)^2\times [0,1]\times [0,\infty))$, for some $p\geq1$, the $p$th-moment of the process $\ZZ$ is finite, i.e. $\sup_{0\leq t\leq T}\E[\|\ZZ_t\|^p]<\infty$.
\end{proposition}

\begin{corollary}\label{cor:sup-norm_moments}
Assume \ref{H1},\ref{H2} and \ref{H4} hold true, and consider $g^{(\cdot)}\in L^p(\nu)$ and $\ZZ_0\in L^p(\Omega; [0,\infty)^2\times [0,1]\times [0,\infty))$, for some $p\geq1$. Then, provided that there exists a strong solution to the system \eqref{eq:system_welldef}, \[\E[\sup_{0\leq t\leq T}\|\ZZ_t\|^p]<\infty.\]
\end{corollary}

As we will show later, the uniform moment control can be easily extended to the whole real line in the case $p=1$.

\paragraph*{Weak existence}
We notice that the difficulty in the direct application of classical results to prove the well-posedness of the SDE lies in the singularity of the diffusion of the compliance process at the boundaries of its domain. We overcome this limitation by means of comparison arguments.

\begin{proposition}\label{prop:weak_existence}
If assumptions \ref{H1}-\ref{hip:linear-growth-phi}-\ref{H3} hold true, then there exists a weak solution for \eqref{eq:full-model}. Furthermore, if there exists a strong solution, then $E_t\in(0,1)$ a.s. for all $t\geq 0.$
\end{proposition}
\begin{proof}

Let $n \ge 1$ be an arbitrary positive integer, and define the stopping time $\tau_n:= \inf\{t \ge 0: E_t (1-E_t) < 1/n\}$. Let $C>0$ be a bound for $\Delta \widetilde{\beta}$, as defined in \eqref{eq:beta_tilde}, and consider the one-dimensional diffusions:

\begin{align*}
dE^{(u)}_t &= \Big( \big(C\, \Psi(E_t^{(u)}) + \bar{\beta}_2\big)\Psi(1-E^{(u)}_t)- \overline \beta_0 \Psi(E^{(u)}_t)\Big)dt \\
    &\quad+ \sigma_E \sqrt{\Psi(E^{(u)}_t)\Psi(1-E^{(u)}_t)}d\WE_t. \\
dE^{(l)}_t &= \Big(\overline \beta_1 \Psi(1-E^{(l)}_t)- \big(\bar{\beta}_1 + C\, \Psi(1-E_t^{(l)})\big) \Psi(E^{(l)}_t)\Big)dt \\
    &\quad+ \sigma_E \sqrt{\Psi(E^{(l)}_t)\Psi(1-E^{(l)}_t)}d\WE_t,
\end{align*}
where $\WE$ is the Brownian motion driving the dynamics of $E$. We assume for $E^{(u)}$ and $E^{(l)}$ the same initial condition as for $E$, i.e. $E^{(u)}_0 = E_0$ and $E^{(l)}_0 = E_0$. 

Let $\tilde \ZZ$ be the solution of \eqref{eq:system_welldef}. Then, for each $n\geq1$, the stopped process $(\tilde\ZZ^{\tau_n}, (E^{(l)})^{\tau_n}, (E^{(u)})^{\tau_n})$ is a Markov process, and the same is true for the component $\tilde \ZZ^{\tau_n}$. Moreover, for every $n\geq1,$ the diffusion coefficient $\sqrt{\Psi(x) \Psi(1-x)}$ associated to the components $E$, $E^{(u)}$ and $E^{(l)}$ are Lipschitz on $\{x\in\R:~x(1-x)> 1/n\}$. Hence, by general results on martingales solutions, there exist a unique solution for the martingale problem $(\mu, \mathcal{L}, \tau_n)$ associated to an initial probability measure $\mu$, the generator $\mathcal{L}$ and the stopped process $\tilde \ZZ^{\tau_n}$.

We claim that the sequence of stopping times $(\tau_n;n\geq1)$ goes to infinity a.s. Then, in virtue of \cite[Theorem 6.3, pp. 219]{ethierkurtz86}, we deduce that there exist a unique solution for the martingale problem associated to $(\mu, \mathcal{L}, +\infty)$. In turn, under Assumptions \ref{H1}-\ref{hip:linear-growth-phi},\ref{H4}, we can apply \cite[Theorem 2.3]{Kurtz2011} implying that a weak solution to the SDE \eqref{eq:system_welldef} exists.

It remains to demonstrate that $\tau_n$ converges to infinity a.s. Indeed, by standard comparison arguments (see e.g., \cite[Theorem 1.1]{ikeda1977}) we have that for any random initial condition $E_0$ with support in $(0,1)$, the component process $E$, starting from $E_0$, is sandwiched between $E^{(l)}$ (from below) and $E^{(u)}$ (from above), whenever these processes are issued from the same initial condition $E_0$. Furthermore, in order to analyze the boundaries of these process we can use the Feller classification of boundaries (see e.g., 
\cite[Proposition 5.22, pp. 345]{karatzas98}) as follows: For any $c\in(0,1)$ consider the scale function for the process $E^{(u)}$
\begin{equation*}
p^{(u)}(x):=\int_c^x \exp\left\{-\frac{2}{\sE^2}\int_c^y\frac{C\,\Psi(z)\Psi(1-z) + \overline{\beta}_2\Psi(1-z) - \overline{\beta}_1\Psi(z)}{\Psi(z)\Psi(1-z)}dz\right\}dy.
\end{equation*}

Since the domain of interest for the scale function is $x\in(0,1)$, it is straightforward that the cutoff has no effect on its behavior and there exists a constant $K>0$, depending on the parameters of the model and $c$, such that
\begin{equation*}
p^{(u)}(x) = K\, \int_c^x e^{-\tfrac{2C\, y}{\sE^2}}\exp\left\{-\frac{2\overline{\beta}_2}{\sE^2}\ln(y) -\frac{2\overline{\beta}_1}{\sE^2}\log(1-y)\right\}dy.
\end{equation*}

Then, it is easy to check that $p^{(u)}(0^+) = -\infty$ and $p^{(u)}(1^-) = +\infty$ under $1-\frac{2\overline{\beta}_2}{\sE^2}<0$ and $1-\frac{2\overline{\beta}_1}{\sE^2}<0$, respectively. A similar remark can be done for the scale function of the process $E^{(l)}$. From this we obtain that, under Assumption \ref{H3}, the boundaries $0$ and $1$ are entrance boundaries for $E^{(u)}$ and $E^{(l)}$. Thus, \[0<E^{(l)}\leq E^{\tau_n}\leq E^{(u)}<1, \mbox{ a.s. for all }n\geq1,\] and $\tau^{(n)}$ necessarily goes to infinity almost surely. We conclude that the SDE \eqref{eq:system_welldef} has weak solution, with $E_t$ having support on $(0,1)$ for all $t\in[0,T]$. In particular, there exist a weak solution to the SDE system \eqref{eq:full-model}, as desired.

\end{proof}

\paragraph*{Strong well-posedness} To prove Theorem \ref{th:well-posedness} we need the following proposition regarding the path-wise uniqueness of solutions for \eqref{eq:full-model}.  
Under boundedness and Lipschitz conditions of the rates in Assumption \ref{H1}, the proof follows closely the ideas of the Theorem by Yamada and Watanabe, with the subtlety we are considering  
a multidimensional model. We postpone the proof to Appendix \ref{sec:proof-sec-well-posedness}.

\begin{proposition}\label{prop:pathwise_uniqueness} Assume \ref{H1}, \ref{H2}, \ref{H3} and \ref{H4} hold true, and Equation \eqref{eq:full-model} has a c\`adl\`ag strong solution. Then, path-wise uniqueness hold for \eqref{eq:full-model}.
\end{proposition}

From Propositions \ref{prop:moments}, \ref{prop:weak_existence} and \ref{prop:pathwise_uniqueness} we can address the proof of  the strong well posedness of our model.

\begin{proof}[Proof of Theorem \ref{th:well-posedness}]

Assume \ref{H1},\ref{H2}, \ref{H3} and \ref{H4}. Then, from Proposition \ref{prop:moments}, \ref{prop:weak_existence} and \ref{prop:pathwise_uniqueness} we can apply Yamada \& Watanabe's result (see e.g., \cite[Proposition 3.20, pp. 309]{karatzas98}) to conclude that, for any initial condition $\ZZ_0\in L^2(\Omega; [0,\infty)^2\times [0,1]\times [0,\infty))$, there exists a unique strong solution of the SDE \eqref{eq:full-model}. Moreover, for any $p\geq2$ such that $g^{(\cdot)}\in L^p(\nu)$ and $\ZZ_0\in L^p(\Omega; [0,\infty)^2\times [0,1]\times [0,\infty))$, there exists a constant $C$, depending on the parameters of the model and terminal time $T$, with 
\[\sup_{0\leq t\leq T}\E[\|\ZZ_t\|^p]<C(T).\]

We claim that under \ref{hip:jumps_positiveness} the components of the process $\K$ are non-negative, i.e. $\ZZ_t\in \mathcal{S}$ for all $t\in[0,T]$ and Theorem \ref{th:well-posedness} holds true.

Indeed, for $\tZZ=(\tXJ,\tXA,E,P)$ the solution of the system without jumps:
\begin{equation}\label{eq:auxiliary-model-without-jumps}
\begin{aligned}
    \tXJ_t
    &=\tXJ_0+\int_0^t\rJ(E_s)\Psi\left(1-\frac{\tXJ_s+\tXA_s}{\KA}\right)\tXA_s   -\left(\PjA + \Fmaxj(E_s) + \mJ\right)\tXJ_sds + \int_0^t\sJ\tXJ_sdW^{J}_s,  \\
    \tXA_t &= \tXA_0 +\int_0^t \PjA\tXJ_sds - \int_0^t \left(\mA + \FmaxA(E_s)\right)\tXA_sds   + \int_0^t\sA\tXA_sdW^{A}_s, \\
  E_t &=E_0 + \int_0^t\left(\Delta \tilde{\beta}(\tXJ_s,\tXA_s,P_s)E_s(1-E_s) +\overline{\beta}_2(1-E_s) - \overline{\beta}_1E_s\right)ds +\int_0^t\sE\sqrt{E_s(1-E_s)}dW^{(E)}_s.
\end{aligned} 
\end{equation}

Then, from the sequence $(\psi)_{n\geq1}\subset \mathcal{C}^2(\R;\R)$ in the proof of Proposition \ref{prop:pathwise_uniqueness}, we define for all $n\geq1$, $\varphi_n(x)= \psi_n(x) \ind{(0,\infty)}(x)$ which is smooth, non-negative and satisfies $\|\psi_n'(x)\|\leq 1,$ and $\varphi_n''(x)\leq \frac{2}{nx}\ind{(0,\infty)}(x)$, with $\lim_{n\to\infty}\varphi_n(x)= (x\lor 0)$. For fixed $\K_0\in \R_+\times \R_+$ we apply It\^o formula, obtaining:
\begin{align*}
&\E[\varphi_n(-\tXJ_t)+\varphi_n(-\tXA_t)] =    \E[\varphi_n(-\tXJ_0)+\varphi_n(-\tXA_0)] \\
&\qquad + \E\left[\int_0^t -\varphi_n'(-\tXJ_s)\rJ(E_s)\Psi\left(1-\frac{\tXJ_s+\tXA_s}{\kA}\right)\tXA_sds\right]\\
&\qquad + \E\left[\varphi_n'(-\tXJ_s)\left(\PjA + \Fmaxj(E_s) + \mJ\right)\tXJ_s ds\right] \\
  & \qquad   + \E\left[\int_0^t\left\{-\varphi_n'(-\tXA_s) \PjA\tXJ_s+\varphi_n'(-\tXA_s)\left(\mA +\FmaxA(E_s)\right)\tXA_s \right\}ds\right] \\
  &   \qquad    +\frac{1}{2} \E\left[\int_0^t \left\{\varphi_n''(-\tXA_s)\sA^2\tXA_s^2+ \varphi_n''(-\tXJ_s) \sJ^2\tXJ_s^2\right\}ds\right].
\end{align*}
Notice that, for all $x\in \R$, we have $\varphi'_n(-x)x\leq 0$. Thus, we can bound the terms using the properties of the approximation sequence $\varphi_n$ and assumption \ref{H1}:
\begin{align*}
&\E[\varphi_n(-\tXJ_t)+\varphi_n(-\tXA_t)]  \leq    \E[\varphi_n(-\tXJ_0)+\varphi_n(-\tXA_0)]  \\
  &  \quad + \E \left[\int_0^t\left\{ -\varphi_n'(-\tXJ_t)\rJ(E_s)\Psi\left(1-\frac{\tXJ_s+\tXA_s}{\kA}\right)\tXA_s -\varphi_n'(-\tXA_t) \PjA\tXJ_s\right\}\ind{\tXA_s\vee \tXJ_s<0}ds\right] \\
  &    \qquad  +\frac{\sA^2\vee \sJ^2}{n}\E \int_0^t \left((-\tXA_s)\vee 0+ (-\tXJ_s)\vee 0\right)ds,
\end{align*}
obtaining
\begin{align*}
&\E[\varphi_n(-\tXJ_t)+\varphi_n(-\tXA_t)]  \leq    \E[\varphi_n(-\tXJ_0)+\varphi_n(-\tXA_0)] \\
& \qquad \leq \E[\varphi_n(-\tXJ_0)+\varphi_n(-\tXA_0)]  + C\left(1+\frac{1}{n}\right)\E \int_0^t \left((-\tXA_s)\vee 0+ (-\tXJ_s)\vee 0\right)ds,
\end{align*}
for some constant $C>0$. 
Therefore, from Fatou's and Gronwall's lemma we can conclude that
\begin{align*}
\E[(-\tXJ_t)\vee 0+(-\tXA_t)\vee0]= 0,
\end{align*}
i.e., for all $t\geq0$, $\tXJ_t$ and $\tXA_t$ are almost surely non-negatives.

Now we incorporate the jumps to the dynamic and show that the positivity still holds. For this we consider the first jump time $\tau_1$. We have shown that, for all $t<\tau_1$, $\XA_t\wedge \XJ_t \geq 0$ almost surely. Then, from \ref{hip:jumps_positiveness}:
\begin{align*}
&\left(x+\phiJ(x,z)\right)\wedge\left(x+\phiA(x,z)\right)\geq 0,\forall \, (x,z)\in [0,\infty)\times \R,
\end{align*}
we get
\( \XJ_{\tau_1^+} = \XJ_{\tau_1^-} + \phiJ(\XJ_{\tau_1^-}, \xi) \geq 0 \), with $\xi\sim \nu(dz)$ and similarly \(\XA_{\tau_1^+}\geq 0\). Thus, if we denote $(\tau_n)_{n\geq 1}$ the sequence of jump times, and consider the initial condition $\XJ_{\tau_n^+}$ we can conclude that, for all $\tau_n^+\leq t<\tau_{n+1}^{-}$, the trajectories are almost surely non-negatives. \end{proof}


\section{Asymptotic behavior of the model}\label{sec:inaccessibility_origin_invariant}

In this section, we establish the existence of a non-trivial invariant measure for the solution $\ZZ$. In particular, under assumption \ref{hip:jumps_positiveness} we have proved the process $\K$ stays in $[0,\infty)^2$ (see Theorem~\ref{th:well-posedness}), however it does not exclude the possibility that the population component $\K$ reaches the absorbing state $(0,0)$ in finite time. The inaccessibility of the origin is not only ecologically relevant, but also essential from a modeling perspective: due to the structure of the dynamics for $E$, the system depends on the presence of biomass —without it, the notion of compliance is superfluous. The following result clarifies this, and serves to motivate our search for invariant measures and conditions for extinction that are not supported on such trivial configurations.

\begin{proposition}[Inaccessibility of the origin for the population component]\label{prop:non-attainability-origin}
Assume \ref{H1},\ref{H3}, \ref{H4} hold true and the jump functions $\phiA,\phiJ$ have linear growth and are Lipschitz in the sense of assumptions \ref{hip:linear-growth-phi}-\ref{hip:lipschitz-condition-phi}. Assume further that there exist constants $\varepsilon, \varepsilon_1, M>0$, with $\varepsilon<\varepsilon_1$, such that 
\begin{align}
&\text{supp}(\phiA(\cdot,z)),\,\text{supp}(\phiJ(\cdot,z)) \subseteq [\varepsilon_1,\infty),\forall z\in \text{supp}(\nu),\label{ass:no_jump_near_0_supp}\\ 
 &x+ \phi_i(x,z)>\varepsilon, \forall (x,z)\in \text{supp}(\phi_i(\cdot,z))\times \text{supp}(\nu),\, i = J,A\label{ass:no_jump_near_0_nokilling}\\
&\int_{\R}\left(| \phi_A (x, z)| + |\phi_J (x, z)|\right)\nu (dz) < M,\mbox{ uniformly on $x$.}\label{ass:no_jump_near_0_no_excess}
\end{align}

Then, the system \eqref{eq:full-model} is well-posed and, with probability one, the process $\K_t$ does not attain the origin in finite time, i.e.
\[
\mathbb{P}_{\bf z} \left( \exists t > 0 : \K_t = (0,0) \right) = 0.
\]
\end{proposition}
\begin{remark} 
It is widely documented that extreme climatic events, such as the warmer phase of 1982/83 ENSO \cite{camus1994recruitment}, can cause local extinction of the kelp populations. However, recovery is still possible, albeit slow, through recolonization of the surviving propagules \cite{martinez2003recovery}. In view of this fact, the assumption of the previous theorem can be interpreted as saying that, when the total biomass is very small, catastrophic events are not accessible to the population (see Assumption \eqref{ass:no_jump_near_0_supp}); hence, no stochastic change due to external events occurs. Similarly, if the total population exceeds a minimal threshold $\varepsilon_1$, then both increases and decreases in biomass may occur (see Assumption \eqref{ass:no_jump_near_0_nokilling}). However, such fluctuations have a limited effect when the population exceeds natural thresholds (e.g., the carrying capacity), due to constraints like maximum biomass, assumption that is supported by condition \eqref{ass:no_jump_near_0_no_excess}.
\end{remark}

\begin{proof}
    
We notice that the assumptions of the statement are a reinforcement for \ref{hip:jumps_positiveness}, and we can guarantee the well-posedness for \eqref{eq:full-model} in $\mathcal{S}$. Now, let $k\geq1$ large enough so that $\K_0 = (J_0,A_0)\in [1/k,k]^2$ and define the random time
\begin{align*}
\eta_k:= \inf\{t \ge 0: |\XA_t+\XJ_t|<  1/k\mbox{, or }|\XA_t+\XJ_t|_1 > k\},
\end{align*}
and the function $g(x,y,e) = x+y-\log(x+y)$. 

Since $\vert \K_0 \vert \in [1/k,k]\subset \R_+$, and $\K$ is c\`adl\`ag, $\eta_k$ is a stopping time  and it corresponds to the first time the trajectories of $\K$ leave $\{(x,y): x+y \in [1/k,k]\}$. From It\^o formula we know that:
\begin{align*}
\E[g(\K_{\eta_k\wedge t},E_{\eta_k\wedge t})] &= \E[g(\K_0,E_0)] + \E\left[\int_0^{\eta_k\wedge t} \mathcal{L}g(\K_s,E_s)ds\right], 
\end{align*}
where 
\begin{align*}
&\mathcal{L} g(x, y,e) \\
&\, = \bigg ( 1-\dfrac{1}{x+y}\bigg)\left\{ \rJ(e)y \ind{0\leq x+y\leq \KA} -(\mJ +\Fmaxj(e))x-(\mA +\FmaxA(e))y\right\} \\
&\, + \dfrac{1}{2} \dfrac{1}{(x+y)^2} (\sigma_A^2 x^2+ \sigma_J^2 y^2)+ \lambda\int_\R\left[ g (x+\phi_A(x, z), y+\phi_J (y, z))- g(x, y) \right]\nu (dz).
\end{align*}

Notice that, for any $(x,y,e,p)\in \mathcal{S}$ we have:
\begin{align*}
\mathcal{L} g(x, y,e) &\leq   \rJ(e)\, y\,  \ind{0\leq x+y\leq \KA}\\
&\qquad + \dfrac{\sA^2\vee \sJ^2}{2} + \lambda\int_\R \left[ g (x+\phi_A(x, z), y+\phi_J (y, z),e)- g(x, y,e)\right] \nu (dz).
\end{align*}

The first term above can be bounded by $\max_{e\in(0,1)}\rJ(e) \, \KA$. Moreover, from our assumptions on $\phiJ$ and $\phiA$, the logarithmic functions in the later integral are well defined and can be estimated in $\mathcal{S}$ as
\begin{align*}
&\int_\R \left[ g (x+\phi_A(x, z), y+\phi_J (y, z),e)- g(x, y,e) \right]\nu (dz)\\
&\qquad=\int_\R \left( \phiA(y, z) + \phiJ(x, z)-\log\Big(1+\frac{\phiA(y, z) + \phiJ(x, z)}{x+y}\Big) \right) \nu (dz)\\
&\qquad \leq \int_{\R}  \left( \phiA(y, z) + \phiJ(x, z)\right)\left(\frac{\phiA(y, z) + \phiJ(x, z) +x+y-1}{\phiA(y, z) + \phiJ(x, z) +x+y}\right) \nu (dz) =: I(x,y).
\end{align*}

Denoting $U_{\pm}(x,y) = \{z\in\R:~\text{sign}(\phiA(y, z) + \phiJ(x, z))=\pm1\}$, and $I_\pm$ the integral $I$ over the set $U_\pm$, we then have

\begin{align*}
&I_+(x,y)\leq \int_{U_{+}}  \left( \phiA(y, z) + \phiJ(x, z)\right) \nu (dz)< M,\\
&I_{-}(x,y)\leq \int_{U_{-}}  \frac{-(\phiA(y, z) + \phiJ(x, z))}{\phiA(y, z) + \phiJ(x, z) +x+y} \nu (dz)\leq \frac{M}\varepsilon.
\end{align*}

In particular, for all $(x,y,e)\in[1/k,k]^2\times(0,1)$, we obtain
\begin{align*}
\mathcal{L} g(x, y,e) &\leq   \sup_{z\in(0,1)}\rJ(z)\, \KA\, + \dfrac{\sA^2\vee \sJ^2}{2} + \lambda M\Big(1+\frac1\varepsilon\Big)=:C,
\end{align*}
which gives us
\begin{align}\label{eq:Dynkin}
\E[g(\K_{\eta_k\wedge t},E_{\eta_k\wedge t})] &\leq  \E[g(\K_0,E_0)] + C t.
\end{align}

By fixing $N \in \N$, and defining $\epsilon_N = \P(\Omega_N)$ for the measurable sets \(
\Omega_N:= \bigcap_{k \ge 1}\{ \omega \in \Omega: \eta_k(\omega) < N\},\)
we observe that the left hand side in \eqref{eq:Dynkin} with $t=N$ admits the lower bounds:
\begin{align*}
    \E[g(\K_{\eta_k\wedge N},E_{\eta_k\wedge N})] &\ge \E[\mathbf{1}_{\Omega_N}g(\K_{\eta_k }, E_{\eta_k })] \ge \operatorname*{ess\,inf}_{\omega\in\Omega}g(\XA_{\eta_k(\omega)}(\omega), \XJ_{\eta_k(\omega)}(\omega),E_{\eta_k (\omega)}(\omega))\,  \epsilon_N.
\end{align*}
The function $g$ is strictly decreasing on $(x+y) \in (0,1]$ and strictly increasing on $[1, +\infty)$, and for $k > 1$:
\begin{align*}
\inf_{(x+y) \in [1/k, k]^c\cap (0,\infty)} g(x, y,e) =\frac1k + \log(k)>\log(k),
\end{align*}
and so there exists a strictly increasing sequence $c_k$ with $c_k \to \infty$ such that $\displaystyle{\inf_{x+y \in [1/k,k]^c\cap (0,\infty)} g(x, y,e)\geq c_k}$.

Now, there are two possible cases: either $\eta_k(\omega)$ is a jump time or not. We denote by $\Omega_J^k\subset \Omega$ the set of events $\omega\in\Omega$ for which $\eta_k(\omega)$ is a jump time of the compound Poisson process, i.e.
\(\Omega_J^k = \{\omega\in\Omega:~\K_{\eta_k(\omega)^-}\neq \K_{\eta_k(\omega)^+}\}\). Then, for $\omega\in \Omega\setminus\Omega_J^k$, $\eta_k(\omega)$ is not a jump time and the process is leaving the region $(\XJ, \XA) \in [1/k, k]$ continuously. In other words, for any $\omega\in \Omega\setminus\Omega_J^k$ we certainly have $g(\XJ_{\eta_k}(\omega),\XA_{\eta_k}(\omega),E_{\eta_k }(\omega)) \ge c_k$. If, on the contrary, $\omega\in \Omega_J^k$, then we have to consider two cases where either $(\XJ, \XA)$ is leaving the box by above or below. In the first case, $\XJ_{\eta_k}(\omega) + \XA_{\eta_k}(\omega)>k$, and we still have $g(\XJ_{\eta_k}(\omega),\XA_{\eta_k}(\omega),E_{\eta_k}(\omega)) \ge c_k$, for some divergent sequence $c_k$. Finally, the case where $\omega\in\Omega_J^k$ and $\XJ_{\eta_k}(\omega)+\XA_{\eta_k}(\omega)<1/k$ is ruled out by the condition \eqref{ass:no_jump_near_0_nokilling} as soon as $k > \dfrac{1}{2\varepsilon}$. We deduce that for $k$ large enough and every $N$:
\begin{align*}
c_k \epsilon_N \le \E[g(\K_0, E_0)] + CN.
\end{align*}
Taking the limit in $k \to \infty$, necessarily we have $\epsilon_N=0$ for every $N \in \N$, and this proves the assertion. \end{proof}

\begin{example}
Let $\varepsilon_1,\varepsilon_2>0$ be given and consider $g_A,g_J$ functions bounded by 
\[|g_A(z)| \leq 1-\varepsilon_2.\]

We define the functions 
\begin{align*}
\phiJ(x,z) &= g_J(z)~(x\wedge \KA-\varepsilon_1)~\ind{x\geq\varepsilon_1},\\
\phiA(y,z) &= g_A(z)~(y\wedge \KA-\varepsilon_1)~\ind{y\geq \varepsilon_1}.
\end{align*}

Then, it is easy to check that $\phiJ,\phiA$ satisfy the hypotheses of Proposition \ref{prop:non-attainability-origin}.
\end{example}

\subsection{Existence of invariant measures}

In order to prove the existence of an invariant measure, we establish the uniform-in-time boundedness of suitable moments which, in turn, ensure tightness of the family of time-averaged laws. 

\begin{proposition}[Existence of invariant measure]\label{th:invariant}
Assume \ref{H1}, \ref{H2}, \ref{H3} and \ref{H4}. Let $\ZZ_0 \in \mathcal{S}$ be an arbitrary initial condition for the SDE \eqref{eq:full-model} and $\ZZ$ its strong solution. If there exist a constant $M>0$ such that \eqref{ass:no_jump_near_0_no_excess} holds, then
\begin{align*}
\sup_{t \in [0, \infty)} \E[\Vert \K_t\Vert_1 ] < \infty.
\end{align*}

Furthermore, under the same assumptions, if $\sup_{t \ge 0} \E (P_t) < \infty$, there exist an invariant probability measure $\mu$ for $\XX$, i.e. if $\XX_0 \sim \mu$, for all $t>0$  we have $\P_{\XX_0} (\XX_t \in A) = \mu (A)$.
\end{proposition}

\begin{proof}
Let $U: \R^2_+ \mapsto \R_+$ be given by $U(x, y,e,p)= x+y$ and let $\mathcal{L}$ be the infinitesimal generator associated with the Markov process $\ZZ$. From our assumptions, it is direct that there exist some constants $C$ and $\alpha > 0$ such that:
\begin{align*}
\mathcal{L} U (x, y, e, p)  &\le \sup_{e\in(0,1)}\rJ(e)\Psi\left(1-\frac{x+y}{\KA}\right)y - \alpha U(x,y)\\
&\qquad+\lambda\int_\R \big(\phiJ(x,z)+\phiA(y,z)\big)\nu(dz)\\
&\leq C -\alpha U (x,y),
\end{align*}
with $\alpha = \mA+\mJ +\displaystyle{\inf_{e\in(0,1)}(\FmaxA(e)+\Fmaxj(e))}$, and \(C = \lambda M+ \displaystyle{\sup_{e\in(0,1)}\rJ(e) \KA}\).

Consider next the function $V(t, x, y)= e^{\alpha t } U (x, y)$. Let $t_k= t \wedge \eta_k$, with \(\eta_k:= \inf\{t \ge 0: \XA_t+\XJ_t >k \}.\) Then, the process:
\begin{align*}
M^{(k)}_{t} := V(t_k,\XJ_{t_k}, \XA_{t_k} )-V(0,\K_0 ) - \int_{0}^{t_k} \left (e^{\alpha s} \mathcal{L} U(\XJ_{s}, \XA_s) + \alpha e^{\alpha s} U(\XJ_{s}, \XA_s)\right) ds
\end{align*}
is a martingale, and so we can use the bound above obtaining:
\begin{align*}
e^{\alpha t}\E[ U(\XJ_{t_k}, \XA_{t_k})\, \ind{t\leq \eta_k}]\leq \E[e^{\alpha t_k} U(\XJ_{t_k}, \XA_{t_k})] \le \XJ_0 + \XA_0 +C\, \E\left[\int_{0}^{t_k} e^{\alpha s} ds\right],
\end{align*}
since $\alpha>0$.

In virtue of Theorem \ref{th:well-posedness}, the right-continuity of the trajectories and Fatou's Lemma, we take $k \to \infty$:
\begin{align*}
\sup_{t\geq 0}\E[ \|\K_t\|_1]\leq \sup_{t\geq0}\left\{\|\K_0\|_1 e^{-\alpha t}+\dfrac{C}{\alpha} (1-e^{-\alpha t})\right\} < \infty,
\end{align*}
concluding the uniform-in-time control of the first moment for the population component of $\ZZ$. Together with the uniform-in-time control of the price process, and the fact that $\P(E_t\in (0,1),\forall t\geq0)=1$, we get that
\begin{align*}
\sup_{t\geq 0}\E[ \|\ZZ_t\|_1] < \infty.
\end{align*}

Denoting $B(r)\subset \R^2$ the ball of radius $r$ centered at the origin, we can use the uniform-in-time moments together with Markov inequality to get:
\begin{align*}
  t- \int_0^t\P_{\ZZ_0}(\ZZ_s \in B(r))ds&\leq \int_0^t\frac{\E_{\ZZ_0}[\|\ZZ_s\|]}{r}ds\leq \sup_{s\geq0}\E_{\ZZ_0}[\|\ZZ_s\|]\frac{t}{r}\\
  1-\frac{\sup_{s\geq0}\E_{\ZZ_0}[\|\ZZ_s\|]}{r} &\leq \frac1t \int_0^t\P_{\ZZ_0}(\ZZ_s \in B(r))ds.
\end{align*}

Consequently, for any $\epsilon>0$ we can choose $r>0$ large enough so that 
\[1-\epsilon\leq \nu_t(B(r)):= \frac1t \int_0^t\P_{\ZZ_0}(\ZZ_s \in B(r))ds,\]
i.e. $(\nu_t)_{t\geq0}$ is tight. The result follows after noticing that the process satisfies the Feller property in $\mathcal{C}_b(\R^4)$ (see, e.g., \cite[Sec. 6.7.1, Note 2 and Note 3]{applebaum2009levy}) and we can apply Krylov-Bogoliubov theorem \cite[Cor. 3.1.2]{da1996ergodicity}.
\end{proof}

\begin{remark}
Notice that, in Proposition \ref{th:invariant} we have assumed the price is controlled uniformly in time.  In simpler models, such as geometric Brownian motion, which has been used previously for modeling algae price~\cite{walsh2018financial}, this control is trivial when the return rate is low, leading to price collapse. In contrast, mean-reverting frameworks (e.g., Ornstein-Uhlenbeck process for modeling the log price of commodities) yield non-trivial behaviors and invariant distributions \cite{schwartz1997stochastic}.
\end{remark}

\subsection{Stability}

We turn to the analysis of the asymptotic behavior of the Lyapunov exponents of the total population. In particular, we aim to determine the conditions under which the species goes extinct. To do this, we study the behavior of the process $(\log(\XJ_t+\XA_t):~t\geq0)$. 

\begin{proposition}\label{prop:extinction} Assuming the hypotheses of Proposition~\ref{prop:non-attainability-origin} together with:
\begin{align}\label{hyp:stability_hyp}
\int_{\R}\!\left(|\phi_A(x,z)|^2 + |\phi_J(x,z)|^2\right)\nu(dz) &< M \quad \text{uniformly in $x$,} \nonumber\\
\frac{\sA^2}{2}\wedge \check{m}_J + \check{m}_A - \hat{r} &> 0, \\
\min\!\left\{\check{m}_J+\tfrac{\sJ^2}{2},\, \check{m}_A - \hat{r} + \tfrac{\sA^2}{2}\right\} &> \tfrac{M\lambda}{\varepsilon_1}, \nonumber
\end{align}
where
\[
\hat{r} := \sup_{e\in(0,1)} \rJ(e), 
\quad \check{m}_A := \min_{e\in(0,1)} \{\mA + \FmaxA(e)\}, 
\quad \check{m}_J := \min_{e\in(0,1)} \{\mJ + \Fmaxj(e)\}.
\]
Then there exists a positive constant $\eta$ such that for any initial condition 
\begin{align*}
\limsup_{t\rightarrow+\infty}\frac{\log(\XJ_t+\XA_t)}{t}&<-\eta,\;\;\P\text{-a.s.}
\end{align*}
In particular the population get asymptotically  extinct exponentially fast. 
\end{proposition}

\begin{proof}
The assumptions in consideration guarantee that there exists a unique solution of the system \eqref{eq:full-model} that satisfies $\K_t\in \R_+^2$ for all times. Denote $U_t = \XJ_t+\XA_t$. With the help of It\^o formula for L\'evy processes and assumption \ref{H1}, we have
\begin{align*}
&\log(U_t)-\log(U_0) \leq \sA\int_0^t\frac{\XA_s}{U_s}d\WA_s + \sJ\int_0^t\frac{\XJ_s}{U_s}d\WJ_s\\
&\qquad + \int_0^t \int_\R \log\left(\frac{U_{s^-} + \phiA(\XA_s,z)+\phiJ(\XJ_s,z)}{U_{s^-}}\right)\Np(dz\otimes ds)\\
&\qquad + \int_0^t \frac{\hat{r}~ \Psi\left(1-\frac{U_s}{\KA}\right)\XA_s U_s - \check{m}_J \XJ_s U_s - \check{m}_A \XA_s U_s - \frac{\sA^2}{2}\XA_s^2- \frac{\sJ^2}{2}\XJ_s^2}{U_s^2}ds.
\end{align*}

If we compensate the Poisson random measure we can write
\begin{align*}
&\frac{\log(U_t)}t-\frac{\log(U_0)}t\leq \frac{M_t^{(A)}}t + \frac{M_t^{(J)}}t + \frac{\tilde{N}_t}t\\
&\qquad +\frac\lambda t\int_0^t \int_\R \log\left(\frac{U_{s} + \phiA(\XA_s,z)+\phiJ(\XJ_s,z)}{U_{s}}\right)\nu(dz)ds\\
&\qquad \frac1{t}\int_0^t \frac{\big(\hat{r}-\frac{\sA^2}{2}-\check{m}_A\big)~ \XA_s^2 + \big(\hat{r}-\check{m}_A-\check{m}_J\big)~ \XA_s\XJ_s - \big(\check{m}_J + \frac{\sJ^2}{2}\big)\XJ_s^2}{U_s^2}ds,
\end{align*}
where the processes $M_t^{(A)} = {\sA}\int_0^t\frac{\XA_s}{U_s}d\WA_s$, $M_t^{(J)} = {\sJ}\int_0^t\frac{\XJ_s}{U_s}d\WJ_s$, and 
\[\tilde{N}_t:=\int_0^t \int_\R \log\left(\frac{U_{s^-} + \phiA(\XA_s,z)+\phiJ(\XJ_s,z)}{U_{s^-}}\right)\big(\Np(dz\otimes ds)-\lambda\nu(dz)ds\big)\]
define martingales.

We notice that $\limsup_{t\rightarrow\infty}\frac{\langle M^{(A)}\rangle_t}t \leq \sA^2<\infty$, and similarly for ${\langle M^{(J)}\rangle_t}$. Thus, from the Law of Large Numbers for martingales (see, e.g., \cite[Th. 1.3.4]{mao2007stochastic}),
\[\limsup_{t\rightarrow\infty}\frac{M_t^{(A)}}t = \limsup_{t\rightarrow\infty}\frac{M_t^{(J)}}t = 0.\]

As for the  martingale $\tilde{N}$, we use the inequality $\log^2(\frac{x}y)\leq\left(\frac{x-y}{y\wedge x}\right)^2$, valid for all $x, y>0$. We get: 
\begin{align*}
\int_0^t \frac{d\langle \tilde{N}\rangle_s}{(1+s)^2}&= \lambda\int_0^t\frac1{(1+s)^2} \int_\R \log^2\left(\frac{U_{s} + \phiA(\XA_s,z)+\phiJ(\XJ_s,z)}{U_{s}}\right)\nu(dz)ds\\
&\leq \lambda\int_0^t\frac1{(1+s)^2}  \int_{\R}\left(\frac{\phiA(\XA_s,z)+\phiJ(\XJ_s,z)}{U_{s}\wedge \big(\phiA(\XA_s,z)+\phiJ(\XJ_s,z)+ U_{s}\big)}\right)^2\nu(dz)ds.
\end{align*}

From Assumption \ref{ass:no_jump_near_0_supp}, we notice that $\int_0^t \frac{d\langle \tilde{N}\rangle_s}{(1+s)^2}$ is non-decreasing and non-zero if and only if there exist some $I\subset[0,t]$ such that $\XJ_s>\varepsilon_1$ or $\XA_s>\varepsilon_1,$ for all $s\in I$. Moreover, from Assumption \ref{ass:no_jump_near_0_nokilling}, we have, for all $z\in \R$:
\begin{align*}
U_{s}\wedge \big(\phiA(\XA_s,z)+\phiJ(\XJ_s,z)+ U_{s}\big)&= U_s - \big(\phiA(\XA_s,z)+\phiJ(\XJ_s,z)\big)_{-}\\
&\geq \XA_s -\big(\phiA(\XA_s,z)\big)_{-} +\XJ_s -\big(\phiJ(\XJ_s,z)\big)_{-}> 2\varepsilon.
\end{align*}

Hence, assuming that $\int_{\R}\phiA^2(x,z)+\phiJ^2(y,z)\nu(dz)\leq M$, uniformly in $x,y\in(0,\infty)^2$, we found
\begin{align*}
0\leq \int_0^t\frac{d\langle \tilde{N}\rangle_s}{(1+s)^2}&< \frac{M\lambda }{2\varepsilon^2} \int_0^t\frac1{(1+s)^2}\longrightarrow \frac{M\lambda }{2\varepsilon^2},\text{ as }t\to\infty.
\end{align*}

From this we conclude that $\P\left(\lim_{t\rightarrow+\infty}\int_0^t\frac{d\langle \tilde{N}\rangle_s}{(1+s)^2}<+\infty\right) = 1$ and, applying again a suitable Law of Large Numbers to the square-integrable local martingale $\tilde{N}$ and the increasing process $(t;\, t\geq0)$ (see e.g. \cite[Th. 1]{liptser1980strong}), we get
\begin{align*}
&\limsup_{t\rightarrow\infty}\frac{\log(U_t)}t \leq  \limsup_{t \rightarrow\infty} \Bigg(\frac{\lambda}{t}\int_0^t \int_\R \frac{(\phiA(\XA_s,z)+\phiJ(\XJ_s,z))(\XA_s+\XJ_s)}{U_{s}^2}\nu(dz)ds\\
&\qquad + \frac1{t}\int_0^t \frac{\big(\hat{r}-\frac{\sA^2}{2}-\check{m}_A\big)~ \XA_s^2 + \big(\hat{r}-\check{m}_A-\check{m}_J\big)~ \XA_s\XJ_s - \big(\check{m}_J + \frac{\sJ^2}{2}\big)\XJ_s^2}{U_s^2}ds\Bigg).
\end{align*}

Using Assumptions \ref{ass:no_jump_near_0_supp}-\ref{ass:no_jump_near_0_no_excess} we can bound the above limit as
\begin{align*}
&\limsup_{t\rightarrow\infty}\frac{\log(U_t)}t\leq  \limsup_{t \rightarrow \infty} \frac1{t}\int_0^t \frac{ M\lambda\XA_s\ind{\varepsilon_1<\XA_s}+M\lambda\XJ_s\ind{\varepsilon_1<\XJ_s}}{U_s^2}ds\\
&\quad +  \limsup_{t \rightarrow \infty} \frac1{t}\int_0^t \frac{\big(\hat{r}-\frac{\sA^2}{2}-\check{m}_A\big)~ \XA_s^2 + \big(\hat{r}-\check{m}_A-\check{m}_J\big)~ \XA_s\XJ_s - \big(\check{m}_J + \frac{\sJ^2}{2}\big)\XJ_s^2 }{U_s^2}ds.
\end{align*}

Finally, we observe that the integrand is negative and bounded provided that \eqref{hyp:stability_hyp} are satisfied.
More precisely, we obtain 
$\limsup_{t\rightarrow\infty}\frac{\log(U_t)}t<-\eta,$ with \newline
$\eta: = \min\left\{ \frac{\sA^2}{2} + \check{m}_A - \hat{r}- \frac{M\lambda}{\varepsilon_1}, \check{m}_J+\frac{\sJ^2}{2} - \frac{M\lambda}{\varepsilon_1}\right\}$.

The last claim follows from the usual arguments.  Fix $\delta \in (0, 1$) and a strictly increasing sequence $t_j \rightarrow \infty$ and define $A_j:=\{ \omega: \vert \K_{t_j}(\omega)\vert > \exp\{-t_j \eta(1-\delta) \}\}.$
Then, 
\begin{align*}
A_j & = \left\{\omega: \dfrac{\log(\vert \K_{t_j}(\omega)\vert)}{t_j}  > -\eta(1-\delta)\right\},
\end{align*}
and the previous theorem implies that $\P ([A_j \text{ i.o.}])=0$. We deduce that there exists a (random but finite) index $\overline j= \overline j(\omega)$ such that for any $j \ge \overline j$ we have $\vert \K_{t_j}\vert  \le \exp(-t_j (\eta (1-\delta)))$, and thus the population decays exponentially fast for exponents arbitrarily close to $\eta$.

\end{proof}

\section{Time-discrete approximation and numerical experiments}\label{sec:numerical_scheme}

In previous sections we have proved analytic properties of our target system \eqref{eq:full-model}. In particular, we have proved that under our running assumptions, there is no absorption in finite time for neither the abundance process $\K$ nor the compliance level process $E$. However, and as usual in stochastic models with non-attainable boundaries, it is still possible that the process gets arbitrarily close to the extinction boundary. This poses the challenge of finding theoretically reliable and numerically feasible   approximation techniques suitable for simulation. Standard numerical methods may fail to accurately capture the behavior of the exact process, particularly in biologically meaningful regimes. 

Several works address positivity-preserving numerical schemes for SDEs with non-globally Lipschitz coefficients~\cite{bossy2021weak, cai2022positivity, cai2024advanced}. One way to overcome the lack of global Lipschitz continuity and ensure convergence is through truncated schemes, as proposed in \cite{cai2022positivity} in the context of epidemiological models. To preserve the positivity of the numerical solution when the coefficients grow polynomially, exponential-type schemes can be used, such as those developed in \cite{bossy2021weak}.

In this paper we introduce a time-discretization process which yields an admissible approximation within the context of biomass and proportions, and we prove its convergence to the exact solution of \eqref{eq:full-model}. 

\subsection{Overview of the construction of the numerical scheme}

We consider a simulation horizon $T>0$, $\dt = T/N $ with $N\in \N$ to be fixed later, a temporal grid $t_k=k\dt$ and $\eta(t):=\sup\{t_k:t_k< t\}$, and $\bZZ_{0}:=\ZZ_0=(\XA_0,\XJ_0,E_0,P_0)$, the initial condition of \eqref{eq:full-model}. Since the process $P$ is exogenous, we will assume that there is a numerical scheme $\bar{P}$, with bounded moments, such that
\begin{equation*}\label{eq:convergence-scheme-for-P}
\begin{aligned}
\sup_{0\leq t\leq T}\E\left[(P_t - \bar P_t)^{2p} \right] &\leq  C\dt^p.
\end{aligned}
\end{equation*}

With the above ingredients, we construct -iteratively- an approximating process $\bZZ_{t}:= (\scXJ, \scXA_t, \bE_t, \bar P_t)$: if we have a value for $\bZZ_{\eta(t)}$, we compute the approximation in $(\eta(t),\eta(t)+\dt)$. The details are given below.

\paragraph{For the compliance approximated-process $E$} We define the approximation as a truncation of the Euler-Maruyama approximation of $E$ denoted as $\hE$, ensuring that it remains within the interval $[0,1]$, or more precisely within $[\dt,1-\dt]$.

\begin{equation}\label{eq:def-E-barra}
    \begin{aligned}
        \hE_t &= \bE_{\eta(t)} +  \int_{\eta(t)}^t\left(\beta_1(\bZZ_{\eta(t)})(1-\bE_{\eta(t)})  -\beta_0(\bZZ_{\eta(t)})\bE_{\eta(t)}\right)ds +\int_{\eta(t)}^t\sE\sqrt{\bE_{\eta(t)}(1-\bE_{\eta(t)})}dW^{(E)}_s,\\
        \bE_t &= \cutDt(\hE_t ).
    \end{aligned}
\end{equation}

%
%
%

\paragraph{For the biomass process $(\K_t)$.}
We adopt a splitting approximation: first, we simulate the continuous (diffusive) part using an exponential-type numerical approximation to preserve positivity between jumps; then, within each subinterval of the time partition, we approximate the jump component by retaining only the first jump, capturing the dominant discrete effect while maintaining numerical tractability. 

Focusing in the time interval $[\eta(t), t]$, and given $\bXA_{\eta(t)},\bXJ_{\eta(t)}$ and $\bE_{\eta(t)}$, we freeze the drift coefficient of the equation for $\XA$ at time $\eta(t)$, and neglect the jumps. Doing so, we obtain the following linear SDE: 
\begin{equation*}\label{eq:def-YA-inv}
  \begin{aligned}
    \hXA_t : =&\bXA_{\eta(t)} 
    + \int_{\eta(t)}^t\left( \rhoA\bXJ_{\eta(t)}-\kappaA(\bE_{\eta(t)})\bXA_{\eta(t)}\right)ds + \int_{\eta(t)}^t \sA \hXA_sd\WA_s,
\end{aligned}  
\end{equation*}
which can be solved explicitly (See \cite[Problem 6.15, Ch. V]{karatzas98}):
\begin{equation*}\label{eq:def-YA-hat}
  \begin{aligned}
    \hXA_t =& \bXA_{\eta(t)} \exp\left(-\frac{\sA^2}{2}\ddt + \sA\deltaWA{\eta(t)}{t}\right) \\
&\quad+     \left(\rhoA\bXJ_{\eta(t)}-\kappaA(\bE_{\eta(t)})\bXA_{\eta(t)}  \right)\int_{\eta(t)}^{t}\exp\left(-\frac{\sA^2}{2}(t-s) + \sA \deltaWA{s}{t}\right) ds.
\end{aligned}  
\end{equation*}

This gives us a time-approximation of the continuous part of the process $\XA$ in the time interval $[\eta(t),t]$.

Next, we define the approximation $  \bXA_t $ by incorporating the first jump occurring within the considered time interval. If $\dt$ is small enough, the probability of having two jumps in $(\eta(t),t)$ is negligible. Then we introduce 
\begin{equation*}\label{eq:def-XA-hat}
\begin{aligned}
  \bXA_t &= 
  \hXA_t  +  \int_{\eta(t)}^{t^+}\int_{-\infty}^\infty \phiA(\bXA_{\eta(t)},z)\ind{\{s^{-}\leq \tau_{\eta(t)}^1\}}\Np(dz\otimes ds)\\
 & = \hXA_t  +   \ind{\{\eta (t) \le \tau^{1}_{\eta(t)} < t\}} \phiA \left(\bXA_{\eta(t)}, \xi_1\right), 
\end{aligned}
\end{equation*}
where $\tau^{1}_{\eta(t)}$ is the first jump time of the Poisson process after $\eta(t)$ and $\{\xi_i:i\geq1\}$ is an i.i.d. sequence with common law $\nu$. 

Notice that $\hXA$ includes a time dependent integral which cannot be implemented computationally. Thus, we define the  numerical scheme for $\XA$, which provides a reasonable and easily implementable approximation:

\begin{equation}\label{eq:explicit-XA-tilde}
\begin{aligned}
\scXA_{t} & = 
  \scXA_{\eta(t)}\exp\left(-\frac{\sA^2}{2}\ddt + \sA\deltaWA{\eta(t)}{t}\right) \\
&\quad+     \left(\rhoA\scXJ_{\eta(t)}-\kappaA(\bE_{\eta(t)})\scXA_{\eta(t)}  \right)\ddt  \exp\left(-\frac{\sA^2}{2}\ddt + \sA\deltaWA{\eta(t)}{t}\right) \\
&\quad+   \ind{\{\eta(t) \le \tau^{1}_{\eta(t)} < t\}} \phiA \left(\scXA_{\eta(t)}, \xi_1\right).
\end{aligned}
\end{equation}


Analogously, we introduce $\scXJ$ the scheme for the juvenile kelp biomass $\XJ$: 

\begin{equation}\label{eq:explicit-XJ-tilde}
\begin{aligned}
&\scXJ_{t}  = \scXJ_{\eta(t)} \exp\left(-\frac{\sJ^2}{2}\ddt  + \sJ\deltaWJ{\eta(t)}{t}\right) \\
&\quad+    \left(\rhoJ(\bZZ_{\eta(t)})\scXA_{\eta(t)} -\kappaJ(\bE_{\eta(t)})\scXJ_{\eta(t)} \right)\ddt \exp\left(-\frac{\sJ^2}{2}\ddt  + \sJ\deltaWJ{\eta(t)}{t}\right) \\
&\quad+   \ind{\{\eta (t) \le \tau^{1}_{\eta(t)} < t\}} \phiJ \left(\scXJ_{\eta(t)}, \xi_1\right).
\end{aligned}
\end{equation}
Notice that $(\scXJ_{\eta(t)} ,\scXA_{\eta(t)} )$ are both non-negative, then $(\scXJ_{t},\scXA_{t})$ are both non-negative for $\dt$ small enough. Indeed, if $\Delta t$ satisfies
\begin{equation}\label{eq:condition_positivity_Dt}
1-\|\kappaJ\|_\infty\vee\|\kappaA\|_\infty\Delta t\geq 0,
\end{equation}
and Assumption \ref{hip:jumps_positiveness} holds, then both $\scXJ$ and $\scXA$ are non-negative. 

\subsection{Convergence of the numerical scheme}

In order to prove the strong convergence of the numerical scheme \eqref{eq:def-E-barra}-\eqref{eq:explicit-XA-tilde}-\eqref{eq:explicit-XJ-tilde} to the exact process $(\K, E)$, we adapt the truncation technique of Cai et al. \cite{cai2022positivity} to handle non-Lipschitz terms (including the superlinear growth drifts and the square-root diffusion term) and preserve positivity and boundedness. 

For fixed $\varepsilon > 0$ and $\delta := \delta(\varepsilon) > 0$ (to be specified later), we define $\ZZ^\delta = (\XA^\delta, \XJ^\delta, E^\delta, P)$ as the solution of a \textit{truncated system}:
\begin{equation}\label{eq:full-model-approximation}
\left\{\begin{aligned}
    \XJ^\delta_t
    =\XJ^\delta_0&+\int_0^t\left(\rrJ(\XJ^\delta_s,\XA^\delta_s,E^\delta_s)\Psi_{1-1/\delta}(\XA_s^\delta) - \kappaJ(E^\delta_s)\Psi_{1-1/\delta}(\XJ_s^\delta)\right) ds\\
    &\quad + \int_{0}^t\sJ \XJ^\delta_s dW^{(J)}_s + \int_{0}^t \int_{\R\setminus\{0\}} \phiJ (\XJ^\delta_{s^-}, z) \Np(dz\otimes ds)\\
    \XA^\delta_t= \XA^\delta_0 &+\int_0^t \left(\PjA\XJ^\delta_s-\kappaA(E^\delta_s)\Psi_{1-1/\delta}(\XA^\delta_s)\right)ds + \int_{0}^t\sA \XA^\delta_s dW^{(A)}_s \\
    &\quad+ \int_{0}^t \int_{\R\setminus\{0\}} \phiA (\XA^\delta_{s^-}, z) \Np(dz\otimes ds) \\
    E^\delta_t =E^\delta_0 &+ \int_0^t\left(\beta_1(\ZZ^\delta_s)(1-\Psi_{\delta}(E^\delta_s))  -\beta_0(\ZZ^\delta_s)\Psi_{\delta}(E^\delta_s)\right)ds\\
    &\quad +\int_0^t\sE\sqrt{\Psi_{\delta}(E^\delta_s)(1-\Psi_{\delta}(E^\delta_s)))}dW^{(E)}_s\\
    P_t = P_0 &+ \int_0^t \mu(P_s)ds + \int_0^t \sP(P_s)d\WP_s,
\end{aligned} 
\right.
\end{equation}
where, with some abuse of notation, $\Psi_{\varepsilon}(x) = \Psi_{0\vee\varepsilon}(x)$.

This truncated process $\ZZ^\delta$ coincides with the original dynamics of $\ZZ$ when $(\K, E)$ lies within a safe domain, i.e. ($\XJ,\XA,E)\in[\delta, +\infty)^2\times[\delta, 1-\delta]$. While keeping the original dynamics in the safe domain, this process modifies the coefficients outside this domain to ensure global Lipschitz continuity and preserving the biological constraints (e.g., non-negativity of populations and proportions).\\

As a companion to $\ZZ^\delta$, we also consider $(\bXA^\delta,\bXJ^\delta,\bE^\delta)$, the numerical approximation for $(\XA^\delta,\XJ^\delta,E^\delta)$. The processes  $(\bXA^\delta,\bXJ^\delta)$ and $(\hXA^\delta,\hXJ^\delta)$ are defined in complete analogy with  $(\bXA,\bXJ)$ and $(\hXA,\hXJ)$ above. Meanwhile,
\begin{align*}
  \bE^\delta_t &=E^\delta_0 
  + \int_0^t\left(\beta_1(\bZZ^{\delta}_{\eta(s)})(1-\Psi_{\delta}(\bE^\delta_{\eta(s)}))  -\beta_0(\bZZ^{\delta*}_{\eta(s)})\Psi_{\delta}(\bE^\delta_{\eta(s)}))\right)ds\\
  &\quad +\int_0^t\sE\sqrt{\Psi_{\delta}(\bE^\delta_{\eta(s)})(1-\Psi_{\delta}(\bE^\delta_{\eta(s)}))}dW^{(E)}_s.
\end{align*}

The core idea leverages the truncated process $\ZZ^\delta$ as an intermediate approximation to quantify the strong convergence error. Since the jump part is regular, the dominant error contribution comes from the continuous dynamics. 

To analyze the convergence of the implemented numerical scheme $\bZZ$ towards the exact process $\ZZ$, we first compare it with a reference method whose simpler structure allows for a more tractable analysis. This intermediate step is necessary to assess the discrepancy between both numerical methods before relating the implemented scheme to the exact solution.

\begin{proposition}\label{lem:control-of-L2-norm-hats-and-bars} Assume \ref{H1},\ref{H2}, \ref{H3}, \ref{H4} and that there exists a suitable approximation of $P$ satisfying \eqref{eq:scheme-for-P}. Then, for any $\dt>0$ satisfying \eqref{eq:condition_positivity_Dt}, we have
\begin{equation*}\label{eq:L2-norm-hats-and-bars}
\begin{aligned}
\sup_{0\leq t\leq T}\E\left[( \bXA_{t} - \scXA_{t} )^2 + ( \bXJ_{t} - \scXJ_{t} )^2  \right] \leq \Cadd\dt^2.
\end{aligned}
\end{equation*}
\end{proposition}

With the framework described above, the convergence error can be analyzed by focusing on $\XJ - \bXJ$, $\XA - \bXA$ and $E - \bE$. This analysis becomes simpler when both the true process and its approximation are confined within admissible domains, and when the diffusion coefficients are globally Lipschitz. This is the framework for the truncated process. In this case, the error corresponds to that of a standard Euler-Maruyama scheme as described in Proposition \ref{prop:convergence-delta-process} below.

\begin{proposition}\label{prop:convergence-delta-process}
Assume \ref{H1},\ref{H2}, \ref{H3}, \ref{H4} and that there exists a suitable approximation of $P$ satisfying \eqref{eq:scheme-for-P}. For any $\delta>0$, there exists a constant $C_\delta$ (possibly exploding when $\delta\to0$) such that for any $\dt>0$ satisfying \eqref{eq:condition_positivity_Dt}, we have

\begin{equation}\label{eq:approximation-error-for-delta-scheme}
\E\left[\sup_{0\leq s \leq T}\left\{ \left|\XA_s^\delta-\bXA_s^{\delta}\right|+\left|\XJ_s^\delta-\bXJ_s^{\delta}\right|+\left|E_s^\delta-\bE_s^\delta\right|\right\}\right] \leq C_\delta \sqrt{\dt},
\end{equation}
\begin{equation}\label{eq:approximation-error-for-delta-scheme-L2}
\sup_{0\leq s \leq T}\E\left[ \left(\XA_s^\delta-\bXA_s^{\delta}\right)^2+\left(\XJ_s^\delta-\bXJ_s^{\delta}\right)^2+\left(E_s^\delta-\bE_s^\delta\right)^2\right] \leq C_\delta \dt.
\end{equation}

\end{proposition}



The proofs of the auxiliary Propositions \ref{lem:control-of-L2-norm-hats-and-bars} and \ref{prop:convergence-delta-process} are postponed to Appendix \ref{sec:delayed-proofs-ns}.

Now, we are in position to prove the main result of this section: Using the a priori estimates from Propositions \ref{lem:control-of-L2-norm-hats-and-bars} and \ref{prop:convergence-delta-process}, we prove strong convergence by pivoting between the original approximation and its regularized version.

\begin{theorem} Assume \ref{H1},\ref{H2}, \ref{H3}, \ref{H4} and that there exists a suitable approximation of $P$ satisfying \eqref{eq:scheme-for-P}.Then, for all $T>0$,
\begin{equation*}\label{eq:estimate-strong-error}
\begin{aligned}
 \lim_{\dt\to0}\sup_{0\leq t\leq T}\left\{\E\left[\left(\XJ_t-\scXJ_t \right)^2\right] + \E\left[\left(\XA_t-\scXA_t \right)^2\right] + \E\left[\left(E_t-\bE_t \right)^2\right] \right\} =0.
\end{aligned}
\end{equation*}

\end{theorem}

\begin{proof}

Following Cai et al. \cite{cai2022positivity}, let us consider $\varepsilon>0$,
\begin{equation*}\label{eq:def_omega_delta_1}
\begin{aligned}
\Omega^\delta_1&:= \left\{\omega\in\Omega: 2\delta<\inf_{0\leq t\leq T}E_t(\omega) \leq \sup_{0\leq t\leq T}E_t(\omega) <1-2\delta\right\}.
\end{aligned}
\end{equation*}
Since $\bar\beta_0$,  $\bar\beta_1$ together with \ref{H3} keep the process $E$ away from the borders, from Proposition \ref{prop:weak_existence} we can choose $\delta_1$ small enough that for all $\delta<\delta_1$ we have
$$
\P\left((\Omega^\delta_1)^c\right)\leq \frac{\varepsilon}{3}.
$$

On the other hand, from Corollary \ref{cor:sup-norm_moments}, and Markov's inequality, there exists $\delta_2$  such that
the event 
$$\Omega^\delta_2:=\left\{\omega\in\Omega: \sup_{0\leq t\leq T}\| \ZZ \| < \frac{1}{\delta}-\delta\right\}$$ 
satisfies, for all $\delta<\delta_2$,
$$
\P\left((\Omega^\delta_2)^c\right)\leq \frac{\varepsilon}{3}.
$$
Next, we consider the event
\begin{equation*}
\label{eq:def_omega_delta_3}
\begin{aligned}
\Omega^\delta_3&:= \left\{\omega\in\Omega: \sup_{0\leq s \leq T}\left\{ \left|\XA_s^\delta-\bXA_s^{\delta}\right|+\left|\XJ_s^\delta-\bXJ_s^{\delta}\right|+\left|E_s^\delta-\bE_s^\delta\right|\right\} <\delta\right\}.
\end{aligned}
\end{equation*}
From Proposition \ref{prop:convergence-delta-process} and Markov's inequality, it follows that
$
\P\left(\left(\Omega^\delta_3\right)^c\right)\leq \frac{C_\delta\sqrt{\dt}}{\delta}.
$
We fix $\dt_0>0$ small enough such that for all $\dt\leq \dt_0$ 
$$
\frac{C_\delta\sqrt{\dt}}{\delta}\leq \frac{\varepsilon}{3}.
$$

Based on the earlier estimations, it is clear that
\begin{equation*}\label{eq:cotrol-prob-Omega-delta}
\begin{aligned}
\P\left((\Omega^\delta)^c \right)\leq \varepsilon,
\end{aligned}
\end{equation*}
with $\Omega^\delta := \Omega^\delta_1\cap \Omega^\delta_2\cap \Omega^\delta_3.$

Furthermore, for any $\omega\in\Omega^\delta$, $(\XA,\XJ,E)(\omega)$ and $(\XA^\delta,\XJ^\delta,E^\delta)(\omega)$ coincide, because they are two solutions for the same equation (provided that $\delta<1$) which has a unique strong solution (see Proposition \ref{th:well-posedness}). Therefore, for  $\omega\in\Omega^\delta$ and all $t\in[0,T]$:
\begin{align*}
 |\XJ_t-\bXJ^\delta_t|+|\XA_t-\bXA^\delta_t|+|E_t-\bE^\delta_t| &=   |\XJ^\delta_t-\bXJ^\delta_t|+|\XA^\delta_t-\bXA^\delta_t|+|E^\delta_t-\bE^\delta_t|    \\
  & \leq \sup_{0\leq t\leq T}\left\{  |\XJ^\delta_t-\bXJ^\delta_t|+|\XA^\delta_t-\bXA^\delta_t|+|E^\delta_t-\bE^\delta_t| \right\}  <\delta, 
\end{align*}
and 
$$E_t(\omega)\in(2\delta,1-2\delta),\;\; \XJ_t<\frac{1}{\delta}-\delta,\;\;\XA_t<\frac{1}{\delta}-\delta.$$ 
Hence,  for all  $\omega\in\Omega^\delta$
$$
\forall\,0\leq t\leq T:\,\bE^\delta_t(\omega)\in\left(\delta,1-\delta\right),\;\; \bXJ^\delta_t(\omega)<\frac{1}{\delta},\;\;\bXA^\delta_t(\omega)<\frac{1}{\delta}.
$$
From here, we conclude that  $(\bXA,\bXJ,\bE)$ and $(\bXA^\delta,\bXJ^\delta,\bE^\delta)$ also coincide for all $\omega \in\Omega^\delta$ as long as $\dt\leq (\delta/C_\delta)^2 \wedge \dt_0$. Indeed, let us assume that for such $\omega$, 
$$
\left(\bXA_{\eta(t)}(\omega),\bXJ_{\eta(t)}(\omega),\bE_{\eta(t)}(\omega)\right)
=
(\bXA^\delta_{\eta(t)}(\omega),\bXJ^\delta_{\eta(t)}(\omega),\bE^\delta_{\eta(t)}(\omega)).
$$

Then, for all $t\in(0,T)$ and all $s\in(\eta(t),\eta(t)+\dt)$, the numerical schemes will produce the same output for  $\bE$ and $\bE^\delta$, because the coefficients coincide in $(\delta,1-\delta)$. Similarly, for such $(s,\omega)$, the coefficients for $\bXJ$ and $\bXJ^\delta$ (or $\bXA$ and $\bXA^\delta$) coincide since the processes live in $(0,1/\delta)$. Therefore, for all $t\in(0,T)$ and all $s\in(\eta(t),\eta(t)+\dt)$ and $\omega\in\Omega^\delta$:
$$
\left(\bXA_{s}(\omega),\bXJ_{s}(\omega),\bE_{s}(\omega)\right)
=
(\bXA^\delta_{s}(\omega),\bXJ^\delta_{s}(\omega),\bE^\delta_{s}(\omega)).
$$


Notice that, adding $\hXA,\hXJ$ as pivots and using \eqref{eq:scheme-for-P}, 
\begin{align*}
	\E\left[\|\ZZ_t - \bar{\ZZ}_t\|^2\right]
	 & \leq 2\E\left[\left\{\left(\XJ_t-\bXJ_t \right)^2 + \left(\XA_t-\bXA_t \right)^2+\left(E_t-\bE_t \right)^2\right\}\ind{\Omega^\delta}\right]	\\
	  & \quad+ 2\E\left[\left\{\left(\XJ_t-\bXJ_t \right)^2 + \left(\XA_t-\bXA_t \right)^2+\left(E_t-\bE_t \right)^2\right\}\ind{(\Omega^\delta)^c}\right]	\\	
	 &\quad+2\E\left[\left(\bXJ_t-\scXJ_t \right)^2 + \left(\bXA_t-\scXA_t \right)^2\right] + C\Delta t.
\end{align*}
From observations above, we have for the first term in the above inequality:
\begin{equation*}\label{eq:error_decomposition}
\begin{aligned}
    \E&\left[\left\{\left(\XJ_t-\bXJ_t \right)^2 + \left(\XA_t-\bXA_t \right)^2+\left(E_t-\bE_t \right)^2\right\}\ind{\Omega^\delta}\right] \\
    &\quad= \E\left[\left\{\left(\XJ_t^\delta-\bXJ_t^\delta \right)^2 + \left(\XA_t^\delta-\bXA_t^\delta \right)^2+\left(E_t^\delta-\bE_t^\delta \right)^2\right\}\ind{\Omega^\delta}\right].
\end{aligned}
\end{equation*}

While for the second term, applying Cauchy-Schwartz inequality, we get:
\begin{align*}
    \E&\left[\left\{\left(\XJ_t-\bXJ_t \right)^2 + \left(\XA_t-\bXA_t \right)^2+\left(E_t-\bE_t \right)^2\right\}\ind{(\Omega^\delta)^c}\right] \\
    &\quad \leq C\sqrt{\left(\sup_{0\leq t\leq T}\E[\|\ZZ_t\|^4] + \sup_{0\leq t\leq T}\E[\|\hat{\ZZ}_t\|^4]\right)\P\left(\left(\Omega^\delta\right)^c\right) }. 
\end{align*}

Thus, from Propositions \ref{lem:control-of-L2-norm-hats-and-bars}, \ref{prop:convergence-delta-process} and \ref{prop:a-priori-estimates-NS} we conclude
\begin{align*}
	\E\left[\|\ZZ_t - \bar{\ZZ}_t\|^2\right]
\leq  C_\delta\dt + C\varepsilon^{1/2} + \Cadd\dt^2 + \Cadd\dt,
\end{align*}
taking supremum over $t\in[0,T]$ and $\dt\to0$, we obtain the result since $\varepsilon$ is arbitrary.
\end{proof}

\subsection{Numerical experiments}

In this section, we present the results of the simulation of model \eqref{eq:full-model}, which was produced using the numerical scheme described in Section \ref{sec:numerical_scheme}. We focus on the evolution of marginal densities, particularly emphasizing the long-term dynamics as indicators of the system's stable states. This aspect is especially relevant for policy evaluation in our ecological context.

Our simulations were calibrated using empirical observations, literature sources, and survey data. Biological parameters were estimated from six years of age-structured kelp abundance data collected at nine temporal sampling points in a no-take zone in the Atacama region \cite{vega2014monitoring}. Socioeconomic parameters were derived from kelp fishers' surveys conducted in the Antofagasta and Atacama regions \cite{avila2025exploring}, which provided data on compliance thresholds and initial behaviors. Price data were informed by historical transaction records and modeled using geometric Brownian motion, with adjustments made based on the average inflation rate of 4\% annually over the past 30 years (1991-2022) \cite{IPCBC}. All model parameters were set on a yearly timescale. In \cite{avila2025}, we provide a detailed description of the origins of the data and the calibration procedure for both the ecological and market pricing parameters in our model. Additionally, we have included further numerical experiments and analyses from a socio-ecological perspective (see also our repository \cite{git-repo}).

In what follows, we fix the random seed for consistency and use $ M=15000$ trajectories per scenario to obtain reliable estimates. The simulations span 30 years, with the first 5 years without harvesting to allow the system to transit from a random uniform initial condition, towards an ecologically meaningful state.

\begin{figure}[h!]
    \centering
    \begin{subfigure}[b]{0.45\textwidth}
     \includegraphics[width=\textwidth]{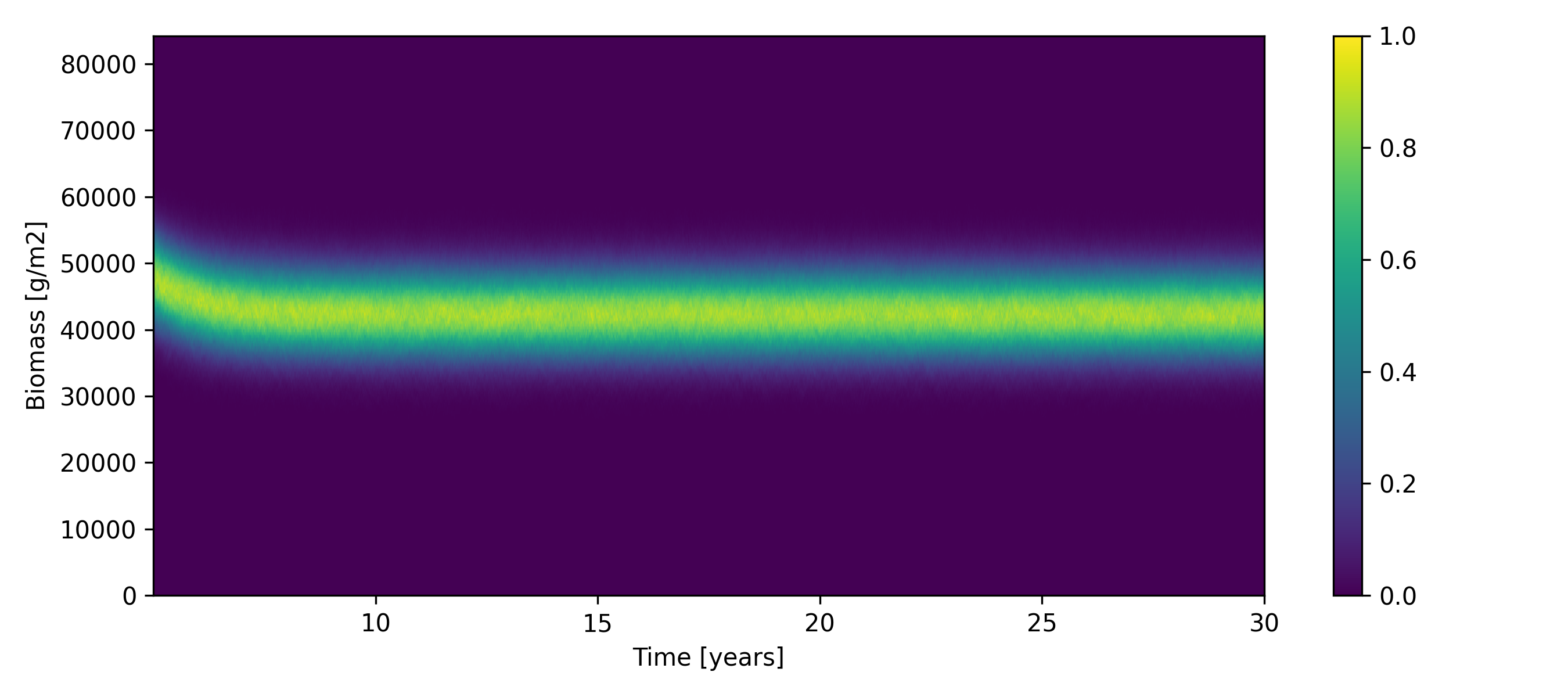}   
     \caption{Harvesting with full compliance.}
    \label{fig:full-compliance}
    \end{subfigure} 
    \begin{subfigure}[b]{0.45\textwidth}
     \includegraphics[width=\textwidth]{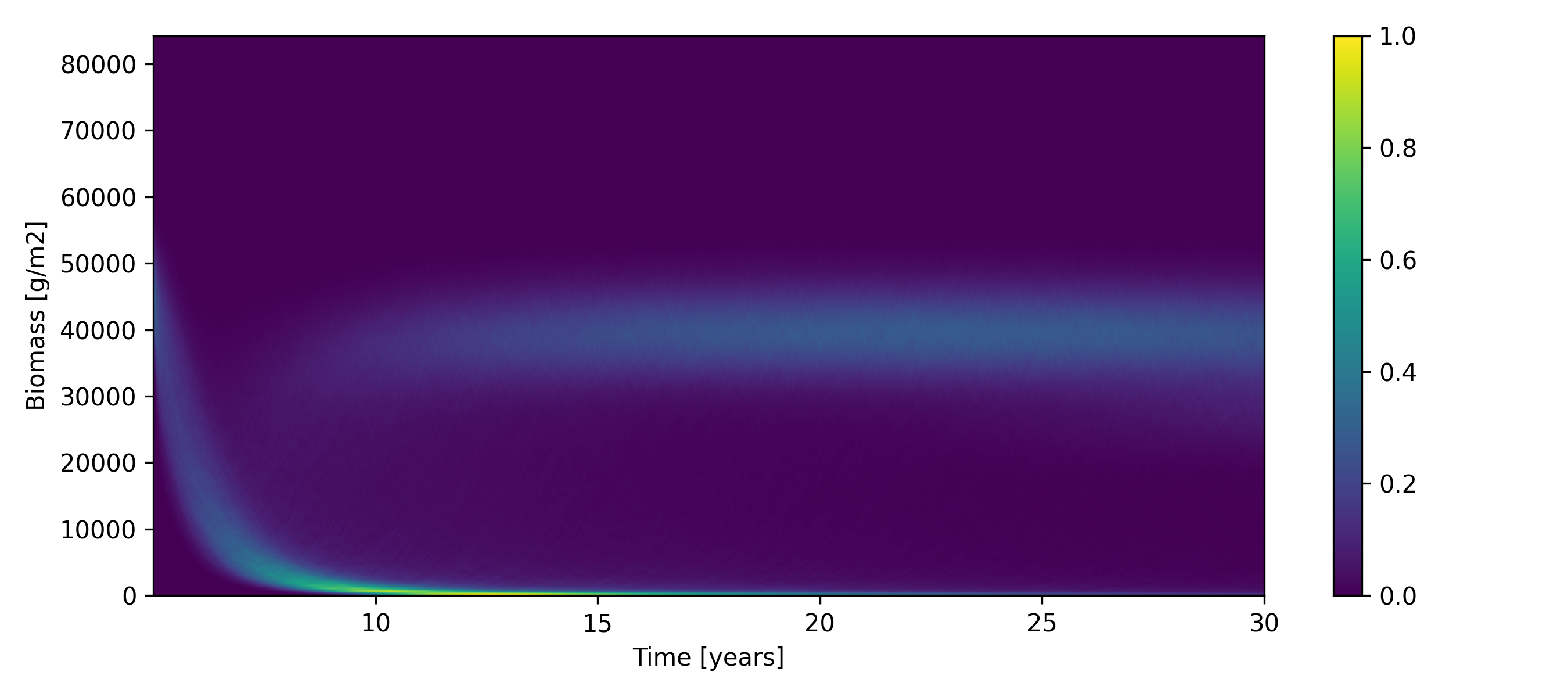}   
     \caption{Harvesting dynamic compliance.}
    \label{fig:dynamic}
    \end{subfigure} 
    \caption{Histograms of the total population (in grams per square meter) over time. }
    \label{fig:histograms-total-population-over-time}
\end{figure}

In Figure \ref{fig:histograms-total-population-over-time}, we observe the empirical distribution of total population over time during the harvesting phase. Fig. \ref{fig:full-compliance} illustrates the scenario under full compliance with harvesting regulations, where the distribution stabilizes around $40000$ grams per square meter. In contrast, Fig.~\ref{fig:dynamic} shows the empirical distribution under dynamic compliance, which appears more dispersed and suggests a potential risk of extinction, as population levels may approach zero with non-negligible probability.

\begin{figure}[h!]
    \centering
    \begin{subfigure}[b]{0.45\textwidth}
     \includegraphics[width=\textwidth]{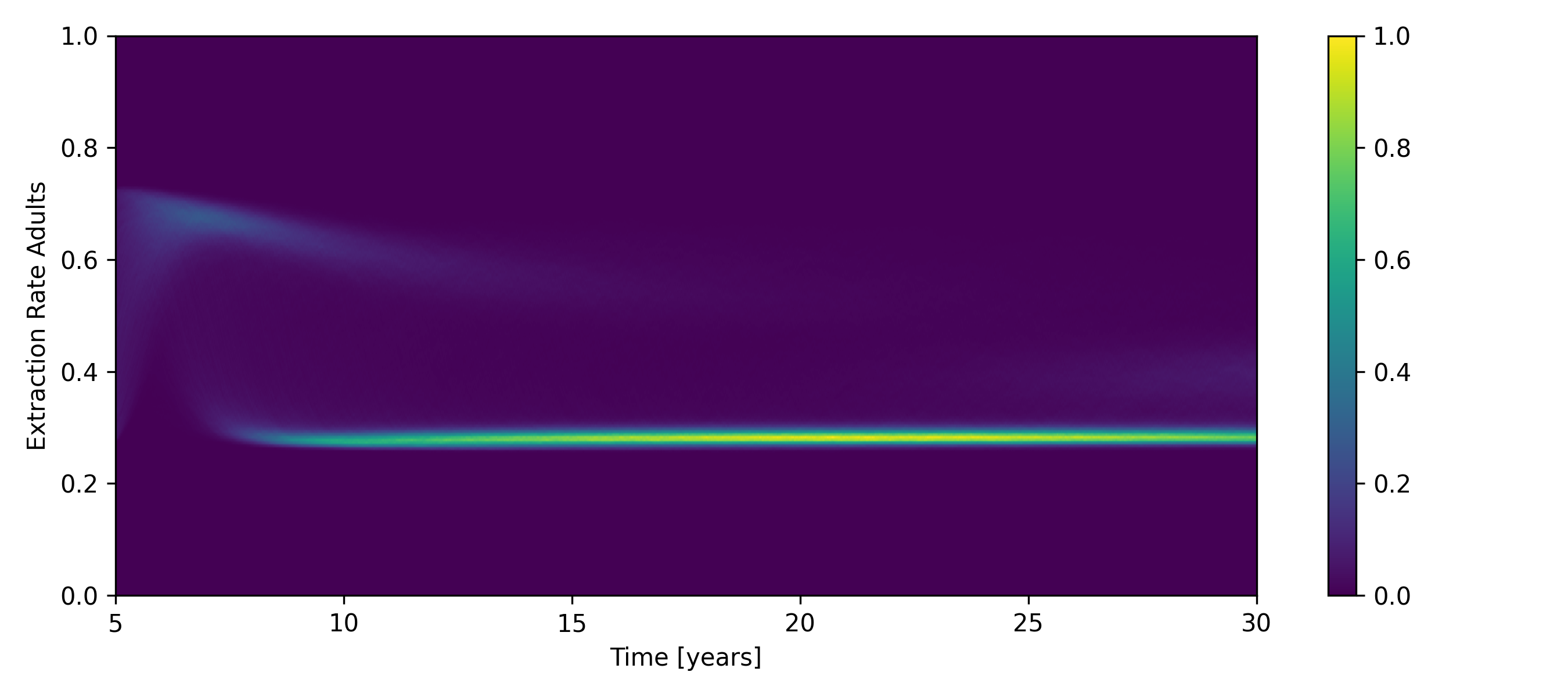}   
     \caption{$\FmaxA$.}
    \label{fig:FA}
    \end{subfigure} 
    \begin{subfigure}[b]{0.45\textwidth}
     \includegraphics[width=\textwidth]{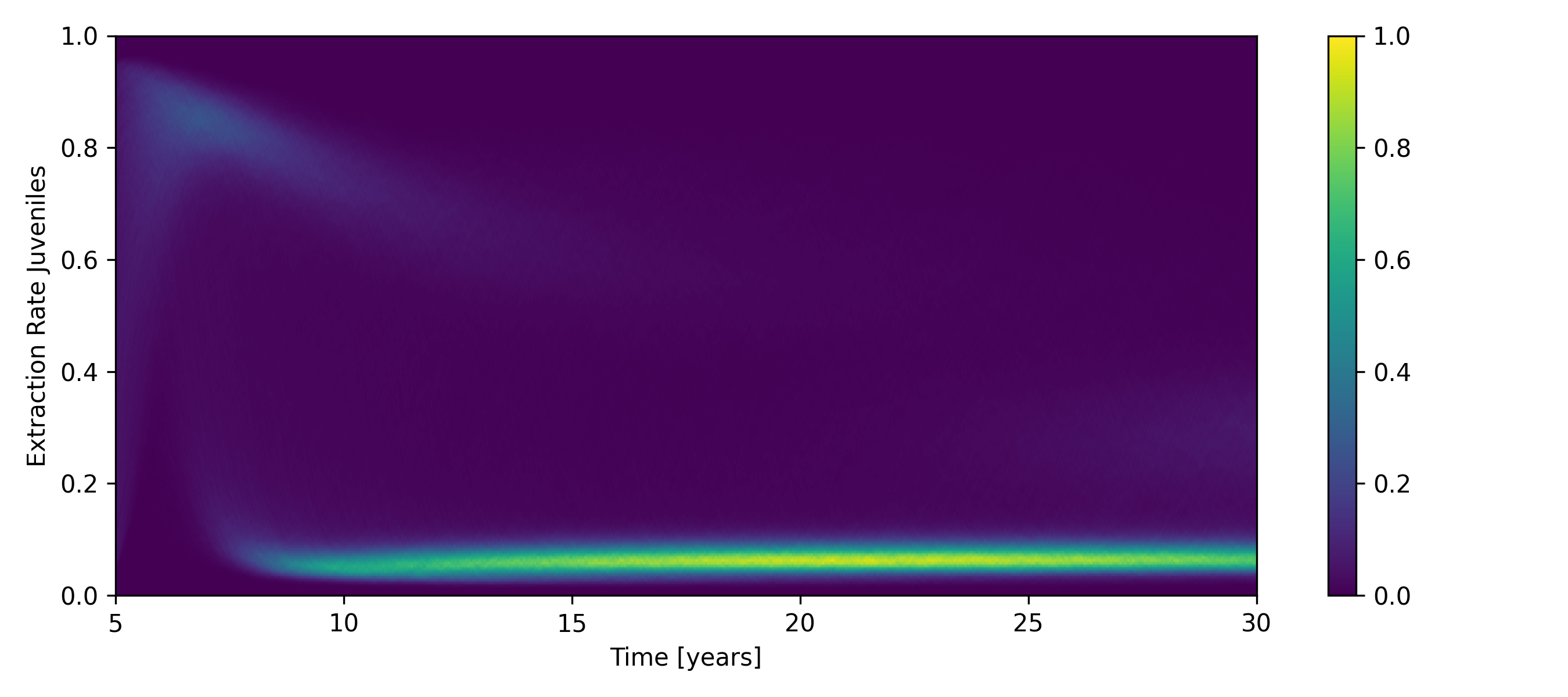}   
     \caption{$\Fmaxj$.}
    \label{fig:FJ}
    \end{subfigure} 
    \caption{Histograms over time for the extraction rate for adults and juveniles for the dynamic model.}
    \label{fig:histograms-E-over-time}
\end{figure}

In Figure~\ref{fig:FA}–\ref{fig:FJ} we present the empirical distribution through time of extraction rates under dynamic compliance for adults and juveniles, respectively. Regulations stipulate that $\FmaxA$ should not exceed $0.25$, and $\Fmaxj$ should be $0$. However, we observe that the probability of violating these thresholds is consistently high, particularly for juveniles. The extraction rate is initially higher but subsequently declines for both populations, likely due to resource abundance rules. Our simulations are consistent with empirical patterns reported across a kelp harvesting gradient \cite{gonzalez2021exploring}: an adult-dominated population in non-take zones (our non-harvesting scenario), a sustainable juvenile-predominated population in management areas under full compliance, and depleted stocks in open-access areas where harvester compliance with regulations is low (our dynamic compliance scenario).

\begin{figure}[h!]
  \centering
  \begin{subfigure}[t]{.4\textwidth} 
    \centering
    \vspace{0pt} 
    \includegraphics[width=\linewidth]{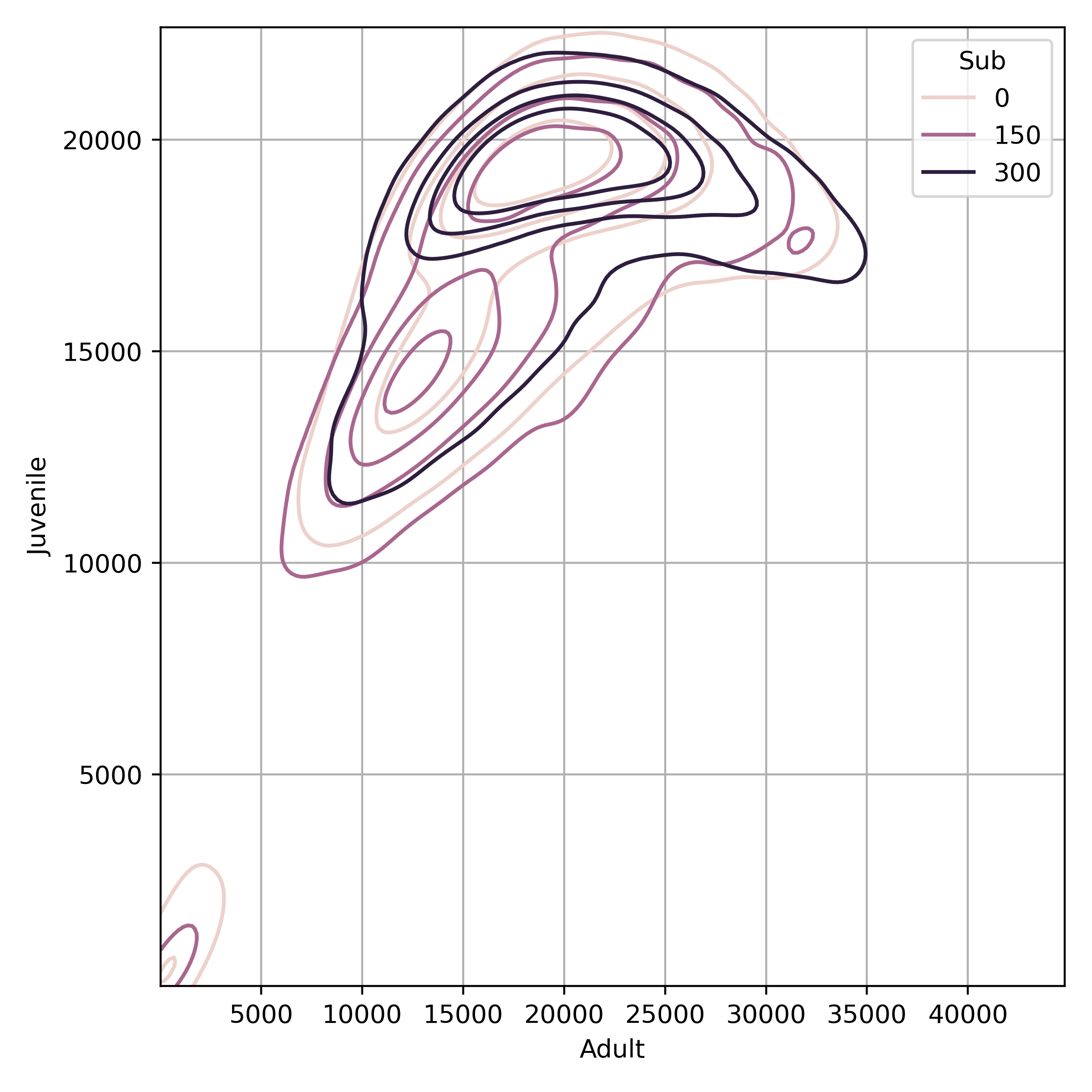}%
    \caption
      {%
        Joint density for juveniles and adults.%
        \label{fig:joint-density-adults-juveniles}%
      }%
  \end{subfigure}\hfill
  \begin{tabular}[t]{@{}c@{}}
    \begin{subfigure}[t]{.5\textwidth} 
      \centering
      \vspace{0pt} 
      \includegraphics[width=\linewidth]{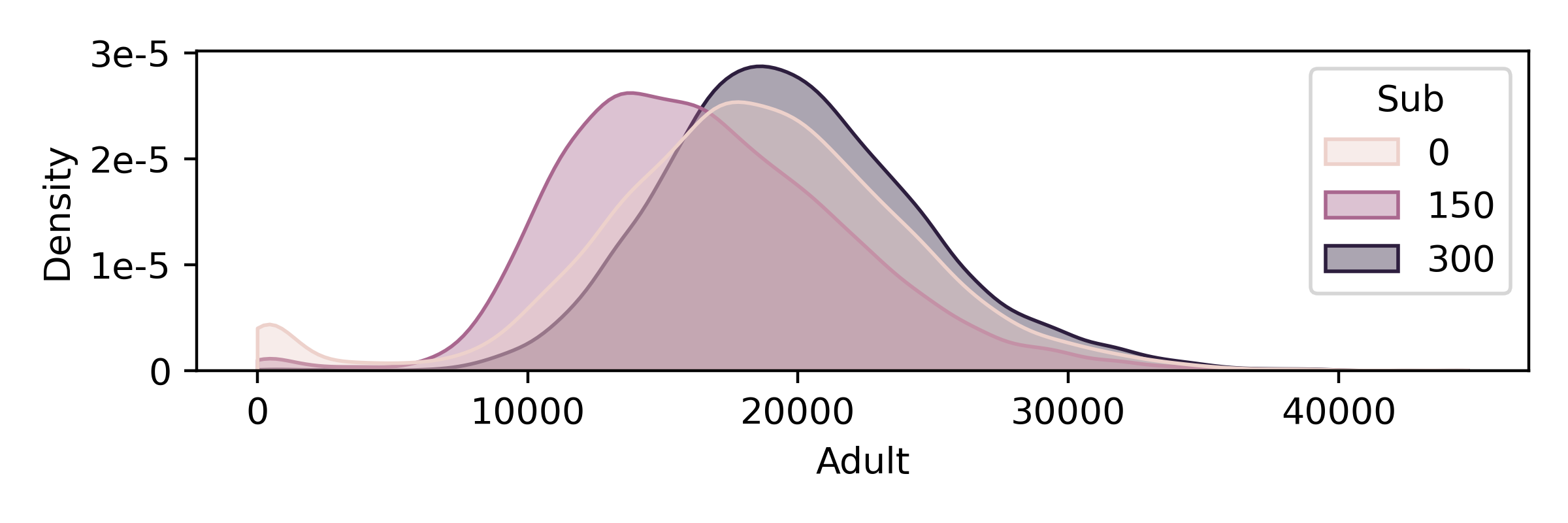}%
      \caption
        {%
          Marginal adult biomass density.%
          \label{fig:marginal-density-adults}%
        }%
    \end{subfigure}\\
    
    \noalign{\smallskip}%
    
    \begin{subfigure}[b]{.5\textwidth} 
      \centering
      \includegraphics[width=\linewidth]{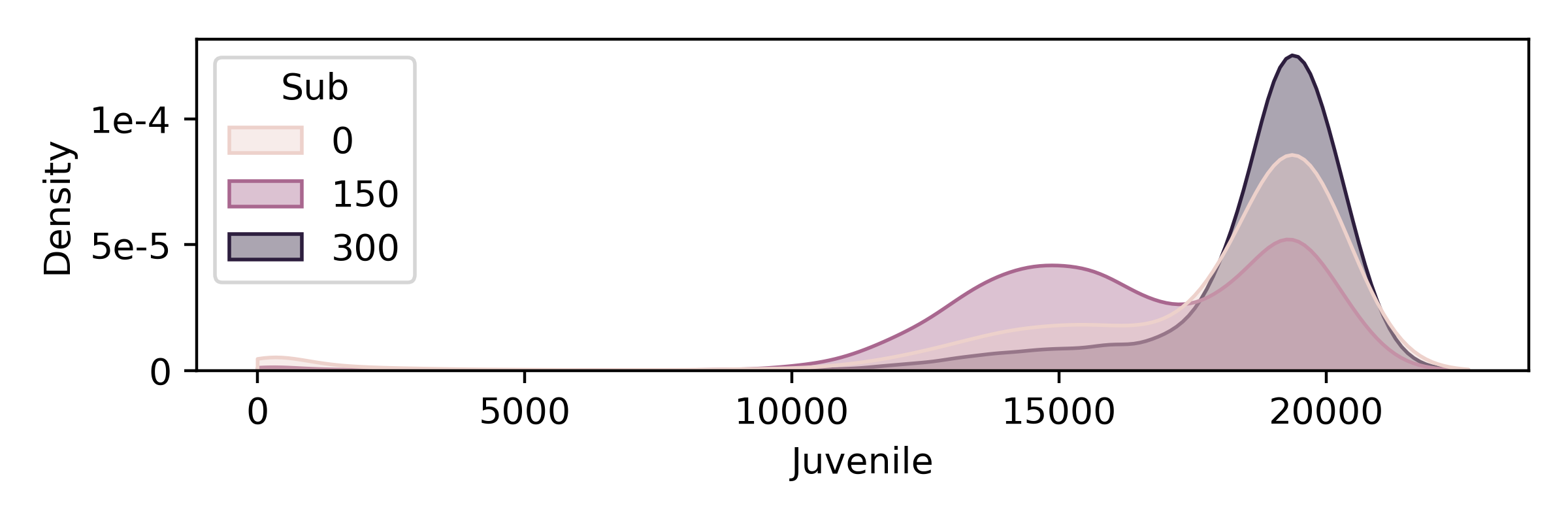}%
      \caption
        {%
          Marginal juvenile biomass density.%
          \label{fig:marginal-density-juveniles}%
        }%
    \end{subfigure}
  \end{tabular}
  
  \caption
    {%
      Estimated joint and marginal densities for $(\XA,\XJ)$ after 25 years of harvesting under different levels of subsidy.%
      \label{fig:joint-and-marginals-A-J}%
    }
\end{figure}

To illustrate how the model can be used to analyze the impact of varying subsidy levels, we compare three scenarios defined by the subsidy parameter $s$ (see Equation~\eqref{eq:def_beta_2}): no subsidy ($s=0$), intermediate subsidy ($s=150$), and high subsidy ($s=300$). Figure~\ref{fig:joint-and-marginals-A-J} presents the model’s joint and marginal empirical densities under these scenarios. The subsidy parameter can be interpreted as a policy instrument to promote compliance with harvesting regulations. Our results indicate that increasing subsidies tends to shift the system toward more sustainable regimes, thereby reducing the risk of kelp collapse and reinforcing compliance. This finding echoes with evidence that integrating economic incentives into regulatory frameworks (such as price premiums from sustainability programs) can enhance fisher compliance, improve marketability, strengthen governance resilience, and support long-term sustainability \cite{fumero2024resilience, jardim2023msc}. Our simulations also help to assess the influence of the level of price premiums needed for compliance and sustainability.

\begin{figure}[h!]
  \centering
  \begin{subfigure}[t]{.4\textwidth} 
    \centering
    \vspace{0pt} 
    \includegraphics[width=\linewidth]{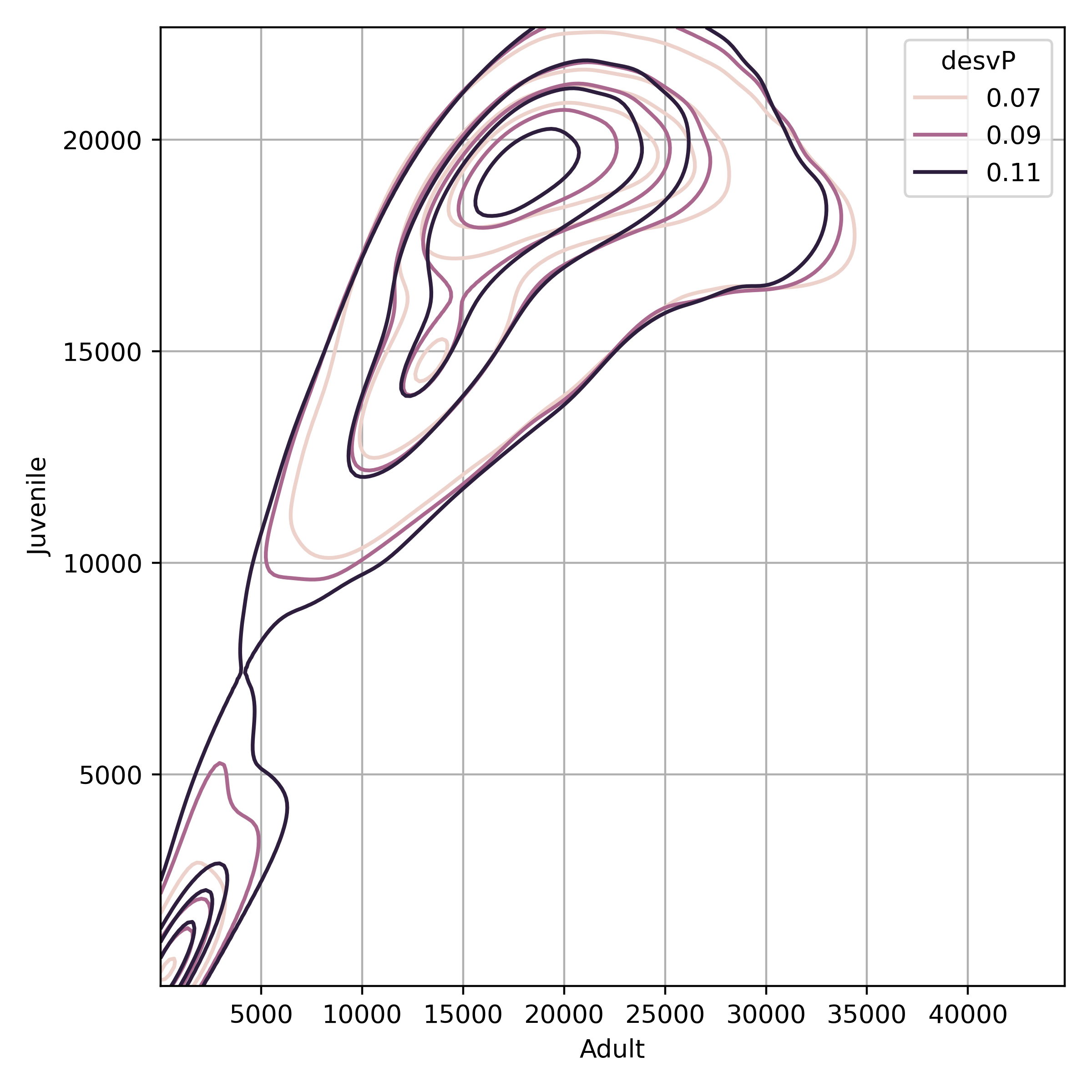}%
    \caption
      {%
        Joint density for juveniles and adults.%
        \label{fig:joint-density-adults-juveniles-sigma}%
      }%
  \end{subfigure}\hfill
  \begin{tabular}[t]{@{}c@{}}
    \begin{subfigure}[t]{.5\textwidth} 
      \centering
      \vspace{0pt} 
      \includegraphics[width=\linewidth]{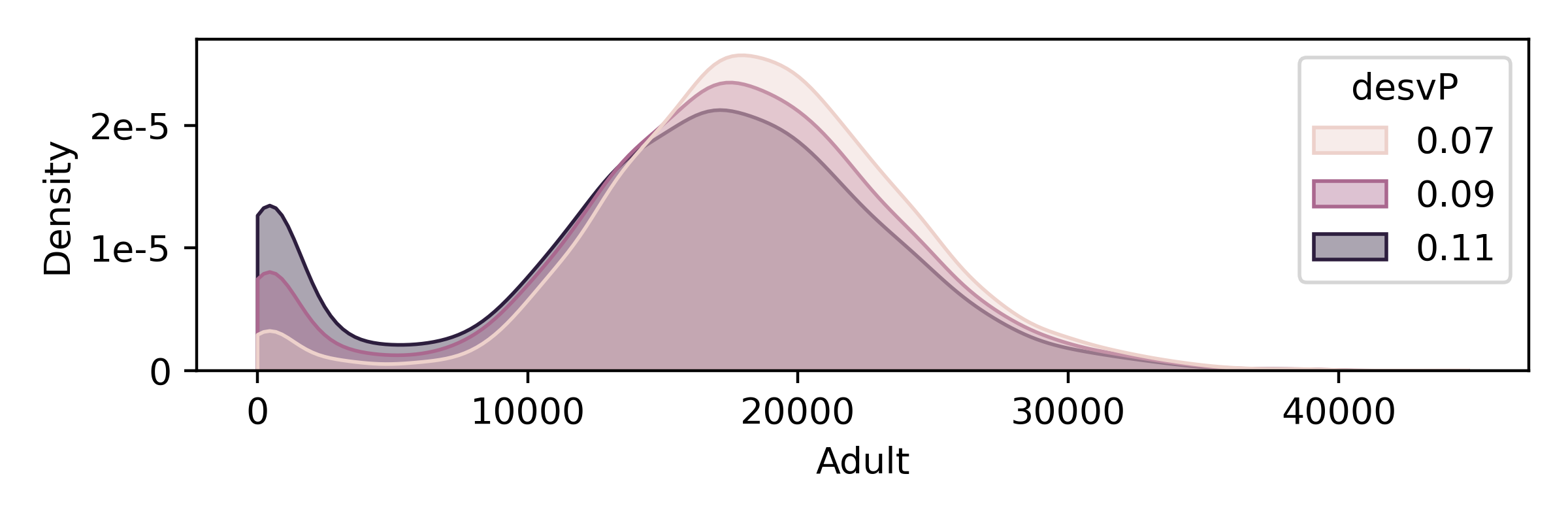}%
      \caption
        {%
          Marginal adult biomass density.%
          \label{fig:marginal-density-adults-sigma}%
        }%
    \end{subfigure}\\
    
    \noalign{\smallskip}%
    
    \begin{subfigure}[b]{.5\textwidth} 
      \centering
      \includegraphics[width=\linewidth]{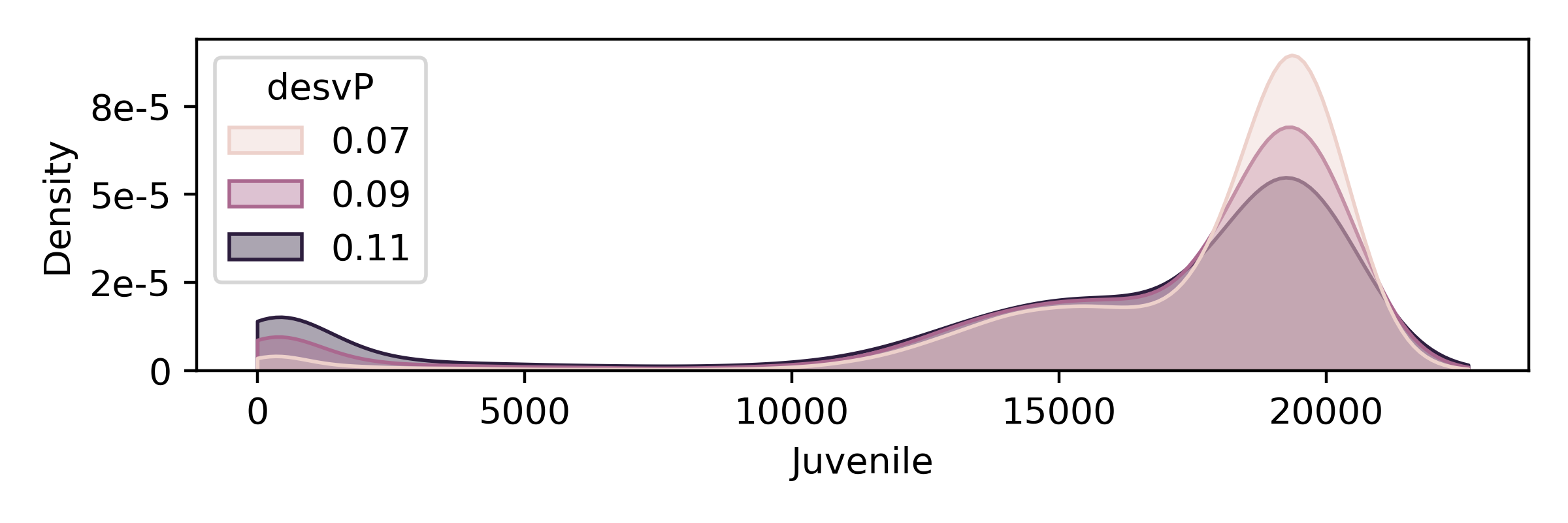}%
      \caption
        {%
          Marginal juvenile biomass density.%
          \label{fig:marginal-density-juveniles-sigma}%
        }%
    \end{subfigure}
  \end{tabular}
  
  \caption
    {%
      Estimated joint and marginal densities for $(\XA,\XJ)$ after 25 years of harvesting under different levels of price volatility.%
      \label{fig:joint-and-marginals-A-J-sigma}%
    }
\end{figure}

Finally, in Figure \ref{fig:joint-and-marginals-A-J-sigma} we compare the estimated joint and marginal densities for $(A,J)$ after 25 five years of exploitation for different levels of price volatility, that is, for $\sigma_P = \{0.07, 0.09, 0.11\}$. The smaller value is the one estimated from price data. Notice that even small variations in this parameter produce considerable changes in the probability to reach the extinction state.

\clearpage

\section{Closing remarks and directions for further research}\label{sec:closing}

The analysis in this work is based on a mean-field assumption for the social network of harvesters, addressing it as a complete graph where each harvester interacts with all others. This simplification may overlook the influence of local peers on harvesters' compliance decisions. Future research could explore more realistic social network structures, such as small-world or scale-free networks, to examine how these properties impact compliance dynamics and the long-term persistence of the kelp population. A deeper understanding of socio-ecological feedbacks could inform more effective management strategies, though analyzing the extended model raises technical challenges, particularly in characterizing the limiting process. 

A second path for future research concerns the rigorous study of stochastic persistence for our system. This paper shows that, under suitable conditions, solutions remain in a biologically meaningful positive domain. Section \ref{sec:inaccessibility_origin_invariant} provides sufficient conditions for invariant measures in the interior and characterizes extinction through exponential stability. A key open problem is to prove that such an invariant measure serves as a global attractor for any initial measure with the same properties. Establishing this would yield theoretical guarantee of long-term kelp persistence under strong stochastic disturbances, thereby reinforcing the ecological relevance of the model.

While the results presented in this article provide a fundamental understanding of the system’s behavior, they offer a generalized view that may not directly translate into actionable policy interventions. In particular, developing management strategies that incorporate subsidies, penalties for non-compliance, collective organization of harvesters, and educational campaigns to raise awareness of the kelp population’s status within an optimal control framework could provide valuable insights into socio-ecological trade-offs and long-term outcomes. Together, these studies would link theoretical insights with practical coastal resource management, providing decision-makers with robust and effective tools. These developments are part of our current research efforts. 

\section*{Acknowledgments and funding}

The authors were supported by ANID-Exploration grant 13220168 ``\textit{Biological and Quantum Open System Dynamics: evolution, innovation and mathematical foundations}''. L.V. has been partially supported by ANID through FONDECYT Iniciaci\'on, project number 11240158-2024 \emph{Adaptive behavior in stochastic population dynamics and non-linear Markov processes in ecoevolutionary modeling}. K.M. and H.O. have been partially supported by INRIA Associated Team SWAM and Programa Regional MathAmdSud AMSUD 240054, \emph{Stochasticity \& Chaos in Multiscale Phenomena}. H.O. has been partially supported by FONDECYT Regular Nº1242001 \emph{Propagation of chaos for mean-field interacting particle systems in mathematical physics and mathematical biology.}
M.I.A-T. has been supported by Fondecyt 3220110 \emph{Uunderstanding non-compliance in kelp fisheries from a socio-ecological perspective} and partially supported by ANID PIA/BASAL AFB240003 \emph{Center of Applied Ecology and Sustainability (CAPES)}, Millennium Science Initiative Program ICN 2019–015 \emph{Instituto Milenio en Socioe-Ecologia Costera (SECOS)}

\appendix

\bibliographystyle{plain} 
\bibliography{biblio}

\section{Appendices: postponed proofs}\label{sec:appendix}

\subsection{Postponed proofs from Section \ref{sec:well-posedness} }\label{sec:proof-sec-well-posedness}
\begin{proof}[Proof of Proposition \ref{prop:moments}]
From the definition of the auxiliary SDE \eqref{eq:system_welldef} and assumptions \ref{H1} and \ref{H4}, it is easy to check that there exists a constant $C$, depending on $\mu,\sigma_{(\cdot)}, \PjA$, such that 
\[2\langle \ZZ, b(\ZZ)\rangle + \|\sigma(\ZZ)\|^2\leq C(\|\ZZ\|^2 + 1).\]

Indeed, for the two first and last components the above statement is straightforward whereas for the third component we have:
\begin{eqnarray*}
2e\Psi(e)\Psi(1-e)\Delta \tilde{\beta}(x,y,p) + \sE^2\Psi(e)\Psi(1-e) &\leq C \Psi(e)\Psi(1-e)(|e|+1)\\
&\leq C (|e|^2+1),
\end{eqnarray*}
and similarly $2e\big(\Psi(1-e)\bar{\beta}_2 - \Psi(e)\bar{\beta}_1\big)$ is bounded from above by $C (|e|^2+1)$ for some non-negative constant $C$.

Therefore, under \ref{hip:linear-growth-phi}, we can apply Theorem 2.2 in \cite{xi2019jump}, obtaining that any solution of the model \eqref{eq:system_welldef} is non-explosive, i.e. for $R>0$, $\tau_R := \inf\{t>0:~ \|\ZZ_t\|\geq R\}$ tends to infinity as $R$ tends to infinity. So it is natural to estimate the $p$th-moments of the solution. 

On the other hand, we notice that $\nabla\|\ZZ\|^p = p \|\ZZ\|^{p-2}\ZZ$ and $\Delta\|\ZZ\|^p = p(p-2)\|\ZZ\|^{p-4}\ZZ\ZZ^T + p \|\ZZ\|^{p-2}I$. Then, 
\begin{eqnarray*}
\texttt{tr}(\Delta \|\ZZ_s\|^p \sigma(\ZZ_s)\sigma^T(\ZZ_s)) &=& p(p-2)\|\ZZ\|^{p-4}_2\texttt{tr}(\ZZ\ZZ^T\sigma(\ZZ)\sigma^T(\ZZ)) + p \|\ZZ\|^{p-2}\|\sigma(\ZZ)\|^2\\
&\leq &  p(p-1)\|\ZZ\|^{p-2}\|\sigma(\ZZ)\|^2,
\end{eqnarray*}
since $\ZZ\ZZ^T$ and $\sigma(\ZZ)\sigma^T(\ZZ)$ are positive semi-definite matrices and we can apply Cauchy-Schwarz inequality.

After a localization argument and the application of the It\^o formula, we have
{
{\begin{align*}
&\E[\|\ZZ_t\|^p] = \E[\|\ZZ_0\|^p]+p \int_0^t\E\left[\langle \|\ZZ_{s}\|^{p-2}\ZZ_{s},b(\ZZ_s)\rangle\right]ds\\
&\quad +\frac12 \int_0^t\E\left[\texttt{tr}(\Delta \|\ZZ_s\|^p\sigma(\ZZ_s)\sigma^T(\ZZ_s))\right]ds+ \lambda\int_0^t\int_\R\E\left[\|\ZZ_{s}+\Phi(\ZZ_s)\|^p-\|\ZZ_{s}\|^p\right]\nu(du)ds.
\end{align*}
}}
Then, from triangular inequality and assumption \ref{H2}, there exists a constant $C$ depending on $p$ such that
{
{\begin{align*}
\E[\|\ZZ_t\|^p]
&\leq \E[\|\ZZ_0\|^p] + \frac{p}2 \int_0^t\E\left[\|\ZZ_{s}\|^{p-2}_2\left(2\langle \ZZ_{s},b(\ZZ_s)\rangle + (p-1)\|\sigma(\ZZ_s)\|^2\right)\right]ds\\
& + C~ \int_{-\infty}^{\infty}\left(\max\{g_{(A)}^p(z), g_{(J)}^p(z)\}+1\right)\nu(dz)~\int_0^t\E\left[\|\ZZ_s\|^p\right]ds\\
&\lesssim  \E[\|\ZZ_0\|^p] + \int_0^t\E\left[\|\ZZ_{s}\|^{p-2}(\|\ZZ_s\|^{2} +1) + \|\ZZ_s\|^{p}\right]ds,
\end{align*}
}}
provided that \ref{H1},\ref{H2} and \ref{H4} are satisfied and $g_{(A)}, g_{(J)}\in L^p(\nu)$.

Finally, from Young's and Gronwall's inequalities we conclude that, for all $t\in[0,T]$, there exists a constant $C(p,T)$, depending on the parameters of the model and the $L^p(\nu)$-norm of $g_{(\cdot)}$, such that:

\begin{align*} \label{eq:p-estimate}
\sup_{t\in[0,T]}\E[\|\ZZ_t\|^p]< C(p,T).
\end{align*}
\end{proof}

\begin{proof}[Proof of Corollary \ref{cor:sup-norm_moments}]
The control of uniform moments is a direct consequence of Proposition \ref{prop:moments}. Indeed, for the process $(\XJ_t,\, 0\leq t\leq T)$ we use assumptions \ref{H1}, \ref{H2}, Jensen, Burkholder-Davis-Gundy and Kunita's inequality for jump process, obtaining:
\begin{align*}
&\E[\sup_{0\leq t\leq T}\XJ_t^{2p}] \lesssim  \E[\XJ_0^{2p}]\\
&\, + \E\left[\left(\sup_{0\leq t\leq T}\int_0^t\left[\rrJ(\K_s, E_s)\XA_s - \kappaJ(\Psi(E_s))\XJ_s + \int_{\R\setminus\{0\}} \phiJ (\XJ_{s}, z)\lambda\nu(dz)\right]ds\right)^{2p}\right] \\
&\,+ \E\left[\left(\sup_{0\leq t\leq T}\int_0^t\XJ_sdW^{(J)}_s\right)^{2p}\right] + \E\left[\left(\sup_{0\leq t\leq T}\int_0^t\int_{\R\setminus\{0\}} \phiJ (\XJ_{s^{-}}, z)\widetilde{\Np}(dz\otimes ds) \right)^{2p}\right]\\
&\qquad \qquad \leq  C_p\E[\XJ_0^{2p}]+ C_p\E\left[\int_0^T\left(\XA_s^{2p} + \XJ_s^{2p} \right)ds\right]  + C_{\text{\tiny{BDG}}}\E\left[\left(\int_0^T\XJ_s^{2}ds\right)^{p}\right]\\
&\quad + C_K\E\left[\left(\int_0^T\int_{\R} \phiJ^2 (\XJ_s, z)\lambda\nu(dz)ds\right)^{p}\right]+C_K\E\left[\int_0^T\int_{\R} \phiJ^{2p} (\XJ_s, z)\lambda\nu(dz)ds\right],
\end{align*}
for some non-negative constants $C_p,C_{\text{\tiny{BDG}}}, C_K$ depending on the parameters of the model, $p$ and the terminal time $T$, and $\widetilde{\Np} = \Np - \lambda\nu $ the compensated Poisson random measure.

From Jensen and Proposition \ref{prop:moments}, we get:
\begin{align*}
\E[\sup_{0\leq t\leq T}\XJ_t^{2p}] &\leq  C \left(\E[\XJ_0^{2p}]+ \sup_{0\leq s\leq T}\E[\XA_s^{2p}] +\sup_{0\leq s\leq T}\E[\XJ_s^{2p}]\right),
\end{align*}
which is finite. 
\end{proof}

\begin{proof}[Proof of Proposition \ref{prop:pathwise_uniqueness} ]
Notice that, by adding the necessary pivots, we can compute some preliminary bounds as follows. Denoting $\Delta x = x_1 - x_2$, $\Delta y = y_1 - y_2$, and $\Delta e = e_1 - e_2$ we observe that from \ref{H1} and the Lipschitz continuity of $\Psi$ 
that, for any $N\geq1$ and all $x,y,e\in B(0,N)$, there exists a constant $C_N>0$, linear in $N$, such that
\[(x_1-x_2)\Big(\langle u_1, b(x_1,y_1,e_1)\rangle - \langle u_1, b(x_2,y_2,e_2)\rangle\Big)\leq C_N((\Delta x)^2 + (\Delta y)^2 + (\Delta e)^2),\]
where $u_i$ is the canonical vector in $\R^4$ in the direction $i$ and $b$ is the drift vector in Equation \eqref{eq:system_welldef} without the cut-off since we already know that the compliance process is in the domain $(0,1)$.  

Analogously, we have 
\[(y_1-y_2)\Big(\langle u_2, b(x_1,y_1,e_1,p_1)\rangle - \langle u_2, b(x_2,y_2,e_2,p_2)\rangle\Big)\leq C_N~ \|\Delta\|^2_2,\]
where $C_N$ is a constant depending linearly on $N$, for all $x,y,e\in B(0,N)$.

Furthermore, using the boundness and Lipschitz property of the map $(x,y,p)\mapsto \Delta \tilde{\beta}(x,y,p)$, we have 
\begin{align*}
&(e_1-e_2)\Big(\langle u_3, b(x_1,y_1,e_1,p_1)\rangle - \langle u_3, b(x_2,y_2,e_2,p_2)\rangle\Big)\leq C~ |\Delta z|\, \|\Delta\|_1\leq C~ \|\Delta\|_2^2,\\
&(p_1-p_2)\Big(\langle u_4, b(x_1,y_1,e_1,p_1)\rangle - \langle u_4, b(x_2,y_2,e_2,p_2)\rangle\Big)\leq C|\Delta p|^2,
\end{align*}
and finally, from the H\"older property of the square root function:
{\small{
\begin{eqnarray*}
|\sE \sqrt{e_1(1-e_1)} - \sE \sqrt{e_2(1-e_2)}|^2 &\leq & 3\sE^2 |\Delta z|.
\end{eqnarray*}
}}

In the same fashion, it is easy to check that there exists a constant $L_N>0$ linear in $N$ such that:
\[\Big|\langle u_i, b(x_1,y_1,z_1,p_1)\rangle - \langle u_i, b(x_2,y_2,z_2,p_2)\rangle\Big|\leq L_N~ \|\Delta\|_1,\]
for $i=1,2,3,4.$

With these inequalities we are ready to prove the strong uniqueness of the system \eqref{eq:system_welldef}.

Let $\psi_n$ be the regularization of the absolute value function introduced by Yamada and Watanabe (see e.g., \cite[Proposition 2.13, pp. 291]{karatzas98}) and define the function
$$H(x,y,p,z):=\|\mathbf{x}\|^2_2+\psi_n(x)+\psi_n(y)+\psi_n(p)+\psi_n(z).$$

Thus, $H_n(\ZZ)\rightarrow \|\ZZ\|_2^2 +\|\ZZ\|_1^1$ as $n\rightarrow+\infty$ monotonically.

For $i=1, 2$, let $\ZZ^{(i)}_t$ be two solutions of the Equation \eqref{eq:full-model} with the same initial condition $\ZZ_0^{(i)} = \ZZ_0$. We also introduce the notations $\Delta_A:=\XA^{(1)}-\XA^{(2)}$, $\Delta_J=\XJ^{(1)}-\XJ^{(2)}$, $\Delta_P =P^{(1)}-P^{(2)} $, $\Delta_E = E^{(1)}-E^{(2)}$, and for a function $G$, $\Delta_{G}:=G(\K^{(1)},E^{(1)})-G(\K^{(2)},E^{(2)})$.  Then
    
    \begin{equation}\label{eq:eq-for-delta}
\begin{aligned}
    d\Delta_J(t) &= \Delta_{\FJ}(t)dt + \sJ\Delta_J(t)d\WJ_t +\int_{\R} \Delta_{\phiJ}(t,z)\Np(dz\otimes dt),\\
    d\Delta_A(t) &= \Delta_{\FA}(t)dt + \sA\Delta_A(t)d\WA_t +\int_{\R} \Delta_{\phiA}(t,z)\Np(dz\otimes dt), \\
    d\Delta_E(t) &= \left(-(\overline{\beta}_2+\overline{\beta}_1)\Delta_E(t)+\Delta_{\widetilde{\beta}}(t)\right)dt \\
    &\qquad +  \sE\left(\sqrt{E^{(1)}_t(1-E^{(1)}_t)}-\sqrt{E^{(2)}_t(1-E^{(2)}_t)}\right)dW^{(E)}_t, \\
    d\Delta_P(t) &= \Delta_\mu(t)dt + \Delta_{\sP}(t) dW_t,
\end{aligned} 
\end{equation}
where $\FJ(x,y,e,p) = \rrJ(x,y,e)y - \kappaJ(e)x$, $\FA(x,y,e,p) = \rho_A x - \kappaA(e)y$ and $\widetilde{\beta}(x,y,e,p) = \Delta \widetilde{\beta}(x,y,p)e(1-e)$

Let $\tau_N$ be the first time any of the components of the processes $\ZZ^{(1)}$ or $\ZZ^{(2)}$ leaves the ball of radius $N$, i.e. $\tau_N=\inf\{t\in[0,T]:~ \max_{i,j}\langle u_j,\ZZ^{(i)}_t\rangle > N, i=1,2, j=1,\ldots,4\}.$ 

Denoting $\Delta \ZZ = [\Delta_J, \Delta_A,\Delta_E,\Delta_P]$ and $\Delta_b = [\Delta_{\FA}, \Delta_{\FJ},-(\overline{\beta}_2+\overline{\beta}_1)\Delta_E + \Delta_{\widetilde{\beta}}, \Delta_{\mu}]$, according to It\^o's Formula, and using that the stochastic integrals with respect to the Brownian motions are local martingales, we have
\begin{equation}\label{eq:terms_uniqueness}
\begin{aligned}
    &\E[H_n\left(\Delta\ZZ_{t \wedge \tau_R}\right)] 
    = 
2\E\left[\int_0^{t\wedge\tau_N}\langle\Delta{\ZZ} , \Delta_b\rangle(s) ds\right] \nonumber\\
&\quad + \E\left[\int_0^{t\wedge\tau_N} \left(\sJ^2\Delta_J^2(s) +\sA^2\Delta_A^2(s) + \Delta_{\sP}^2(s) \right.\right.\\
&\qquad \qquad\qquad \left.\left.+ \sE^2 \left(\sqrt{E^{(1)}_s(1-E^{(1)}_s)}-\sqrt{E^{(2)}_s(1-E^{(2)}_s)}\right)^2\right)ds\right]
    \nonumber\\
    &\quad +
    \E\left[\int_0^{(t\wedge\tau_N)^+}\int_{\R} \Big((\Delta_A(s-) + \Delta_{\phiA}(s,z) )^2 - \Delta_A(s-)^2   \Big)\Np(dz\otimes ds)\right]\nonumber\\
    &\quad +
    \E\left[\int_0^{(t\wedge\tau_N)^+}\int_{\R} \Big((\Delta_J(s-) + \Delta_{\phiJ}(s,z) )^2 - \Delta_J(s-)^2   \Big)\Np(dz\otimes ds)\right]\nonumber\\
    &\quad+\E\left[\int_0^{t\wedge\tau_N} \psi_n'(\Delta_A(s))\Delta_{\FA}(s)+\psi_n'(\Delta_J(s))\Delta_{\FJ}(s)+\psi_n'(\Delta_P(s))\Delta_\mu(s)ds\right]
    \nonumber\\
    &\quad 
    +  \E\left[\int_0^{t\wedge\tau_N} \psi_n'(\Delta_E(s))\left(-(\overline{\beta}_2+\overline{\beta}_1)\Delta_E(s)+\Delta_{\widetilde{\beta}}(s)\right)ds\right] 
    \nonumber\\
    &\quad 
    +\frac12\E\left[\int_0^{t\wedge\tau_N} \left(\psi_n''(\Delta_A(s))(\sA\Delta_A(s))^2 + \psi_n''(\Delta_J(s))\right.\right.\\
    &\qquad\qquad\left.\left.\times (\sJ\Delta_J(s))^2 + \psi_n''(\Delta_P(s))(\Delta_{\sP}(s))^2\right)ds\right]
    \nonumber\\
    &\quad + \frac{\sE^2}2\E\left[ \int_0^{t\wedge\tau_N} \psi_n''(\Delta_E(s))\Bigg(\sqrt{E^{(1)}_s(1-E^{(1)}_s)}-\sqrt{E^{(2)}_s(1-E^{(2)}_s)} \Bigg)^2ds\right]\nonumber\\
    &\quad +\E\left[\int_0^{(t\wedge\tau_N)^+}\int_{\R} \Big(\psi_n\big(\Delta_A(s-) + \Delta_{\phiA}(s,z) \big) - \psi_n(\Delta_A(s-))  \Big)\Np(dz\otimes ds)\right]\nonumber\\
   &\quad +\E\left[\int_0^{(t\wedge\tau_N)^+}\int_{\R} \Big(\psi_n\big(\Delta_J(s-) + \Delta_{\phiJ}(s,z) \big) - \psi_n(\Delta_J(s-))  \Big)\Np(dz\otimes ds)\right]\nonumber\\
   &\qquad =: \sum_{k=1}^{10}T_k.
\end{aligned}
\end{equation}

We bound now each term $T_k$, for $k\in\{1,\ldots,10\}$. First, using the preliminary bounds and recalling that $\|\psi_n'\|_\infty\leq 1$, we can easily bound the Lebesgue integrals in $T_1+T_2+T_5+T_6$:
\begin{align*}
    &T_1+T_2+T_5+T_6\\
    &\qquad \leq C_N~\E\left[\int_0^{t\wedge\tau_N}(\|\Delta{\ZZ}(s)\|_2^2 + |\Delta_E(s)|)ds\right] + L_N\E\left[\int_0^{t\wedge\tau_R} \|\Delta_{\ZZ}(s)\|_1 ds\right],
\end{align*}
where, for any $N\geq1$, the constants $C_N,L_N>0$ depend on the model parameters and are linear in $N$.

Now, we recall that, for all $n\geq1$,$\psi_n''(x)=\rho_n(x)$ with $\rho_n$ supported in $(a_{n+1},a_n)$, for some $a_n>0$. Moreover, $0\leq \rho_n(x) \leq 2/nx$ for all $x\geq0$, thus:
\begin{align*}
&\E\int_0^{t\wedge\tau_N} \frac{1}{2}\psi_n''(\Delta_A(s))(\sA\Delta_A(s))^2ds\\
&\qquad \leq \E\int_0^{t\wedge\tau_R} \frac{1}{n\Delta_A(s)}\ind{(a_{n+1},a_n)}(\Delta_A(s))(\sA\Delta_A(s))^2ds \lesssim \frac{1}{n}\int_0^{t} \E\bigg|\Delta_A(s\wedge\tau_N)\bigg|ds.
\end{align*}
The same argument and the H\"older continuity of the square root shows that:
\begin{align*}
   T_7 + T_8 \lesssim\frac{1}{n}\left(1+\int_0^{t} \E\big\|\Delta_{\ZZ}(s\wedge\tau_N)\big\|_1ds\right).
\end{align*}
It remains to bound the terms coming from the jumps, for which we use Assumption \ref{hip:lipschitz-condition-phi} and Young inequality:
\begin{align*}
    T_3 
    &\leq \lambda\E\int_0^{t\wedge\tau_N} \left(2|\Delta_A(s)|\int_{\R}|\Delta_{\phiA}(s,z)|\nu(dz) + \int_{\R}(\Delta_{\phiA}(s,z) )^2   \nu(dz)\right)ds\\
    & \lesssim \int_0^{t} \E|\Delta_A(s\wedge\tau_R)|^2 ds.
\end{align*}

In the same way, since $\psi_n$ is globally $1$-Lipschitz:
\begin{align*}
T_9=\E&\int_0^{(t\wedge\tau_N)^+}\int_{\R} \Big(\psi_n\big(\Delta_A(s-) + \Delta_{\phiA}(s,z) \big) - \psi_n(\Delta_A(s-))  \Big)\Np(dz\otimes ds)\\
&\lesssim \int_0^{t} \E|\Delta_A(s\wedge\tau_R)| ds.
\end{align*}

With the same arguments, we can estimate the remaining terms $T_4$ and $T_{10}$, obtaining that
\begin{align*}
     &\E \left[H_n\left(\Delta\ZZ_{t \wedge\tau_N}\right) \right]\\
     &\quad \leq C_N\int_0^{t}\left(\E \|\Delta\ZZ_{s\wedge\tau_N}\|_2^2+ \E\|\Delta\ZZ_{s\wedge\tau_N}\|_1\right)ds 
     +\frac{C_N}{n}\int_0^{t} \E\big\|\Delta{\ZZ}_{s\wedge\tau_N}\big\|_1ds +\frac{C_N}{n},
\end{align*}
for some constant $C_N$ depending on the parameters of the model, $T$ and $N$.

Taking limit when $n$ goes to infinity, and recalling that $H_n(\cdot)$ converges monotonically to $\|\cdot\|_2^2+\|\cdot\|_1$ it follows that, for $N\geq1$ fixed:
\begin{align*}
     \E &\left[\left\|\Delta\ZZ_{t \wedge \tau_N}\right\|_2^2  + \left\|\Delta\ZZ_{t \wedge \tau_N}\right\|_1 \right]\leq C_N\int_0^t  \E \left[\left\|\Delta\ZZ_{s \wedge \tau_N}\right\|_2^2+\left\|\Delta\ZZ_{s \wedge \tau_N}\right\|_1 \right] ds.
\end{align*}
We finish the proof applying Gronwall's lemma and then Fatou's Lemma to take the limit $N\to\infty$. Then, we conclude from the c\`adl\`ag property of the trajectories and the fact that $\tau_N\rightarrow\infty$, that $\ZZ^{(1)}_t = \ZZ^{(2)}_t$, for $\P$-almost all $t\geq0$.
\end{proof}

\subsection{Postponed proofs of Section \ref{sec:numerical_scheme}}\label{sec:delayed-proofs-ns}

\subsubsection{Preliminary estimations}

\begin{proposition}\label{prop:control-of-brownian-exponentials}
Let us consider $\sigma>0$ and $\eta(t)\leq s_1\leq s_2\leq t$:
$$
{\bf e}(s_1,s_2):=\exp\left(-\frac{\sigma^2}{2}(s_2-s_1) + \sigma \deltaWA{s_1}{s_2}\right),
$$
and
$$
I(\eta(t),t):=\int_{\eta(t)}^{t}\left({\bf e}(s,t) -{\bf e}(\eta(t),t)\right)ds.
$$

Then, ${\bf e}(s_1,s_2)$ is measurable with respect to $\F_{s_2}$, and we have 
\begin{align*}
&\E\left[{\bf e}(s_1,t){\bf e}(s_2,t) \big|  \F_{s_1} \right] =e^{\sigma^2(t-s_2)},     \\
&\E\left[I(\eta(t),t){\bf e}(\eta(t),t) \big| \F_{\eta(t)} \right] = \frac{e^{\sigma^2\ddt}-1}{\sigma^2}-\ddt e^{\sigma^2\ddt} = \mathcal{O}\left(\ddt^2 \right),\\
&\E\left[I^2(\eta(t),t) \big| \F_{\eta(t)} \right] =  \frac{2e^{\sigma^2\ddt}}{\sigma^4}-\frac{2}{\sigma^4} - \left(\frac{2}{\sigma^2} -\left(t-\eta(t) \right)\right)\left(t-\eta(t) \right)e^{\sigma^2\ddt},
\end{align*}
and thus 
\(\E\left[I^2(\eta(t),t) \big| \F_{\eta(t)} \right]= \mathcal{O}\left(\ddt^3 \right).\)
\end{proposition}

\begin{proof} 
Just from its definition, 
\begin{align*}
{\bf e}(s_1,t)
     & =
  \exp\left(-\frac{\sigma^2}{2}(t-s_2)-\frac{\sigma^2}{2}(s_2-s_1) + \sigma (\WA_{t}- \WA_{s_2})+\sigma (\WA_{s_2}- \WA_{s_1})\right)     \\
  & = {\bf e}(s_2,t){\bf e}(s_1,s_2), 
\end{align*}
from where it follows:
\begin{align*}
{\bf e}(s_2,t){\bf e}(s_1,t) = {\bf e}^2(s_2,t){\bf e}(s_1,s_2), 
\end{align*}
for any $\eta(t)\leq s_1\leq s_2\leq t$. In particular, for any $\eta(t)\leq s\leq t$ we have
\begin{align*}
{\bf e}(s,t){\bf e}(\eta(t),t)
  & = {\bf e}^2(s,t){\bf e}(\eta(t),s).
\end{align*}

Furthermore, from the independence of Brownian increments, and the form of the moment generating function for a Gaussian random variable we obtain
$$
\E\left[ {\bf e}(s,t) / \F_s\right] = 1,\;\; \E\left[ {\bf e}^2(s,t) / \F_s\right]  = e^{\sA^2(t-s)}.
$$

Combining the previous computations, we get for any $\eta(t)\leq s\leq t$:
\begin{align*}
\E\left[{\bf e}(s,t){\bf e}(\eta(t),t) \big|  \F_{\eta(t)} \right] 
&=  \E\left[ {\bf e}^2(s,t){\bf e}(\eta(t),s)\big|  \F_{\eta(t)} \right]    \\
&=  \E\left[ \E\left[ {\bf e}^2(s,t){\bf e}(\eta(t),s)\big|  \F_{s} \right] \big|  \F_{\eta(t)} \right]    \\
&=  \E\left[ {\bf e}(\eta(t),s) \E\left[ {\bf e}^2(s,t)\big|  \F_{s} \right] \big|  \F_{\eta(t)} \right]    \\
&=  \E\left[ {\bf e}(\eta(t),s) e^{\sA^2(t-s)}\big|  \F_{\eta(t)} \right]    \\
  & =e^{\sA^2(t-s)},     
\end{align*}
and more generally, 
\begin{align}\label{eq:exponential_product_expectation}
\E\left[{\bf e}(s_1,t){\bf e}(s_2,t) \big|  \F_{s_1} \right]& =e^{\sigma^2(t-s_2)}.     
\end{align}

Moreover,
\begin{align*}
&\E\left[I(\eta(t),t){\bf e}(\eta(t),t) \big| \F_{\eta(t)} \right]  \\
&\qquad =   \E\left[\int_{\eta(t)}^{t}({\bf e}(s,t){\bf e}(\eta(t),t) -{\bf e}^2(\eta(t),t))ds \big|  \F_{\eta(t)} \right]    \\
&\qquad=  
\int_{\eta(t)}^{t}\E\left[{\bf e}(s,t){\bf e}(\eta(t),t) \big|  \F_{\eta(t)} \right] ds  -\ddt\E\left[{\bf e}^2(\eta(t),t) \big|  \F_{\eta(t)} \right]    \\
&\qquad =  
\int_{\eta(t)}^{t}e^{\sA^2(t-s)} ds  -\ddt e^{\sA^2\ddt}   \\
  &\qquad  =   \frac{1}{\sA^2}\left(e^{\sA^2\ddt}-1-\sA^2\ddt e^{\sA^2\ddt} \right),
\end{align*}
which is of order $\ddt^2$.

Finally
\begin{align*}
&\E\left[I^2(\eta(t),t) \big| \F_{\eta(t)} \right] \\
  &\qquad   = \E\left[  \int_{\eta(t)}^{t}\int_{\eta(t)}^{t}\left({\bf e}(s_1,t) -{\bf e}(\eta(t),t) \right)\left({\bf e}(s_2,t) -{\bf e}(\eta(t),t) \right)ds_1ds_2  \big| \F_{\eta(t)} \right]   \\
    & \qquad  = \E\left[  \int_{\eta(t)}^{t}\int_{\eta(t)}^{t}{\bf e}(s_1,t) {\bf e}(s_2,t)ds_1ds_2  \big| \F_{\eta(t)} \right]    
  \\
  & \qquad \quad- \E\left[  \int_{\eta(t)}^{t}\int_{\eta(t)}^{t}{\bf e}(s_1,t) {\bf e}(\eta(t),t)ds_1ds_2  \big| \F_{\eta(t)} \right]         \\
  & \qquad \quad  - \E\left[  \int_{\eta(t)}^{t}\int_{\eta(t)}^{t}{\bf e}(\eta(t),t)\left({\bf e}(s_2,t) -{\bf e}(\eta(t),t) \right)ds_1ds_2  \big| \F_{\eta(t)} \right]   \\   
  &\qquad = I_1-I_2-I_3.
\end{align*}

We compute term by term with the help of \eqref{eq:exponential_product_expectation}
\begin{align*}
 I_1
&    =   \int_{\eta(t)}^{t}\int_{\eta(t)}^{s_2}  \E\left[{\bf e}(s_1,t) {\bf e}(s_2,t) \big| \F_{\eta(t)} \right]ds_1ds_2\\
&\qquad\ + \int_{\eta(t)}^{t}\int_{s_2}^{t}\E\left[{\bf e}(s_1,t) {\bf e}(s_2,t) \big| \F_{\eta(t)} \right] ds_1ds_2   \\
&    =   \int_{\eta(t)}^{t}\int_{\eta(t)}^{s_2} e^{\sA^2(t-s_2)}ds_1ds_2   + \int_{\eta(t)}^{t}\int_{s_2}^{t}e^{\sA^2(t-s_1)} ds_1ds_2.
\end{align*}

Integrating we get
\begin{align*}
 I_1
&  =  \frac{-2\sA^2\ddt + 2e^{\sA^2\ddt}-2}{\sA^4}.    
\end{align*}

Similarly,
\begin{align*}
I_2
& = \int_{\eta(t)}^{t}\int_{\eta(t)}^{t}  \E\left[ {\bf e}(s_1,t) {\bf e}(\eta(t),t) \big| \F_{\eta(t)} \right]  ds_1ds_2     \\
& = \int_{\eta(t)}^{t}\int_{\eta(t)}^{t}e^{\sA^2(t-s_1)}  ds_1ds_2     \\
&   =  (t- \eta(t))\frac{e^{\sA^2\ddt} - 1}{\sA^2}.
\end{align*}
\begin{align*}
I_3
 & = \ddt\E\left[ I(\eta(t),t){\bf e}(\eta(t),t)\big| \F_{\eta(t)} \right] = \ddt\frac{e^{\sA^2\ddt}-1}{\sA^2}-\ddt^2e^{\sA^2\ddt}.
\end{align*}

Hence,
\begin{align*}
 \E\left[I^2(\eta(t),t) \big| \F_{\eta(t)} \right] & =\frac{-2}{\sA^2}\ddt+\frac{2}{\sA^2}\frac{e^{\sA^2\ddt}-1}{\sA^2} \\
 &\quad - 2(t- \eta(t))\frac{e^{\sA^2\ddt} - 1}{\sA^2} +\ddt^2e^{\sA^2\ddt}      \\
  &=  \frac{2}{\sA^4}\left(e^{\sA^2\ddt}-1 - {\sA^2}\ddt e^{\sA^2\ddt} +\frac{\left(\sA^2\ddt \right)^2}{2}e^{\sA^2\ddt}    \right), 
\end{align*}
which is of order $\ddt^3$.

\end{proof}


\begin{proposition}\label{prop:a-priori-estimates-NS}Assume \ref{H1},\ref{H2}, \ref{H3}, \ref{H4} and there exists a suitable approximation of $P$ satisfying \eqref{eq:scheme-for-P}. Then, for $\dt$ such that \eqref{eq:condition_positivity_Dt} is satisfied, there exists a positive constant $\Cadd_\delta$, depending on $p$ and $\delta$, such that:
\begin{enumerate}
\item Existence of moments: For any $p\geq1$:
\begin{equation}\label{eq:moments}
\begin{aligned}
\sup_{0\leq t\leq T} \E\left[(\bXJ^\delta_t)^{2p} + (\bXA^\delta_t)^{2p}+(\XJ^\delta_t)^{2p}+(\XA^\delta_t)^{2p}\right] 	  &\leq \mathbf{C}_\delta,\\
\sup_{0\leq t\leq T} \E\left[(\bXJ_t)^{2p}+(\bXA_t)^{2p}\right] 	  &\leq \CP .
\end{aligned}
\end{equation}
\item Local Error: For all $p\geq 1/2$,
\begin{equation}\label{eq:local-errors}
\begin{aligned}
\sup_{0\leq t\leq T}\E\left[\left(\bXA^\delta_t-\bXA^\delta_{\eta(t)} \right)^{2p}\right] +\sup_{0\leq t\leq T}\E\left[\left(\bXJ^\delta_t-\bXJ^\delta_{\eta(t)} \right)^{2p}\right] &\leq \mathbf{C}_\delta \dt^{p\wedge 1},\nonumber\\
\sup_{0\leq t\leq T}\E\left[\left(\bE^\delta_t-\bE^\delta_{\eta(t)} \right)^{2p} \right]&\leq \mathbf{C}_\delta \dt^p.
\end{aligned}
\end{equation}

\item Jump contribution:
\begin{equation}\label{eq:jump_contribution}
\begin{aligned}
\sup_{0\leq t\leq T} \E\left[\left(\hXA^\delta_t-\bXA^\delta_t \right)^2\right] 	+\sup_{0\leq t\leq T} \E\left[\left(\hXJ^\delta_t-\bXJ^\delta_t \right)^2\right] 	  &\leq \Cadd_\delta\dt.
\end{aligned}
\end{equation}

\end{enumerate}

\end{proposition}

\begin{proof}

We start by proving the finiteness of the moments focusing on the process $\XA^\delta$. The proof for the process $\XJ$ will is analogous. 

We apply It\^o's formula to the process $\bXJ$ in \eqref{eq:explicit-XJ-tilde} in the time interval $[0,t]$:
\begin{equation*}
  \begin{aligned}
    &(\bXJ^\delta_t)^{2p} =\bXJ_0^{2p} + 2p\int_{0}^t (\bXJ^\delta_s)^{2p-1}\rhoJ(\bXJ_{\eta(s)}^\delta ,\bXA_{\eta(s)}^\delta,\bE_{\eta(s)}^\delta)\Psi_{1-1/\delta}(\bXA_{\eta(s)})\\
    &\qquad - 2p\int_{0}^t (\bXJ^\delta_s)^{2p-1}\kappaJ(\bE_{\eta(s)}^\delta)\Psi_{1-1/\delta}(\bXJ_{\eta(s)})ds \\
    &\qquad + 2p\int_{0}^t \sJ (\bXJ^\delta_s)^{2p-1}\hXJ^\delta_sd\WJ_s + p(2p-1)\int_{0}^t \sJ^2 (\bXJ^\delta_s)^{2p-2}(\hXJ^\delta_s)^2ds\\
    &\qquad + \int_{0}^{t+}\int_{-\infty}^\infty \left((\bXJ^\delta_{s^-} + \phiJ(\bXJ^\delta_{\eta(s)},z)\ind{\{s\leq \tau_{\eta(s)}^1\}})^{2p} - (\bXJ^\delta_{s^-})^{2p}\right)\Np(dz\otimes ds).
\end{aligned}  
\end{equation*}

Taking expectations --after a localization argument if needed-- and using that the numerical approximations are non-negative for $\dt$ satisfying \eqref{eq:condition_positivity_Dt}, we get:
\begin{equation*}
  \begin{aligned}
\E[(\bXJ^\delta_t)^{2p}] \leq&\E[\bXJ_0^{2p}]+\frac{2p}\delta\int_{0}^t \E[(\bXJ^\delta_s)^{2p-1}\rhoJ(\bXJ_{\eta(s)}^\delta ,\bXA_{\eta(s)}^\delta,\bE_{\eta(s)}^\delta)]ds\\
&+ p(2p-1)\sJ^2\int_{0}^t  \E[(\bXJ^\delta_s)^{2p-2}(\hXJ^\delta_s)^2]ds\\
    & + \E\left[\int_{0}^{t+}\int_{\R} \left((\bXJ^\delta_{s^-} + \phiJ(\bXJ^\delta_{\eta(s)},z)\ind{\{s\leq \tau_{\eta(s)}^1\}})^{2p} - (\bXJ^\delta_{s^-})^{2p}\right)\Np(dz\otimes ds)\right].
\end{aligned}  
\end{equation*}

Then, applying Young inequality together with assumptions \ref{H1}-\ref{H2}, there exists a non-negative constant $C_p$, depending on $p$, such that
\begin{equation*}
  \begin{aligned}
\E[(\bXJ^\delta_t)^{2p}] \leq&\E[\bXJ_0^{2p}] +\frac{C_p}\delta\int_{0}^t \E[(\bXJ^\delta_s)^{2p-1}]ds + C_p\int_{0}^t  \E[(\bXJ^\delta_s)^{2p}] + \E[(\hXJ^\delta_s)^{2p}]ds \\
&\qquad + C_p\int_{0}^{t}\int_{-\infty}^\infty\E\left[ \phiJ^{2p}(\bXJ^\delta_{\eta(s)},z)\ind{\{s\leq \tau_{\eta(s)}^1\}}\right]\lambda \nu(dz)ds\\
\leq&\E[\bXJ_0^{2p}] + C_\delta\int_{0}^t \E[1 + (\bXJ^\delta_s)^{2p}+(\hXJ^\delta_s)^{2p}]ds \\
&\qquad+ C_p\int_{0}^{t}\E\left[ (\bXJ^\delta_{\eta(s)})^{2p}\ind{\{s\leq \tau_{\eta(s)}^1\}}\right]ds.
\end{aligned}  
\end{equation*}

From independence between $\tau_{\eta(t)}^1$,  the time until the first atom of $\Np$ with ``time'' coordinate larger than $\eta(t)$ and $\F_{\eta(t)}$ we can condition the last integral and get
\begin{align*}
\E\left[ (\bXJ^\delta_{\eta(s)})^{2p}\ind{\{s\leq \tau_{\eta(s)}^1\}}\right] &=\E\left[ (\bXJ^\delta_{\eta(s)})^{2p}\right]\P\left(s-\eta(s)\leq \tau_{\eta(s)}^1-\eta(s)\right)\\
&=\E\left[ (\bXJ^\delta_{\eta(s)})^{2p}\right]e^{-\lambda(s-\eta(s))}.
\end{align*}

Hence, for some constant $C_\delta$ depending on $\delta$:
\begin{equation*}
  \begin{aligned}
\E[(\bXJ^\delta_t)^{2p}] 
\leq&C_\delta\int_{0}^t \E[1 + (\bXJ^\delta_s)^{2p}+(\hXJ^\delta_s)^{2p} + (\bXJ^\delta_{\eta(s)})^{2p}]ds.
\end{aligned}  
\end{equation*}

The same can be done with the processes $\hXJ, \bXA$ and $\hXA$, deducing from Gronwall inequality that, for some constant $C_\delta>0$ we have
\begin{equation*}
  \begin{aligned}
\sup_{0\leq t\leq T}\E[(\bXJ^\delta_t)^{2p} + (\hXJ^\delta_t)^{2p} + (\bXA^\delta_t)^{2p} + (\hXA^\delta_t)^{2p}] 
\leq&C_\delta.
\end{aligned}  
\end{equation*}

Repeating the procedure above we estimate
\begin{equation*}
  \begin{aligned}
\E[(\XJ^\delta_t)^{2p}] \leq&\E[\XJ_0^{2p}]+\frac{C_p}\delta\int_{0}^t \E[(\XJ^\delta_s)^{2p-1}]ds + p(2p-1)\sJ^2\int_{0}^t  \E[(\XJ^\delta_s)^{2p}]ds\\
    & + \E\left[\int_{0}^{t+}\int_{-\infty}^\infty \left((\XJ^\delta_{s^-} + \phiJ(\XJ^\delta_{s},z))^{2p} - (\XJ^\delta_{s^-})^{2p}\right)\Np(dz\otimes ds)\right]\\
    &\leq C_\delta +C_\delta \int_{0}^t  \E[(\XJ^\delta_s)^{2p}]ds,
\end{aligned}  
\end{equation*}
and in the same way for $\XA^\delta$, concluding
\begin{equation*}
  \begin{aligned}
\sup_{0\leq t\leq T}\E[(\XJ^\delta_t)^{2p} + (\XA^\delta_t)^{2p}] 
\leq&C_\delta.
\end{aligned}  
\end{equation*}

In the case of $\bXA$ and $\bXJ$ we repeat the same arguments with It\^o formula, Young inequality, \ref{H1}-\ref{H2}, proving that there exists a non-negative constant $C_p$ depending on $p$ such that
\begin{equation*}
  \begin{aligned}
\sup_{0\leq t\leq T}\E[(\XJ_t)^{2p} + (\XA_t)^{2p}] 
\leq&\mathbf{C}_p.
\end{aligned}  
\end{equation*}

We now proceed with the estimation of the local error for the process $\bXJ$:
\begin{equation*}
  \begin{aligned}
    &(\bXJ^\delta_t - \bXJ^\delta_{\eta(t)})^{2p} =\\
    &\qquad  \Big(\int_{\eta(t)}^t \left(\rhoJ(\bXJ_{\eta(s)}^\delta ,\bXA_{\eta(s)}^\delta,\bE_{\eta(s)}^\delta)\Psi_{1-1/\delta}(\bXA_{\eta(s)}) -\kappaJ(\bE_{\eta(s)}^\delta)\Psi_{1-1/\delta}(\bXJ_{\eta(s)})\right)ds \\
    &\qquad + \int_{\eta(t)}^t \sJ \hXJ^\delta_sd\WJ_s + \int_{\eta(t)}^{t+}\int_{-\infty}^\infty  \phiJ(\bXJ^\delta_{\eta(s)},z)\ind{\{s^{-}\leq \tau_{\eta(s)}^1\}}\Np(dz\otimes ds)\Big)^{2p},
    \end{aligned}  
\end{equation*}

\begin{equation*}
\begin{aligned}
   \E\left[(\bXJ^\delta_t - \bXJ^\delta_{\eta(t)})^{2p}\right] &\leq \frac{C_p}{\delta^{2p}} \ddt^{2p} + C_p\E\Big(\int_{\eta(t)}^t \hXJ^\delta_sd\WJ_s\Big)^{2p} \\
   &\qquad +C_p\E\Big( \int_{\eta(t)}^{t+}\int_{-\infty}^\infty  \phiJ(\bXJ^\delta_{\eta(s)},z)\ind{\{s^{-}\leq \tau_{\eta(s)}^1\}}\Np(dz\otimes ds)\Big)^{2p}.
\end{aligned}  
\end{equation*}

Using Burkholder-Davis-Gundy and Jensen inequalities, and the moments of the numerical scheme,
\begin{equation*}
  \begin{aligned}
   \E\left[(\bXJ^\delta_t - \bXJ^\delta_{\eta(t)})^{2p}\right] &\leq \frac{C_p}{\delta^{2p}} \ddt^{2p} + C_{\text{\tiny{BDG}}}\E\left[\Big(\int_{\eta(t)}^t (\hXJ^\delta_s)^2ds\Big)^{p}\right] \\
   &\qquad +C_p\E\Big( \int_{\eta(t)}^{t+}\int_{-\infty}^\infty  \phiJ(\bXJ^\delta_{\eta(s)},z)\ind{\{s^{-}\leq \tau_{\eta(s)}^1\}}\Np(dz\otimes ds)\Big)^{2p}\\
   &\leq \frac{C_p}{\delta^{2p}} \ddt^{2p} + C_{\text{\tiny{BDG}}}\sup_{0\leq t\leq T}\E\left[ (\hXJ^\delta_t)^{2p}\right]\, \ddt^{p}   +C_p\E\Big[\ind{\{\eta (t) \le \tau^{1}_{\eta(t)} < t\}} \phiJ^{2p} \left(\bXJ_{\eta(t)}^\delta, \xi_1\right)\Big].
\end{aligned}  
\end{equation*}

From assumption \ref{hip:linear-growth-phi}, the finiteness of the moments of $\bXJ$, and independence between $\tau_{\eta(t)}^1$,  the time until the first atom of $\Np$ with ``time'' coordinate larger than $\eta(t)$ and $\xi_1$ and $\F_{\eta(t)}$:
\begin{align*}
\E\left[\left(  \ind{\{\eta (t) \le \tau^{1}_{\eta(t)} < t\}} \phiJ \left(\bXJ_{\eta(t)}^\delta, \xi_1\right)\right)^2 \right] 
& =   \P\left(\eta (t) \le \tau^{1}_{\eta(t)} < t \right)\E\left[ \phiJ^2 \left(\bXJ_{\eta(t)}^\delta, \xi_1\right)\right] \\
&\leq C_\delta \left(1-e^{-\lambda\dt}\right) \leq C_\delta \dt.
\end{align*}

This also gives us
\[\E\left[\left(\hXJ_t^\delta-\bXJ_t^\delta \right)^2\right]\leq  C_\delta \dt.\]

Summarizing the two last estimations:
\begin{align*}
 \sup_{0\leq t\leq T}\E\left[\left(\bXJ_t^\delta-\bXJ_{\eta(t)}^\delta \right)^{2p}\right]\leq C_\delta \dt^{p\wedge 1},\text{ and }\qquad \sup_{0\leq t\leq T}\E\left[\left(\hXJ_t^\delta-\bXJ_t^\delta \right)^2\right] 	  	& \leq C_\delta \dt,
\end{align*}
which will hold also for $\bXA$, concluding the proof of assertion \eqref{eq:jump_contribution} and \eqref{eq:local-errors}.

We conclude the proof by analyzing the local error of the compliance process. For this we come back to the definition of the numerical scheme $\bE^\delta$:
\begin{align*}
  (\bE^\delta_t - \bE^\delta_{\eta(t)})^{2p} &= \left(\int_{\eta(t)}^t\left(\beta_1(\bZZ^{\delta}_{\eta(t)})(1-\Psi_{\delta}(\bE^\delta_{\eta(t)}))  -\beta_0(\bZZ^{\delta}_{\eta(t)})\Psi_{\delta}(\bE^\delta_{\eta(t)}))\right)ds\right.\\
  &\qquad +\left.\int_{\eta(t)}^t\sE\sqrt{\Psi_{\delta}(\bE^\delta_{\eta(t)})(1-\Psi_{\delta}(\bE^\delta_{\eta(t)}))}dW^{(E)}_s\right)^{2p}.
\end{align*}

Hence, it is straightforward that
\begin{align*}
(\bE^\delta_t - \bE^\delta_{\eta(t)})^{2p} &\leq C_\delta\ddt^{2p} + C_p\left(\int_{\eta(t)}^t\sqrt{\Psi_{\delta}(\bE^\delta_{\eta(t)})(1-\Psi_{\delta}(\bE^\delta_{\eta(t)}))}dW^{(E)}_s\right)^{2p},
\end{align*}
for some constant $C_\delta>0$ depending on $\delta$.

Taking expectations and using the moments of a Gaussian random variable, we have for all $t\in[0,T]$:
\begin{align}\label{eq:localerror_E_estim}
\E[(\bE^\delta_t - \bE^\delta_{\eta(t)})^{2p}] &\leq C_\delta\ddt^{2p} + C_\delta \E\left[(W^{(E)}_t- W^{(E)}_{\eta(t)})^{2p}\right]\nonumber\\
&\leq C_\delta\ddt^{p},
\end{align}
and in particular 
\[\E[(\bE^\delta_t - \bE^\delta_{\eta(t)})^{2p}]\leq C_\delta\dt^{p},\]
concluding the proof.
\end{proof}


\begin{proof} [Proof of Proposition \ref{lem:control-of-L2-norm-hats-and-bars} ]
Notice that, following the notations of Proposition \ref{prop:control-of-brownian-exponentials}:
\begin{align*}
\bXA_t - \scXA_t&= \left(\bXA_{\eta(t)} - \scXA_{\eta(t)} \right)\eA(\eta(t),t) + \left(\rhoA\bXJ_{\eta(t)}-\kappaA(\bE_{\eta(t)})\bXA_{\eta(t)}  \right)I(\eta(t),t) \\
&\quad+    \left(\rhoA\left(\bXJ_{\eta(t)} -\scXJ_{\eta(t)} \right)-\kappaA(\bE_{\eta(t)})\left(\bXA_{\eta(t)} -\scXA_{\eta(t)}  \right) \right)\ddt\eA(\eta(t),t) \\
&\quad
+   \ind{\{\eta(t) \le \tau^{1}_{\eta(t)} < t\}} \left(\phiA \left(\bXA_{\eta(t)}, \xi_1\right)-\phiA \left(\scXA_{\eta(t)}, \xi_1\right) \right).
\end{align*}

Then,
\begin{align*}
&( \bXA_{t} - \scXA_{t} )^2\\
& \leq
\left(\bXA_{\eta(t)} - \scXA_{\eta(t)} \right)^2 \eA^2(\eta(t),t)) \\
&\quad+     
3\left(\rhoA\bXJ_{\eta(t)}-\kappaA(\bE_{\eta(t)})\bXA_{\eta(t)}  \right)^2I^2(\eta(t),t) \\
&\quad+    3\left(\rhoA\left(\bXJ_{\eta(t)} -\scXJ_{\eta(t)} \right)-\kappaA(\bE_{\eta(t)})\left(\bXA_{\eta(t)} -\scXA_{\eta(t)}  \right) \right)^2\ddt^2\eA^2(\eta(t),t) \\
&\quad
+   3\cdot\ind{\{\eta(t) \le \tau^{1}_{\eta(t)} < t\}} \left(\phiA \left(\bXA_{\eta(t)}, \xi_1\right)-\phiA \left(\scXA_{\eta(t)}, \xi_1\right) \right)^2 \\
&\quad+
2\left(\bXA_{\eta(t)} - \scXA_{\eta(t)} \right)\eA(\eta(t),t) \\
&\qquad\times     
\Big\{\left(\rhoA\bXJ_{\eta(t)}-\kappaA(\bE_{\eta(t)})\bXA_{\eta(t)}  \right)I(\eta(t),t) \\
&\qquad\qquad+ \ind{\{\eta(t) \le \tau^{1}_{\eta(t)} < t\}} \left(\phiA \left(\bXA_{\eta(t)}, \xi_1\right)-\phiA \left(\scXA_{\eta(t)}, \xi_1\right) \right)\\
&\qquad\qquad+    \left(\rhoA\left(\bXJ_{\eta(t)} -\scXJ_{\eta(t)} \right)-\kappaA(\bE_{\eta(t)})\left(\bXA_{\eta(t)} -\scXA_{\eta(t)}  \right) \right)\ddt\eA(\eta(t),t)  \Big\}.
\end{align*}

Therefore,
\begin{equation}\label{eq:Initial-bound-hat-bar}
\E\left[( \bXA_{t} - \scXA_{t} )^2  \right]  \leq  E_1 + E_2+E_3+E_4+E_5.
\end{equation}

In what follows we will treat each term in the right-hand side of \eqref{eq:Initial-bound-hat-bar} separately. We start with $E_1$ with the use of the tower property of conditional expectations, and applying Proposition  \ref{prop:control-of-brownian-exponentials}:
\begin{align*}
E_1
  &  = \E\left[\left(\bXA_{\eta(t)}-\scXA_{\eta(t)}\right)^2 \E\left[ \exp\left(-\sA^2\ddt + 2\sA\deltaWA{\eta(t)}{t}\right) / \F_{\eta(t)} \right]\right]    \\
  & =      \E\left[\left(\bXA_{\eta(t)}-\scXA_{\eta(t)}\right)^2\right]e^{\sA^2\ddt}.
\end{align*}

Now we estimate $E_2$, by applying Proposition \ref{prop:control-of-brownian-exponentials} and Proposition \ref{prop:a-priori-estimates-NS}:
\begin{align*}
E_2&=3\E\left[\left(\rhoA\bXJ_{\eta(t)}-\kappaA(\bE_{\eta(t)})\bXA_{\eta(t)}  \right)^2\E\left[I^2(\eta(t),t)\Big|\F_{\eta(t)} \right] \right]\\
& \leq \Cadd \E\left[\left(\rhoA\bXJ_{\eta(t)}-\kappaA(\bE_{\eta(t)})\bXA_{\eta(t)}  \right)^2\ddt^3\right]\\ 
  &  \leq  \Cadd \ddt^3.    
\end{align*}

Analogously, we estimate $E_3$:
\begin{align*}
E_3
& = 3 \ddt^2 e^{\sA^2\ddt} \E\left[   \left(\rhoA(\bXJ_{\eta(t)}-\scXJ_{\eta(t)})-\kappaA(\bE_{\eta(t)})(\bXA_{\eta(t)}-\scXA_{\eta(t)})  \right)^2\right]\\
&\leq C \ddt^2 e^{\sA^2\ddt} \E\left[  (\bXJ_{\eta(t)}-\scXJ_{\eta(t)})^2+(\bXA_{\eta(t)}-\scXA_{\eta(t)})^2\right].
\end{align*}

For the term $E_4$ we use independence and \ref{H2} as follows:
\begin{align*}
E_4 & = 3\E\Big[\Big( \ind{\{\eta(t) \le \tau^{1}_{\eta(t)} < t\}} \left(\phiA \left(\bXA_{\eta(t)}, \xi_1\right)-\phiA \left(\scXA_{\eta(t)}, \xi_1\right) \right)
  \Big)^2 \Big]       \\
  &   \leq C\E\left[\left( \bXA_{\eta(t)}-\scXA_{\eta(t)}\right)^2 \right]\P\left(\eta(t) \le \tau^{1}_{\eta(t)} < t \right)   \\
  &    \leq C\E\left[\left( \bXA_{\eta(t)}-\scXA_{\eta(t)}\right)^2 \right]\left(1-e^{-\lambda\ddt}\right)\\
  &\leq C\E\left[\left( \bXA_{\eta(t)}-\scXA_{\eta(t)}\right)^2 \right]\ddt.
\end{align*}

We continue with the estimation of $E_5$. Notice that, keeping the non-negative terms:
\begin{equation}\label{eq:e2-splitting}
\begin{aligned}
&E_5\\
&\le
2\E\Big[\left(\bXA_{\eta(t)} - \scXA_{\eta(t)} \right)\left(\rhoA\bXJ_{\eta(t)}-\kappaA(\bE_{\eta(t)})\bXA_{\eta(t)}  \right)I(\eta(t),t)\eA(\eta(t),t)  \Big]\\
&\;
+2\rhoA\E\Big[\left(\bXA_{\eta(t)} - \scXA_{\eta(t)} \right)\left(\bXJ_{\eta(t)} -\scXJ_{\eta(t)} \right)\ddt\eA^2(\eta(t),t) \Big]\\
&\;
+2\E\Big[\left(\bXA_{\eta(t)} - \scXA_{\eta(t)} \right)\eA(\eta(t),t)\times     
\Big\{ \ind{\{\tau^{1}_{\eta(t)} < t\}} \left(\phiA \left(\bXA_{\eta(t)}, \xi_1\right)-\phiA \left(\scXA_{\eta(t)}, \xi_1\right) \right) \Big\}\Big]\\
&= 2\E\Big[\left(\bXA_{\eta(t)} - \scXA_{\eta(t)} \right)\left(\rhoA\bXJ_{\eta(t)}-\kappaA(\bE_{\eta(t)})\bXA_{\eta(t)}  \right)\E\left[I(\eta(t),t)\eA(\eta(t),t) \big| \F_{\eta(t)} \right]  \Big]\\
&\; +  2\rhoA\E\Big[\left(\bXA_{\eta(t)} - \scXA_{\eta(t)} \right)\left(\bXJ_{\eta(t)} -\scXJ_{\eta(t)} \right)\ddt\E\Big[\eA^2(\eta(t),t)\Big|\F_{\eta(t)}\Big] \Big]    \\
&\;
+2\E\Big[\left(\bXA_{\eta(t)} - \scXA_{\eta(t)} \right)\eA(\eta(t),t)\times     
\Big\{ \ind{\{\tau^{1}_{\eta(t)} < t\}} \left(\phiA \left(\bXA_{\eta(t)}, \xi_1\right)-\phiA \left(\scXA_{\eta(t)}, \xi_1\right) \right) \Big\}\Big]\\
&= E_5^1+E_5^2+E_5^3.
\end{aligned}
\end{equation}

Then, from Proposition \ref{prop:control-of-brownian-exponentials}, assumption \ref{H1} and Young's inequality we have 
\begin{align*}
E_5^1 & =  2\E\left[\left(\bXA_{\eta(t)} - \scXA_{\eta(t)} \right)\left(\rhoA\bXJ_{\eta(t)}-\kappaA(\bE_{\eta(t)})\bXA_{\eta(t)}  \right)\left( \frac{e^{\sA^2\ddt}-1}{\sA^2}-\ddt e^{\sA^2\ddt} \right)  \right]     \\
  & \leq C\E\Big[\left|\bXA_{\eta(t)} - \scXA_{\eta(t)} \right|\left|\rhoA\bXJ_{\eta(t)}-\kappaA(\bE_{\eta(t)})\bXA_{\eta(t)}  \right|\ddt^2  \Big] 
      \\
   & \leq C\ddt\E\Big[\left(\bXA_{\eta(t)} - \scXA_{\eta(t)} \right)^2 \Big] + C \ddt^3\E\left[(\bXJ_{\eta(t)})^2+(\bXA_{\eta(t)})^2\right].
\end{align*}

Similarly,
\begin{align*}
E_5^2
&=2\rhoA\E\Big[\left(\bXA_{\eta(t)} - \scXA_{\eta(t)} \right)\left(\bXJ_{\eta(t)} -\scXJ_{\eta(t)} \right)\ddt e^{\sA^2\ddt}\Big]    \\
  &   \leq 
    C\E\Big[\left(\bXA_{\eta(t)} - \scXA_{\eta(t)} \right)^2+\left(\bXJ_{\eta(t)} -\scXJ_{\eta(t)} \right)^2\Big]\ddt e^{\sA^2\ddt}.      
\end{align*}

For the third term in the right-hand side of \eqref{eq:e2-splitting}, we first use the conditional independence of $\eA(\eta(t),t)$, $\tau^{1}_{\eta(t)}$ and $\xi_1$ and then the Lipschitz property of $\phiA$ in \ref{H2} to obtain,
\begin{align*}
 &E_5^3\\
 & \leq 2\E\Big[\left|\bXA_{\eta(t)} - \scXA_{\eta(t)} \right|\E\left[\left|\phiA \left(\bXA_{\eta(t)}, \xi_1\right)-\phiA \left(\scXA_{\eta(t)}, \xi_1\right) \right|\eA(\eta(t),t) \ind{\{ \tau^{1}_{\eta(t)} < t\}}\Big| \F_{\eta(t)}\right]\Big]      \\
  &   \leq  C \E\Big[\left(\bXA_{\eta(t)} - \scXA_{\eta(t)} \right)^2 \Big]\P\left(\eta(t) \le \tau^{1}_{\eta(t)} < t \right)        \\
  &   =     C \E\Big[\left(\bXA_{\eta(t)} - \scXA_{\eta(t)} \right)^2 \Big](1-e^{\lambda\ddt})\\
    &   \leq     C\ddt \E\Big[\left(\bXA_{\eta(t)} - \scXA_{\eta(t)} \right)^2 \Big].
\end{align*}

With all the previous computations, we have
$$
E_5\leq  C\E\Big[\left(\bXA_{\eta(t)} - \scXA_{\eta(t)} \right)^2+\left(\bXJ_{\eta(t)} -\scXJ_{\eta(t)} \right)^2\Big]\ddt e^{\sA^2\ddt} +  \Cadd\ddt^3.
$$

Summarizing,
\begin{equation*}
\begin{aligned}
\E\left[( \bXA_{t} - \scXA_{t} )^2  \right]  &\leq  \E\left[\left(\bXA_{\eta(t)}-\scXA_{\eta(t)}\right)^2\right]e^{\sA^2\ddt}\\
&\quad
+  C\ddt e^{\sA^2\ddt}\E\left[\left( \bXA_{\eta(t)}-\scXA_{\eta(t)}\right)^2+\left( \bXJ_{\eta(t)}-\scXJ_{\eta(t)}\right)^2 \right] +  \Cadd\ddt^3.
\end{aligned}
\end{equation*}

In the same way, for $\XJ$ we show that
\begin{equation*}
\begin{aligned}
\E\left[( \bXJ_{t} - \scXJ_{t} )^2  \right]  &\leq  \E\left[\left(\bXJ_{\eta(t)}-\scXJ_{\eta(t)}\right)^2\right]e^{\sJ^2\ddt}\\
&\quad
+  C\ddt e^{\sJ^2\ddt}\E\left[\left( \bXA_{\eta(t)}-\scXA_{\eta(t)}\right)^2+\left( \bXJ_{\eta(t)}-\scXJ_{\eta(t)}\right)^2 \right] +  \Cadd\ddt^3.
\end{aligned}
\end{equation*}

Hence
\begin{equation}\label{eq:iteration-of-schemes}
\begin{aligned}
&\E\left[( \bXA_{t} - \scXA_{t} )^2+( \bXJ_{t} - \scXJ_{t} )^2  \right] \\
& \leq  e^{\sigma^2\ddt} (1+C\ddt)\E\left[\left( \bXA_{\eta(t)}-\scXA_{\eta(t)}\right)^2+\left( \bXJ_{\eta(t)}-\scXJ_{\eta(t)}\right)^2 \right]  +  \Cadd\ddt^3,
\end{aligned}
\end{equation}
where $\sigma^2 = \sA^2\vee \sJ^2$.

If we denote by $u_k:=\E\left[( \bXA_{t_k} - \scXA_{t_k} )^2+( \bXJ_{t_k} - \scXJ_{t_k} )^2  \right] $, we obtain from the previous inequality,
$$
u_{k+1}\leq e^{\sigma^2\dt} (1+C\dt)u_k+  \Cadd\dt^3,
$$
which iterating leads to
$$
u_{k+1}\leq e^{\sigma^2\dt(k+1)} (1+\Cmult\dt)^{k+1}u_0+  \Cadd\dt^3\frac{e^{\sigma^2\dt(k+1)} (1+\Cmult\dt)^{k+1}-1}{e^{\sigma^2\dt} (1+\Cmult\dt)-1}.
$$

Since $u_0=0$, for any $k=1,\ldots,N-1$, we obtain:
$$
u_{k+1}\leq \Cadd \dt^2.
$$
Plugging this estimate in \eqref{eq:iteration-of-schemes}, we obtain
$$
\E\left[( \bXA_{t} - \scXA_{t} )^2+( \bXJ_{t} - \scXJ_{t} )^2  \right]  \leq\Cadd \dt^2,
$$
as desired.
\end{proof}


\begin{proof}[Proof of Proposition \ref{prop:convergence-delta-process}]
To simplify the notation bellow, we define the approximation error
\begin{equation*}
\err^\delta_t:=\E\left[(\XJ^\delta_t - \bXJ^\delta_t)^2+(\XA^\delta_t - \bXA^\delta_t)^2+(E^\delta_t - \bE^\delta_t)^2\right],
\end{equation*}
for all $0\leq t\leq T.$

Denoting $\ZZ^{\delta}=(\XJ^\delta,\XA^\delta,E^\delta)$, applying It\^o's formula in $[\eta(t),t]$ and taking expectation we obtain:
\begin{align*}
&\err^\delta_t  =  \err^\delta_{\eta(t)}   +\E \int_{\eta(t)}^{t} 2\DXJ_s 
  \left\{ \rrJ(\ZZ^{\delta}_s)\Psi_{1-1/\delta}(\XA_s^\delta)- \rrJ(\bXJ^\delta_{\eta(t)},\bXA^\delta_{\eta(t)},\bE^\delta_{\eta(t)})\Psi_{1-1/\delta}(\bXA_{\eta(t)}^\delta) \right.\\
&\hspace{2in}\; 
  \left. - \kappaJ(E^\delta_s)\Psi_{1-1/\delta}(\XJ_s^\delta)  + \kappaJ(\bE^\delta_{\eta(t)})\Psi_{1-1/\delta}(\bXJ_{\eta(t)}^\delta)\right\}ds  \\
  & \; +\E \int_{\eta(t)}^{t} 2\DXA_s\left\{\PjA(\XJ^\delta_s- \bXJ^\delta_{\eta(s)})-\kappaA(E^\delta_s)\Psi_{1-1/\delta}(\XA^\delta_s) +\kappaA(\bE^\delta_{\eta(s)})\Psi_{1-1/\delta}(\bXA^\delta_{\eta(s)}) \right\}ds     \\
    &\;
    + \E\left\{\int_{\eta(t)}^{t}\sA^2\left(\XA^\delta_s- \hXA^\delta_s\right)^2ds  + \int_{\eta(t)}^{t}\sJ^2\left(\XJ^\delta_s- \hXJ^\delta_s\right)^2ds\right\}\\
    &\; 
    +  \E \int_{\eta(t)}^{t+}\int_{\R}\left\{\left( \DXA_{s^-} + \phiA(\XA^\delta_{s^-},z) - \phiA(\bXA^\delta_{\eta(t)},z)\ind{\{s\leq \tau_{\eta(t)}^1\}} \right)^2\right.\\
    &\hspace{3in}\left.- \left( \DXA_{s^-}\right)^2\right\} \Np(dz\otimes ds)\\
    &\; 
     +   \E\int_{\eta(t)}^{t+}\int_{\R}\left\{\left( \DXJ_{s^-} + \phiJ(\XJ^\delta_{s-},z) - \phiJ(\bXJ^\delta_{\eta(t)},z)\ind{\{s\leq \tau_{\eta(t)}^1\}} \right)^2\right.\\
     &\hspace{3in}\left.- \left( \DXJ_{s^-} \right)^2\right\}\Np(dz\otimes ds)\\
     & \;
    +\E \int_{\eta(t)}^t2\DE_s \Big\{\Delta\tilde\beta(\bXX_{\eta(t)}^{*\delta})\Psi_\delta(\bE_{\eta(t)}^\delta)(1-\Psi_\delta(\bE_{\eta(t)}^\delta))     \\
    & \qquad-\Delta\tilde\beta(\XX_{s}^{*\delta})\Psi_\delta(E_s^\delta)(1-\Psi_\delta(E_{s}^\delta))
    -(\bar\beta_1+\bar\beta_2)\Psi_\delta(\bE_{\eta(t)}^\delta)+(\bar\beta_1+\bar\beta_2)\Psi_\delta(E_{s}^\delta)  \Big\}ds    \\
  &   \; + \E \int_{\eta(t)}^t\sE^2\left[\sqrt{\Psi_\delta(\bE_{\eta(t)}^\delta)(1-\Psi_\delta(\bE_{\eta(t)}^\delta))} -\sqrt{\Psi_\delta(E_{s}^\delta)(1-\Psi_\delta(E_{s}^\delta))}  \right]^2ds\\
  &  =   \err^\delta_{\eta(t)}+ R^\delta_1+R^\delta_2+R^\delta_3+R^\delta_4+R^\delta_5+R^\delta_6+R^\delta_7.
\end{align*}

We proceed term by term using the Lipschitz and boundedness properties of  $\Psi_{1-1/\delta}$ and rate coefficients in \ref{H1}. For instance, for $ R^\delta_1$, by adding some pivots and using the boundedness and Lipschitz properties of $\kappaJ$, $\rJ$, $\Psi$ and $\Psi_{1-1/\delta}$ (with Lipschitz constant $1$), 
\begin{align*}
 &R^\delta_1=\\
 &\, 
\E \int_{\eta(t)}^{t} 2\DXJ_s 
  \Big\{  \left(\rJ(E^\delta_s)-\rJ(\bE^\delta_{\eta(t)})\right)\Psi\left(1-\frac{\XA^\delta_s+\XJ^\delta_s}\KA\right)\Psi_{1-1/\delta}(\XA_s^\delta)\\
  &\; - \rJ(\bE^\delta_{\eta(t)})\Psi\left(1-\frac{\bXA^\delta_{\eta(t)}+\bXJ^\delta_{\eta(t)}}\KA\right)\left(\Psi_{1-1/\delta}(\bXA_{\eta(t)}^\delta)-\Psi_{1-1/\delta}(A_{s}^\delta)\right)   \\
  &\; - \rJ(\bE^\delta_{\eta(t)})\left(\Psi\left(1-\frac{\bXA^\delta_{\eta(t)}+\bXJ^\delta_{\eta(t)}}\KA\right)-\Psi\left(1-\frac{\XA^\delta_s+\XJ^\delta_s}\KA\right)\right)\Psi_{1-1/\delta}(\XA_s^\delta)   \\
&\; 
 + \kappaJ(\bE^\delta_{\eta(t)})\left(\Psi_{1-1/\delta}(\bXJ_{\eta(t)}^\delta)-\Psi_{1-1/\delta}(\XJ_s^\delta)\right)+\Psi_{1-1/\delta}(\XJ_s^\delta)\left(\kappaJ(\bE^\delta_{\eta(t)})- \kappaJ(E^\delta_s)\right)\Big\}ds  \end{align*}
 
 \begin{align*}
R^\delta_1&\lesssim
\E \int_{\eta(t)}^{t} |\DXJ_s| 
  \Big\{  \frac1\delta |E^\delta_s-\bE^\delta_{\eta(t)}|+|\bXJ_{\eta(t)}^\delta-\XJ_s^\delta|+ |\bXA_{\eta(t)}^\delta-\XA_s^\delta|   \\
  &\qquad +\frac1\delta\left|\frac{\bXA^\delta_{\eta(t)}+\bXJ^\delta_{\eta(t)}}\KA-\frac{\XA^\delta_s+\XJ^\delta_s}\KA\right| \Big\}ds,
\end{align*}
for some non-negative constant $C$.

Then, applying Young's inequality, and the control for local errors in Proposition \ref{prop:a-priori-estimates-NS}, we have
\begin{align*}
 R^\delta_1&\lesssim 
\E \int_{\eta(t)}^{t}\frac{ |\DXJ_s| }\delta
  \Big( |\bE^\delta_{\eta(s)}-\bE^\delta_s| + |\DE_s|+ |\bXJ^\delta_{\eta(s)}-\bXJ^\delta_s| +|\DXJ_s|+|\bXA^\delta_{\eta(s)}-\bXA^\delta_s| +|\DXA_s|\Big)ds 
\\
& \leq
C_\delta \int_{\eta(t)}^{t}\err^\delta_s ds+C_\delta \dt^2.
\end{align*}

An analogous computation for $R^\delta_2$ shows that
\begin{align*}
 R^\delta_2
  &\leq
C_\delta  \int_{\eta(t)}^{t}\err^\delta_s ds+C_\delta \dt^2.
\end{align*}

For $R^\delta_3$, we just have to introduce the control for the jump contribution of the scheme and apply Proposition \ref{prop:a-priori-estimates-NS}. Then, we obtain
$$
R^\delta_3 \leq C \int_{\eta(t)}^{t}\err^\delta_s ds+ C_\delta \dt^2.
$$
The terms $R^\delta_4$ and $R^\delta_5$, the integrals with respect to the Poisson measure are completely analogous. For $R^\delta_4$ we have:
\begin{align*}
R^\delta_4&= \E\int_{\eta(t)}^{t+}\int_{\R}\left(  \phiA(\XA^\delta_{s-},z) - \phiA(\bXA_{\eta(t)}^\delta,z)\ind{\{s^{-}\leq \tau_{\eta(t)}^1\}}\right)^2\Np(dz\otimes ds)\\
&\qquad +\E\int_{\eta(t)}^{t+}\int_{\R}
+2\DXA_s\left(\phiA(\XA^\delta_{s-},z) - \phiA(\bXA_{\eta(t)}^\delta,z)\ind{\{s^{-}\leq \tau_{\eta(t)}^1\}}\right)\Np(dz\otimes ds)\\
 &\leq 2\E\int_{\eta(t)}^{t}\int_{\R}\left\{\left(  \phiA(\XA^\delta_{s},z) - \phiA(\bXA_{\eta(t)}^\delta,z)\ind{\{s^{-}\leq \tau_{\eta(t)}^1\}}\right)^2+\left(\DXA_s\right)^2\right\}\lambda\nu(dz)ds\\
 &\lesssim 
 \E\int_{\eta(t)}^{t}\int_{\R} 
 \left\{\left(  \phiA(\XA^\delta_{s},z)-\phiA(\bXA^\delta_{s},z)\right)^2 + \left(\phiA(\bXA^\delta_{s},z) - \phiA(\bXA_{\eta(t)}^\delta,z)\right)^2\right\}\nu(dz)ds\\
 &\qquad 
 +\E\int_{\eta(t)}^{t}\int_{\R} \Big(
\phiA(\bXA_{\eta(t)}^\delta,z)^2\ind{\{s> \tau_{\eta(t)}^1\}}
  + \left(\DXA_s\right)^2\Big)\nu(dz)ds.
  \end{align*}

Now we use assumptions \ref{hip:lipschitz-condition-phi}, \ref{hip:linear-growth-phi} and the estimation of local error in Proposition \ref{prop:a-priori-estimates-NS}, obtaining that
\begin{align*}
R^\delta_4&\leq
 C\E\int_{\eta(t)}^{t}
 \left\{\left(\bXA_{s}^\delta-\bXA_{\eta(t)}^\delta\right)^2
 +\left(\DXA_s \right)^2\right\}ds
+ C\E\int_{\eta(t)}^{t} (\bXA_{\eta(t)}^\delta)^2\ind{\{s> \tau_{\eta(t)}^1\}}ds\\
&\leq   C \int_{\eta(t)}^{t}\err^\delta_sds+ \Cadd_\delta \dt^2
+ C\E\left[(\bXA_{\eta(t)}^\delta)^2\int_{\eta(t)}^{t}\E\left[\ind{\{s> \tau_{\eta(t)}^1\}} / \F_{\eta(t)} \right] ds\right]\\
&\leq   C \int_{\eta(t)}^{t}\err^\delta_sds+ \Cadd_\delta \dt^2
+ C\E\left[(\bXA_{\eta(t)}^\delta)^2\right]\int_{\eta(t)}^{t}\E\left[\ind{\{s> \tau_{\eta(t)}^1\}} \right] ds,
\end{align*}
where in the last bound we have used the independence between $\tau_{\eta(t)}^1$,  the time until the first atom of $\Np$ with ``time'' coordinate larger than $\eta(t)$ and $\F_{\eta(t)}$. Notice that the random variable $(\tau_{\eta(t)}^1-\eta(t))$ follows an exponential distribution with rate $\lambda$, and thus:
\begin{align*}
\int_{\eta(t)}^{t}\E\left[\ind{\{s> \tau_{\eta(t)}^1\}}\right]     ds
& =    \int_{\eta(t)}^{t}\P\left(s-\eta(t)> \tau_{\eta(t)}^1-\eta(t)\right)  ds \\
  &     =    \int_{\eta(t)}^{t}\left(1 - e^{-\lambda(s-\eta(t))}\right)ds   \\
  &= t-\eta(t) + \frac{1}{\lambda}\left(e^{-\lambda\ddt} -1\right) \leq \frac{1}{\lambda}\frac{\left(\lambda( t-\eta(t)) \right)^2}{2}.
\end{align*}

From this and the moments of the scheme in Proposition \ref{prop:a-priori-estimates-NS} we obtain the estimation for $R^\delta_4$ as:
$$ R^\delta_4
\leq  C_\delta \int_{\eta(t)}^{t}\err^\delta_sds+ C_\delta \dt^2.
$$

Similarly, we prove that
$$ R^\delta_5
\leq  C_\delta\int_{\eta(t)}^{t}\err^\delta_sds+ C_\delta \dt^2.
$$

For $R^\delta_6$, we proceed as for $R^\delta_1$ and $R^\delta_2$, pivoting and using the boundedness and Lipschitz properties of $\Psi_\delta$ and $\Delta\tilde\beta$, to obtain from Young's inequality and local error in Proposition \ref{prop:a-priori-estimates-NS}:
\begin{align*}
R^\delta_6
  & =2\E \int_{\eta(t)}^t \Delta E^\delta_s \Big\{\left(\Delta\tilde\beta(\bXX_{\eta(t)}^{*\delta}) -\Delta\tilde\beta(\XX_{s}^{*\delta}) \right)\Psi_\delta(\bE_{\eta(t)}^\delta)(1-\Psi_\delta(\bE_{\eta(t)}^\delta))     \\
    & \quad
    +\Delta\tilde\beta(\XX_{s}^{*\delta})\left(\Psi_\delta(\bE_{\eta(t)}^\delta)(1-\Psi_\delta(\bE_{\eta(t)}^\delta)) -\Psi_\delta(E_s^\delta)(1-\Psi_\delta(E_{s}^\delta)) \right) \Big\}ds
      \\
  & \quad
  +2(\bar\beta_1+\bar\beta_2)\E \int_{\eta(t)}^t\Delta E^\delta_s \left(\Psi_\delta(E_{s}^\delta) -\Psi_\delta(\bE_{\eta(t)}^\delta) \right)ds     \\
  &      
  \leq
C \E \int_{\eta(t)}^t|\Delta E_s^\delta| \left\|\bXX_{\eta(t)}^{*\delta} -\XX_{s}^{*\delta} \right\|    ds 
  + C \E \int_{\eta(t)}^t|\Delta E^\delta_s |\left|E_{s}^\delta -\bE_{\eta(t)}^\delta \right| ds\\
    &\leq  C_\delta  \int_{\eta(t)}^t\err_s^{\delta}ds + C_\delta \dt^2.
\end{align*}

Finally, since the map $x\mapsto \sqrt{x(1-x)}$ is Lipschitz in $[\delta, 1-\delta]$, with constant of order $\frac1{\sqrt{\delta}}$, it follows :
$$
R^\delta_7 \leq C_\delta  \int_{\eta(t)}^t\E\left[(\bE^\delta_{\eta(t)} - E^\delta_s)^2\right]ds \leq C_\delta  \int_{\eta(t)}^t\err_s^{\delta}ds + C_\delta\dt^2.
$$

Summarizing the estimations above, we obtain
$$
\err^\delta_t \leq \err^\delta_{\eta(t)} +  C_\delta  \int_{\eta(t)}^t\err_s^{\delta}ds + C_\delta \dt^2.
$$

Therefore, by Gronwall's inequality
\begin{equation}\label{eq:iterative-bound-for-L2-error-delta-scheme}
\begin{aligned}
\err^\delta_t \leq  \left[  \err^\delta_{\eta(t)} +  C_\delta\dt^{2} \right]e^{C_\delta \dt}.
\end{aligned}
\end{equation}

Evaluating in $t=t_{k+1}$ we obtain, for any $k\in\{0,\ldots, N-1\}$:
$$\err^\delta_{t_{k+1}}  \leq  e^{C_\delta\dt}\err^\delta_{t_k} + C_\delta\dt^{2}e^{C_\delta\dt}, $$
which is a recurrence inequality of the form $a_{k+1}\leq Aa_k+B$ with estimation
$$a_{k+1}\leq B\frac{A^{k+1}-1}{A-1}\leq B\frac{A^{k+1}}{A-1}.$$

In our case,
$$\err^\delta_{t_{k+1}}  \leq C_\delta\dt^{2}\frac{\left(e^{C_\delta\dt}\right)^{k+2}}{e^{C_\delta\dt}-1}\leq   C_\delta \dt,\forall k\in\{0,\ldots, N-1\},$$
where the constant $C_\delta>0$ is changing from line to line. 

Introducing this last bound in \eqref{eq:iterative-bound-for-L2-error-delta-scheme}, we obtain for all $t\in[0,T]$,
$$
\err^\delta_t \leq  C_\delta \dt  +  C_\delta\dt^{2}  \leq  C_\delta\dt,
$$
obtaining the estimation \eqref{eq:approximation-error-for-delta-scheme-L2}, as desired.

\paragraph{The Burkholder-Davis-Gundy part:}
It remains to prove \eqref{eq:approximation-error-for-delta-scheme-L2}. From the definition of $\bXA^\delta_t$ and $\XA^\delta_t$ we compute, for any $t\in[0,T]$:
\begin{align*}
\bXA^\delta_t - \XA^\delta_t  &= \int_{0}^t\left(\rhoA\left(\bXJ^\delta_{\eta(s)}-\XJ^\delta_s\right)-\kappaA(\bE^\delta_{\eta(s)})\Psi_{1-1/\delta}(\bXA^\delta_{\eta(s)})+\kappaA(E^\delta_s)\Psi_{1-1/\delta}(\XA^\delta_s)\right)ds \\
&\qquad + \sA \int_{0}^t \left(\hXA^\delta_s - \XA^\delta_s\right)d\WA_s\\
&\qquad + \int_{0}^{t+}\int_{-\infty}^\infty \left(\phiA(\bXA^\delta_{\eta(s)},z)\ind{\{s^{-}\leq \tau_{\eta(s)}^1\}}-\phiA (\XA^\delta_s, z)\right)\Np(dz\otimes ds).
\end{align*}
Then, using  triangle inequality, we have for any $t\in[0,T]$:
\begin{align*}
&\E\left[\sup_{u\leq t}|\DXA_u| \right]\\
&\quad \leq  \E\left[\int_{0}^{t}| \rhoA\left(\XJ^\delta_s - \bXJ^\delta_{\eta(s)} \right)  - \left(\kappaA(E^\delta_s )\Psi_{1-1/\delta}(\XA^\delta_s)  -\kappaA(\bE^\delta_{\eta(s)})\Psi_{1-1/\delta}(\bXA^\delta_{\eta(s)})\right)|ds\right]\\
 &\qquad+ \sA\E \left[\sup_{u\leq t}\left|\int_{0}^{u}\left(\XA^\delta_s-\hXA^\delta_s\right)d\WA_s\right| \right] \\
 &\qquad+ \E\left[\int_{0}^{t+}\int_{\R} |\phiA(\XA^\delta_{s-},z) - \phiA(\bXA_{\eta(s)}^\delta,z)\ind{\{s^{-}\leq \tau_{\eta(s)}^1\}}|\Np(dz\otimes ds)\right]\\
 &\quad = E^A_1+E^A_2+E^A_3.
\end{align*}

Similar to the $L^2$ estimations above, we use assumptions \ref{H1}-\ref{H2} and local error from Proposition \ref{prop:a-priori-estimates-NS}:
\begin{align*}
E^A_1 &\leq C_\delta \int_{0}^{t}\E\left[ |\DXA_s|+|\DXJ_s| + |\DE_s| \right]ds + C_\delta \sqrt{\dt}\\  
&\leq \Cadd_\delta \sqrt{\dt},
\end{align*}
where the last inequality comes from the $L^2$ error in \eqref{eq:approximation-error-for-delta-scheme-L2}.

Analogously, we have
\begin{align*}
E^A_3    & =      \E\int_{0}^{t}\int_{-\infty}^\infty |\phiA(\XA^\delta_{s},z) -\phiA(\bXA^\delta_{s},z)|\lambda\nu(dz)ds\\
&\quad +\E\int_{0}^{t}\int_{-\infty}^\infty |\phiA(\bXA^\delta_{s},z) - \phiA(\bXA_{\eta(s)}^\delta,z)|\lambda\nu(dz)ds\\    
    &\quad
  +\E\int_{0}^{t}\int_{-\infty}^\infty |\phiA(\bXA^\delta_{\eta(s)},z)|\ind{\{s> \tau_{\eta(s)}^1\}}\nu(dz)ds\\    
  &\leq
  C_\delta\E\int_{0}^{t}\left(|\DXA_s|+|\bXA^\delta_{s} - \bXA_{\eta(s)}^\delta|\right)ds +\int_{0}^{t}\E\left[|\bXA^\delta_{\eta(s)}|\ind{\{s> \tau_{\eta(s)}^1\}} \right]ds\\    
  &\leq
  C_\delta\sqrt{\dt}
  +C_\delta\int_{0}^{t}\E\left[|\bXA^\delta_{\eta(s)}| \right]\P\left(s> \tau_{\eta(s)}^1 \right)ds\\    
      &\leq
  \Cadd_\delta\sqrt{\dt}
  +\Cmult_\delta\int_{0}^{t}1 - e^{-\lambda(s- {\eta(s)})}ds\\    
        &\leq
  C_\delta\sqrt{\dt}
  +C_\delta\int_{0}^{t}( s- {\eta(s)})ds \leq
  C_\delta\sqrt{\dt}.
\end{align*}

In the case of $E^A_2$, we use the Burkholder-Davis-Gundy and Jensen inequalities together with Proposition \ref{prop:a-priori-estimates-NS} and \eqref{eq:approximation-error-for-delta-scheme-L2}, from which there exists a non-negative constant $C_{\text{\tiny{BDG}}}$ to obtain
\begin{align*}
E^A_2 & \leq C_{\text{\tiny BDG}}\E\left[\sqrt{\int_{0}^{t}\left(\XA^\delta_s-\hXA^\delta_s\right)^2ds} \right]      \\
  &    \leq C_{\text{\tiny BDG}}\sqrt{\int_{0}^{t}\E\left[\left(\XA^\delta_s-\bXA^\delta_s\right)^2+\left(\bXA^\delta_s-\hXA^\delta_s\right)^2 \right]ds}  \leq C_\delta\sqrt{\dt}.
\end{align*}

Summarizing the estimations for $E_i^A$ for $i=1,2,3$, we conclude:
$$\E\left[\sup_{t\leq T}|\DXA_t| \right]\leq C_\delta \sqrt{\dt}.$$

A completely analogous computation shows
$$\E\left[\sup_{t\leq T}|\DXJ_t| \right]\leq C_\delta \sqrt{\dt}.$$

Regarding the approximation of the compliance process, we proceed similarly, using Lipschitz properties of the rates, BDG and Jensen inequalities: 
\begin{align*}
&\E\left[\sup_{0\leq s \leq t} \left|\DE_s\right|\right]\\
 &\; \lesssim
\int_{0}^t\E\left|\Delta \tilde{\beta}(\ZZ^{*\delta}_s)\Psi_{\delta}(E^\delta_s)(1-\Psi_{\delta}(E^\delta_s))  
-\Delta \tilde{\beta}(\bZZ^{*\delta}_{\eta(s)})\Psi_\delta(\bE^\delta_{\eta(s)})(1-\Psi_\delta(\bE^\delta_{\eta(s)}))  \right|ds \\
&\quad +\left(\bar\beta_1+\bar\beta_2\right)\int_{0}^t\E\left|
\Psi_{\delta}(E^\delta_s)-\Psi_\delta(\bE^\delta_{\eta(s)})
\right|ds \\
&\quad+\E\left[\sup_{0\leq u\leq t}\left|\int_{0}^u\left(\sqrt{\Psi_{\delta}(E^\delta_s)(1-\Psi_{\delta}(E^\delta_s))}
-\sqrt{\Psi_{\delta}(\bE^\delta_{\eta(s)})(1-\Psi_{\delta}(\bE^\delta_{\eta(s)}))}\right)dW^{(E)}_s\right|\right]\\
 &\;\leq
 C\int_{0}^t
\E\left[\left|E^\delta_s-\bE^\delta_{\eta(s)} \right|\right]
+\E\left[\left\|\ZZ^{*,\delta}_s-\bZZ^{*,\delta}_{\eta(s)} \right\|_1\right]ds \\
&\quad+C_{\text{\tiny{BDG}}}\sqrt{\int_{0}^t\E\left[\left(\sqrt{\Psi_{\delta}(E^\delta_s)(1-\Psi_{\delta}(E^\delta_s))}
-\sqrt{\Psi_{\delta}(\bE^\delta_{\eta(s)})(1-\Psi_{\delta}(\bE^\delta_{\eta(s)}))}\right)^2\right]ds}.
\end{align*}

We already know from Proposition \ref{prop:a-priori-estimates-NS} and \eqref{eq:approximation-error-for-delta-scheme-L2} that
$$
\E\left[\left|E^\delta_s-\bE^\delta_{\eta(s)} \right|\right]
+\E\left[\left\|\ZZ^{*,\delta}_s-\bZZ^{*,\delta}_{\eta(s)} \right\|_1\right]
\leq C_\delta \sqrt{\dt}.
$$

On the other hand, since the map $x\mapsto \sqrt{x(1-x)}$ is Lipschitz in $[\delta, 1-\delta ]$ we get
\begin{align*}
&\E\left[\left(\sqrt{\Psi_{\delta}(E^\delta_s)(1-\Psi_{\delta}(E^\delta_s))}
-\sqrt{\Psi_{\delta}(\bE^\delta_{\eta(s)})(1-\Psi_{\delta}(\bE^\delta_{\eta(s)}))}\right)^2\right]\\
&\qquad \leq  C_\delta\E\left[\left|E^\delta_s-\bE^\delta_{\eta(s)} \right|^2 \right]\leq C_\delta \dt,
\end{align*}
and putting all together,
\begin{align*}
\E\left[\sup_{0\leq s \leq T} \left|\DE_s\right|\right] &\leq C_\delta \sqrt{\dt},
\end{align*}
concluding the proof of \eqref{eq:approximation-error-for-delta-scheme}.
\end{proof}

\end{document}